\documentclass[final,3p,times]{elsarticle}

\usepackage{amssymb}

\pdfoutput=1

\journal{X}

\usepackage{lipsum}
\usepackage{amsfonts}
\usepackage{graphicx}
\usepackage{epstopdf}
\usepackage{algorithmic}
\ifpdf
  \DeclareGraphicsExtensions{.eps,.pdf,.png,.jpg}
\else
  \DeclareGraphicsExtensions{.eps}
\fi

\usepackage{natbib}
\usepackage{graphicx}
\usepackage{float}
\usepackage{subfig}
\usepackage{amsopn}
\usepackage{amsmath}
\usepackage{hyperref}
\usepackage{cleveref}
\usepackage{color,xcolor}
\usepackage{algorithm}
\usepackage{ifpdf}

\newcommand{\tpr}{\boldsymbol{\otimes}}   
\newcommand{\dpr}{\boldsymbol{\cdot}}     

\newcommand{\vc}[1]{\boldsymbol{#1}}      

\newcommand{\dudx}[2]{\frac{\partial{#1}}{\partial{#2}}}

\newtheorem{remark}{Remark}
\newtheorem{definition}{Definition}
\newtheorem{proof}{Proof}
\newtheorem{proposition}{Proposition}

\begin{document}

\begin{frontmatter}



\title{Diffuse interface relaxation model for two-phase compressible flows with diffusion processes}

%
%
%
%
%
%
%

\author{Chao Zhang$^{a,b}$}
\author{Igor Menshov$^{c,d}$}
\author{Lifeng Wang$^{a,e}$}
\author{Zhijun Shen$^{a,e}$}

\address[org1]{Institute of Applied Physics and Computational Mathematics, Beijing 100094, China}
\address[org2]{Lomonosov Moscow State University, Moscow 119991, Russia}
\address[org3]{Keldysh Institute for Applied Mathematics RAS, Moscow 125047, Russia}
\address[org4]{SRISA RAS, Moscow 117218, Russia}
\address[org5]{Center for Applied Physics and Technology, HEDPS, Peking University, Beijing 10087, China}

\begin{abstract}
The paper addresses a two-temperature model for simulating compressible two-phase flow taking into account diffusion processes related to the heat conduction and viscosity of the phases. This model is reduced from the two-phase Baer-Nunziato model in the limit of complete velocity relaxation and consists of the phase mass and energy balance equations, the mixture momentum equation, and a transport equation for the volume fraction. Terms describing effects of mechanical relaxation, temperature relaxation, and thermal conduction on volume fraction evolution are derived and demonstrated to be significant for heat conduction problems. The thermal conduction leads to instantaneous thermal relaxation so that the temperature equilibrium is always maintained  in the interface region with meeting the entropy relations. A numerical method is developed to solve the model governing equations that ensures the pressure-velocity-temperature (PVT) equilibrium condition in its high-order extension. We solve the hyperbolic part of the governing equations with the Godunov method with the HLLC approximate Riemann solver. The non-linear parabolic part is solved with an efficient Chebyshev explicit iterative method without dealing with large sparse matrices. To verify the model and numerical methods proposed, we demonstrate numerical results of several numerical tests such as the multiphase shock tube problem, the multiphase impact problem, and the planar  ablative Rayleigh–Taylor instability problem.
\end{abstract}

%
\newpage

\end{frontmatter}


\section{Introduction}
Numerical {\color{black}modeling} of compressible multiphase flow have found many applications in various natural, industrial and technological areas. Typical applications include bubble dynamics   \cite{Ranjan2011,quirk_karni_1996}, underwater explosion \cite{MILLER2013132,Kedrinsky2000,Holt1977,Shyue2014An}, cavitation flows  \cite{Lemartelot2014,saurel2009simple,saurel_petitpas_abgrall_2008},  multiphase flows in the porous rock  \cite{Balashov2019Numerical}, inertial confinement fusion \cite{2017Theoretical,Rinderknecht2018}, Rayleigh–Taylor \cite{sharp1983overview,kull1991theory,youngs1984numerical} and Richtmyer–Meshkov instabilities \cite{krechetnikov_2009,brouillette2002richtmyer,zhou2017rayleigh} and so on. In some problems where steep distributions of flow parameters occur, diffusion processes such as the heat conduction and viscous stress may have significant impact.  How to properly take into account these processes in multiphase hydrodynamics with resolved interfaces is the main issue of the present paper.

Numerical methods for simulating compressible multiphase flows can be generally classified into two categories  depending on the approach to resolve material interfaces: Diffuse interface methods (DIM)  \cite{Galina2021,SCHMIDMAYER2017468,Kemm2020A,F2019On,Dumbser2013,Dumbser2011A,saurel2018diffuse,saurel1999multiphase,saurel2009simple,saurel1999simple,abgrall1996prevent,Menshov2015On,menshov2018generalized,Coralic2014Finite,shyue1998efficient,SHYUE199943,Chiapolino2018} and the sharp interface methods (SIM) \cite{Kenamond2021A,Kikinzon2018Establishing,dyadechko2005moment,dyadechko2008reconstruction,kucharik2011hybrid,hirt1974arbitrary,glimm2003conservative,glimm1998three,mulder1992computing,osher2002level,fedkiw1999non}. 
The present work is done in the framework of the former -- DIM. Instead of explicitly tracking sharply resolved material interfaces as in SIMs, {\color{black}material interfaces in DIMs are captured by allowing a numerical diffusion zone of mixture flow that is modeled} as physical one. Thanks to these numerical diffusion, different components can be described with a unique set of partial differential equations and equation of state (EOS). Therefore, one can perform throughout computations on the Eulerian grid without specifying concrete interface locations. Moreover, DIMs avoid dealing with complicated grid movements and non-conservativeness issues.

The models for multiphase flows with resolved interfaces generally fall into two groups: {\color{black}One is based on the generalization of the conventional one-fluid Euler equations to multicomponent cases  \cite{abgrall1996prevent,shyue1998efficient,SHYUE199943,Johnsen2012Preventing,Alahyari2015,allaire2002five},  the other is} based on the reduction of non-equilibrium multi-phase flow models  \cite{BAER1986861,kapila2001two,saurel2009simple,pelanti2014mixture,murrone2005five}.  

{\color{black}The first group is more concerned with numerical aspects, in particular, the property to preserve the pressure-velocity equilibrium (the PV property), and also additionally temperature equilibrium (the PVT property) when the thermal conduction is also considered. These properties are used as important numerical condition or criterion to derive such models.} The definitions of these properties are given in \cref{subsec:Evolution_PVT}. For these models, material interfaces are represented by variable EOS parameters or by a characteristic function such as the Heaviside function that is interpreted as volume fraction in the context of multiphase flows.
A representative of {\color{black}these} models is the following model  \cite{shyue1999fluid,shyue1998efficient,2018Capuano} based on the PV property, which is formulated as
\begin{subequations}\label{eq:Abg}
\begin{align}
\dudx{\rho}{t} + \nabla\dpr(\rho \vc{u}) = 0,\label{eq:Abg:mass} \\
\dudx{\rho\vc{u}}{t} + \nabla\dpr\left(\rho\vc{u}\tpr\vc{u}\right) + \nabla  p  = 0, \label{eq:Abg:mom}\\
\dudx{ \rho E }{t} + \nabla\dpr\left[\left( \rho E + p \right) \vc{u} \right] = 0, \label{eq:Abg:en}\\
\frac{\partial}{\partial t} \left( \frac{1}{\gamma-1} \right) + \vc{u} \dpr \nabla \left( \frac{1}{\gamma-1} \right) = 0, \label{eq:Abg:gamf}\\
\frac{\partial}{\partial t} \left( \frac{\gamma p_{\infty}}{\gamma-1} \right) + \vc{u} \dpr \nabla \left( \frac{\gamma p_{\infty}}{\gamma-1} \right) = 0,\label{eq:Abg:gambf} \\
\frac{\partial \rho q}{\partial t} + \nabla \dpr \left( \rho q \vc{u} \right) = 0, \label{eq:Abg:q}
\end{align}
\end{subequations}
where $\rho, \; \vc{u}, \; p, \; E$ are the mixture density, velocity, pressure and specific total energy, respectively. The parameters $\gamma, \; p_{\infty}, \; q$ come from the EOS. Here, we consider the stiffened gas (SG) EOS for the $k$-th component that takes the following form:
\begin{subequations}
\begin{align}
\rho_k e_k = \frac{p_k + \gamma_k p_{\infty,k}}{\gamma_k - 1} + \rho_k q_k,\label{eq:eos_stiffened} \\
\rho_k e_k = \rho_k C_{v,k} T_k +  p_{\infty,k} + \rho_k q_k, \label{eq:eos_stiffened2}
\end{align}
\end{subequations}
where $C_{v,k}$ is the specific heat at constant volume. The parameters $\gamma_k$, $p_{\infty,k}$ and $q_k$ are constants characterizing the thermodynamic behaviours of the $k$-th  phase.

When thermal conduction is considered, the temperature becomes continuous at interfaces. However, Johnsen et al. \cite{Alahyari2015,Johnsen2012} pointed out that the {\color{black}system of equations} (\ref{eq:Abg}) does not preserve temperature equilibrium. 
Based on similar ideas as in designing the model (\ref{eq:Abg}) with the PV property, they proposed a method for defining the mixture EOS that ensures the PVT property. They add the following evolution equations for $C_v, \; p_{\infty}$ to the model \cref{eq:Abg}
\begin{subequations}
\begin{align}  
  \frac{\partial \rho C_v}{\partial t} + \nabla \dpr \left( \rho C_v \vc{u} \right) = 0,\label{eq:Abg:Cv}  \\
  \frac{\partial p_{\infty}}{\partial t} + \vc{u} \dpr \nabla p_{\infty} = 0.\label{eq:Abg:Pinf}
\end{align}
\end{subequations}

The evolved parameters obtained from \cref{eq:Abg:gamf,eq:Abg:gambf,eq:Abg:q} are used to compute the pressure, while those obtained from \cref{eq:Abg:Cv,eq:Abg:Pinf} to compute the temperature. This model can also be formulated in volume fraction framework by replacing all the evolution \cref{eq:Abg:gamf,eq:Abg:gambf,eq:Abg:q,eq:Abg:Cv,eq:Abg:Pinf} for EOS parameters with 
{\color{black}
\begin{subequations}
\begin{align}  
  \frac{\partial \rho Y_2}{\partial t} + \nabla \dpr \left( \rho Y_2 \vc{u} \right) = 0,\label{eq:Abg:Y}  \\
  \frac{\partial \alpha_2}{\partial t} + \vc{u} \dpr \nabla \alpha_2 = 0.\label{eq:Abg:alp}
\end{align}
\end{subequations}
where $\alpha_2$ and $Y_2$ are the volume fraction and mass fraction of the second component, respectively.
}

When velocity is spatially uniform, the internal energy is purely advected, 
\begin{equation}\label{eq:rhoe_adv}
\frac{\text{D} \rho e}{\text{D} t} = 0,
\end{equation}
where $\rho e = \sum{\alpha_k \rho_k e_k}$, {\color{black}the operator ${\text{D}{ \cdot }}/{\text{D} t}$ denotes the material derivative.}

By using \cref{eq:rhoe_adv}, the following mixture rules are proposed in  \cite{Alahyari2015,Johnsen2012}  to maintain the PVT property:
\begin{enumerate}
\item[•]
To maintain pressure equilibrium, the mixture EOS parameters are defined as
\begin{equation}\label{eq:mixeos1}
\frac{1}{\gamma-1} = \sum{\frac{\alpha_k}{\gamma_k-1}}, \quad \frac{\gamma p_{\infty}}{\gamma-1} = \sum{\frac{\alpha_k \gamma_k p_{\infty,k}}{\gamma_k-1}}, \quad  \rho q = \sum {\alpha_k \rho_k q_k}.
\end{equation}

\item[•]While to maintain temperature equilibrium, the following mixture EOS parameters should be defined as
\begin{equation}\label{eq:mixeos2}
\rho C_v = \sum {\alpha_k \rho_k C_{v,k}}, \quad p_{\infty} = \sum {\alpha_k p_{\infty,k}}, \quad  \rho q = \sum {\alpha_k \rho_k q_k}.
\end{equation}
\end{enumerate}

As can be noted, two different mixture rules are used for computing pressure and temperature, {\color{black}resulting in two different definitions for $p_{\infty}$ (and interface location when the fluid distribution is represented by their own $p_{\infty}$).} This ambiguity in mixture EOS definition also leads to difficulties in defining some thermodynamic variables, such as the mixture entropy. Therefore, the issue of consistency with the second law of thermodynamics is {\color{black}a key point to cause controversy}. In fact, the volume fraction based model consisting of \cref{eq:Abg:mass,eq:Abg:mom,eq:Abg:en,eq:Abg:Y,eq:Abg:alp} formally coincides with the five-equation model  \cite{allaire2002five} that {\color{black}lacks a mathematical entropy}. {\color{black}In the following we refer to this model with the mixture rules \cref{eq:mixeos1,eq:mixeos2} as the one-temperature five-equation model.}

Most of the second group models for simulating compressible multiphase flows come from the seven-equation Baer-Nunziato {\color{black}one}  \cite{BAER1986861}. In the original Baer-Nunziato model, each component is described by their own velocity, temperature, and pressure. However, for certain application scenarios such as the multiphase flows where each phase occupies its own volume, the physics included in the Baer-Nunziato model is not always necessary. Therefore, 
{\color{black}a variety of reduced models are proposed}, for example, the six-equation model with equilibrium velocity  \cite{kapila2001two,saurel2009simple,pelanti2014mixture}, the five-equation model with equilibrium velocity and equilibrium pressure  \cite{kapila2001two,murrone2005five} and the {\color{black}four-equation} model with equilibrium velocity, pressure and temperature  \cite{Lemartelot2014}. A complete hierarchy of {\color{black}these} models are formulated in  \cite{lund2012hierarchy}. Since these models are compatible with the complete Baer-Nunziato {\color{black}one}, they are more physically sound and reasonable. Besides, in  \cite{Balashov2018Quasi} a one-temperature quasi-hydrodynamic multiphase model with viscosity and heat conduction has been derived with the Coleman-Noll procedure \cite{Coleman1963TheTO}.

Among these models, the model with equilibrium temperature \cite{Lemartelot2014} is most appropriate to consider heat conduction process, however, it fails to ensure the PV or the PVT condition. Moreover, it does not provide topological information of the material interface, {\color{black}nor does it describe the evolution of the volume-fraction averaged material properties such as thermal conductivity and viscosity}.  Therefore, we are more interested in the temperature non-equilibrium models  \cite{kapila2001two,murrone2005five}. To the best of the author's knowledge, the work on implementing heat conduction in the framework of multi-temperature model is absent in literature. We aim to fill this gap in the present work. 

We build a two-temperature model based on the reduction of the Baer-Nunziato {\color{black}one}. {\color{black}The obtained model consists of two energy equations including thermal relaxation between phases driving temperatures into equilibrium. It includes viscosity, heat conduction and external energy source in each phase.} Note that the  heat conduction process is accompanied with instantaneous thermal relaxation so that temperature equilibrium is maintained. We demonstrate that the impact of these thermal relaxations  (which are usually neglected in the first group models) on volume fraction is significant. {\color{black}The obtained model ensures the pressure and the temperature equilibria during the heat conduction.}
We prove that the model agrees with the second law of thermodynamics.  Numerically, our model is proved to satisfy the PVT  property with a uniquely defined EOS.

We use the fractional step method  to solve the model. The solution procedure can be divided into four steps, i.e., the hyperbolic step, the viscous step, the thermal relaxation step and the heat conduction step. The homogeneous hyperbolic part is solved with the Godunov method coupled with the HLLC Riemann solver. The diffusion process (viscous step and heat conduction step) are governed by a set of parabolic partial differential equations. They are solved with an efficient method of local iterations, that allows much larger time step than the traditional explicit scheme and quite straightforward for parallel implementation. The thermal relaxation procedure is realized by solving a non-linear system with two variables (equilibrium temperature and volume fraction). We prove that the thermal relaxation procedure does not undermine the PVT property.

The rest of this article is organized as follows. In \Cref{sec:model}, we deduce a five-equation model and  a six-equation model, with more attention being devoted to the latter as it is more convenient for considering thermal processes in the multiphase system with phase energy equations. In \Cref{sec:numer_met}, we design numerical methods for solving the proposed model and prove some relevant properties. In \Cref{sec:numer_res}, numerical results of our model are presented and compared with those of other models.

\section{Model formulation}
\label{sec:model}
\subsection{The Baer-Nunziato type model}
The starting point of our model formulation is the complete Baer-Nunziato model \cite{BAER1986861} or its variant for compressible two-phase flows  \cite{saurel1999multiphase}. In this model each phase is assumed to behave as a pure fluid except when it interacts with the other fluid through relaxation terms.
Including viscosity, heat conduction and external energy source to the Baer-Nunziato model, we obtain the following formulation:
\begin{subequations}\label{eq:bn}
\begin{align}
  \label{eq:bn:mass}
  \dudx{\alpha_k\rho_k}{t} + \nabla\dpr(\alpha_k\rho_k\vc{u}_k) = 0, \\
  \label{eq:bn:mom}
  \dudx{\alpha_k\rho_k\vc{u}_k}{t} + \nabla\dpr\left(\alpha_k\rho_k\vc{u}_k\tpr\vc{u}_k + \alpha_k p_k \overline{\overline{I}} - \alpha_k \overline{\overline{\tau}}_k \right) = p_I \nabla \alpha_k \nonumber 
  \\ -  \overline{\overline{\tau}}_I \dpr \nabla {\alpha_k}
    + \mathcal{M}_k, \\
   \label{eq:bn:en}
  \dudx{\alpha_k \rho_k E_k}{t} + 
  \nabla\dpr\left[
\alpha_k\left( \rho_k E_k + p_k \right) \vc{u}_k   - \alpha_k \overline{\overline{\tau}}_k \dpr \vc{u}_k
  \right] 
  = p_I \vc{u}_I \dpr \nabla \alpha_k \nonumber \\
  - \vc{u}_{I} \dpr \left( \overline{\overline{\tau}}_I \dpr \nabla {\alpha_k} \right) 
  -  {\color{black}{p}_I} \mathcal{F}_k + \vc{u}_I \mathcal{M}_k 
 + \mathcal{Q}_k + q_k + \mathcal{I}_k, \\
   \label{eq:bn:vol}
  \dudx{\alpha_2 }{t} + 
  \vc{u}_I \dpr \nabla\alpha_2 = \mathcal{F}_2,
\end{align}
\end{subequations}
where the notations used are standard: $\alpha_k, \; \rho_k, \; \vc{u}_k, \; p_k, \; \overline{\overline{\tau}}_k, \; E_k$ are the volume fraction, density, velocity, pressure, viscous stress, and specific total energy of $k$-th component.

For viscous stress we use the  Newtonian approximation
\begin{equation}\label{eq:newton_vis}
\overline{\overline{\tau}}_k = 2\mu \overline{\overline{D}}_k + \left(\mu_{b,k} - \frac{2}{3}\mu_k \right) \nabla \dpr \vc{u}_k,
\end{equation}
where $\mu_k > 0$ is the coefficient of shear viscosity and  $\mu_{b,k} > 0$ is the coefficient of bulk viscosity, $\overline{\overline{D}}_k$ is defined as
\[
\overline{\overline{D}}_k = \frac{1}{2} \left[ \nabla \vc{u}_k + \left(  \nabla \vc{u}_k \right)^{\text{T}} \right].
\]

The total energy is $E_k  = e_k + \mathcal{K}_k$ where $e_k$, and $\mathcal{K}_k = \frac{1}{2}\vc{u}_k \dpr \vc{u}_k$ are the specific  internal energy and kinetic energy, respectively.

{\color{black}The inter-phase exchange terms include the velocity relaxation $\mathcal{M}_k$, the pressure relaxation $\mathcal{F}_k$, and the temperature relaxation $\mathcal{Q}_k$,}
\begin{equation}\label{eq:relaxations}
\begin{split}
\mathcal{M}_k=\vartheta\left( \vc{u}_{k^*}-\vc{u}_k \right), \quad \mathcal{F}_k=\eta\left( {p}_{k} - {p}_{k^*} \right), \quad \mathcal{Q}_k=\varsigma\left( T_{k^*} - T_{k} \right).
\end{split}
\end{equation}
where $k^{*}$ denotes the conjugate component of the $k$-th component,  i.e., $k=1,\; k^{*} =2$ or $k=2,\; k^{*} =1$. The relaxation rates are all positive $\;\vartheta > 0, \; \eta > 0, \; \varsigma > 0$.

The variables with the subscript ``I'' represent the variables at interfaces, for which there are several possible definitions  \cite{perigaud2005compressible,saurel2018diffuse}. Whatever the definitions we choose, 
{\color{black}\[\lim_{\eta\to\infty} p_{I} = \lim_{\eta\to\infty}p_{k} = p, \quad \lim_{\vartheta\to\infty}\vc{u}_{I} = \lim_{\vartheta\to\infty}\vc{u}_{k} = \vc{u}, \quad \lim_{\vartheta\to\infty} {\overline{\overline{\tau}}}_{I} = \lim_{\vartheta\to\infty}{\overline{\overline{\tau}}}_{k} = {\overline{\overline{\tau}}}.\] }

{\color{black}The heat conduction term is given as:
\begin{equation}\label{eq:q}
q_k = \nabla \dpr \left( \alpha_k \lambda_k \nabla T_k \right),
\end{equation}
and the external heat source term is written as:
\begin{equation}\label{eq:I}
\mathcal{I}_k = \alpha_k I_k,
\end{equation}
where $I_k$ denotes the the intensity of the external heat source released in the $k$-th phase, and $I_k \geq 0$.
}

For future use we can deduce the corresponding balance equations for phase internal energies and phase entropies from \cref{eq:bn}. The deduction procedure is similar to that in  \cite{zhangPHD2019,murrone2005five,kapila2001two,Zhang2019Eulerian, Zhang2020An} with the exception that we include viscosity, heat conduction, and external energy source here. We directly give the equation for the phase internal energy as follows:
\begin{equation} \label{eq:phase_int_en_eqn}
\begin{aligned}
\frac{\partial \alpha_{k} \rho_{k} e_{k}}{\partial t}+\nabla \cdot(\left.\alpha_{k} \rho_{k} e_k \boldsymbol{u}_{k} \right) = -\alpha_{k} p_{k} \nabla \cdot \boldsymbol{u}_{k}- {\color{black}{p}_{I}} {\mathcal{F}}_{k}+p_{I}\left(\boldsymbol{u}_{I}-\boldsymbol{u}_{k}\right) \cdot \nabla \alpha_{k} \\
+\left(\boldsymbol{u}_{I}-\boldsymbol{u}_{k}\right) \cdot \mathcal{M}_{k}+\mathcal{Q}_{k}+q_{k}+\mathcal{I}_{k} \\
+ \left( \boldsymbol{u}_{k} - \boldsymbol{u}_{I} \right)\dpr \left( \overline{\overline{\tau}}_I \dpr \nabla \alpha_k  \right) 
+ \alpha_k \overline{\overline{\tau}}_{k} : \overline{\overline{D}}_k.
\end{aligned}
\end{equation}

By using the Gibbs relation,
\begin{equation}\label{eq:gibbs}
T_k \text{d} s_k = \text{d} e_k - \frac{{\color{black}{p}_k}}{\rho_k^2} \text{d} \rho_k
\end{equation}
we further obtain
\begin{equation} \label{eq:phase_entropy_eqn}
\begin{array}{r}
T_k \left[ \dudx{\alpha_k \rho_k s_k}{t} + \nabla \dpr \left( \alpha_k \rho_k \boldsymbol{u}_k s_k \right) \right]  =  {\color{black}\left({p}_{k}-{p}_{I}\right)} {\mathcal{F}}_{k}   +  \left(p_{I}-p_{k}\right)\left(\boldsymbol{u}_{I}-\boldsymbol{u}_{k}\right) \cdot \nabla \alpha_{k} \\
+ \left(\boldsymbol{u}_{I}-\boldsymbol{u}_{k}\right) \cdot \mathcal{M}_{k}+ \mathcal{Q}_{k}+q_{k} + \mathcal{I}_k  \\
+ \left( \boldsymbol{u}_{k} - \boldsymbol{u}_{I} \right)\dpr \left( \overline{\overline{\tau}}_I \dpr \nabla \alpha_k  \right) +
\alpha_k \overline{\overline{\tau}}_{k} : \overline{\overline{D}}_k.
\end{array}
\end{equation}

{\color{black}Even though \cref{eq:bn} is the most complete model including relaxations in pressure, velocity and temperature, however, practical implementation of this model is rather complicated because of its complex wave structure and stiff relaxation procedures. Therefore, we will consider two possible reductions of this model that are given in the following sections.}

The Baer-Nunziato model is deduced by using the Coleman-Noll procedure  \cite{coleman1974thermodynamics,Coleman1963TheTO,BAER1986861}, keeping the second law of thermodynamics. Maintaining the physical consistency with the Baer-Nunziato model, the reduced models should also satisfy the second law of thermodynamics, as we demonstrate below.

\begin{remark}
{\color{black}The thermodynamically compatible two-phase compressible flow model proposed in  \cite{Romenski2007Conservative,Romenski2010Compressible,Romenski2012,Romenski2010} can also be reformulated in the form of Baer-Nunziato model with additional source terms describing the lift forces.}
\end{remark}

\begin{remark}
{\color{black}For turbulent bubbly flows, the viscous pressure has been proposed to consider the pulsation damping of the bubbles  \cite{saurel2003,perigaud2005compressible,Gavrilyuk2002}. Including this viscous pressure, the relaxing pressure is
\begin{equation}
\widetilde{p}_{k} = {p}_{k} + {p}_{\mu,k},
\end{equation}
with ${p}_{\mu,k}$ being the viscous pressure  \cite{Gavrilyuk2002}:
\begin{equation}\label{eq:vis_pres_def}
{p}_{\mu,k} = z_k (\alpha_k) \frac{\text{D}_{I} \alpha_k}{\text{D} t} = z_k (\alpha_k) \mathcal{F}_k,
\end{equation}
where $z_k$ is a function of $\alpha_k$, $\frac{\text{D}_{I} \cdot}{\text{D} t}$ denotes the material derivative related to the interface velocity $\mathbf{u}_{I}$, and $\mathcal{F}_k = \eta\left( \widetilde{p}_{k} - \widetilde{p}_{k^*} \right)$.

With the viscous pressure, the terms including $\mathcal{F}_k$ on r.-h.s. of \cref{eq:phase_int_en_eqn} and \cref{eq:phase_entropy_eqn} should be replaced by 
${\widetilde{p}_{I}} {\mathcal{F}}_{k}$ and 
$\left( {p}_{k}-\widetilde{p}_{I}\right) {\mathcal{F}}_{k}$, respectively.
In order that the term $\left( {p}_{k}-\widetilde{p}_{I}\right) {\mathcal{F}}_{k}$  makes a non-negative contribution to the phase entropy in \cref{eq:phase_entropy_eqn}, $\widetilde{p}_{I}$ should be a convex combination of $\widetilde{p}_{k}$ and  $z_k$ be non-positive.

It can be seen from \cref{eq:vis_pres_def} that the viscous pressure ${p}_{\mu,k}$ vanishes when pressure equilibrium is reached, thus, it has no impact on the solution of the reduced models (in \Cref{subsec:five_eqn,subsec:six_eqn}) derived in the limit of instantaneous mechanical relaxation. Therefore, we temporarily omit this term in the following discussions.
}
\end{remark}

\subsection{The reduced five-equation model}
\label{subsec:five_eqn}
By performing asymptotic analysis of the Baer-Nunziato model in the limit of instantaneous mechanical relaxations with the method  similar to  \cite{kapila2001two}, one can obtain the following system of equations:
\begin{subequations} \label{eq:reduced_five}
\begin{align}
\dudx{\alpha_k \rho_k}{t} + \nabla \dpr \left( {\alpha_k \rho_k \vc{u}} \right) = 0,  \label{eq:reduced_five:mass} \\
\dudx{\rho\vc{u}}{t} + \nabla\dpr\left(\rho\vc{u}\tpr\vc{u} + p \overline{\overline{I}} \right)  = 
\nabla\overline{\overline{\tau}}, \label{eq:reduced_five:mom} \\
 \dudx{\rho E}{t} +  \nabla\dpr \left( \rho E \vc{u} +  p \vc{u} \right) = \nabla \dpr \left( \overline{\overline{\tau}} \dpr \vc{u} \right) + \sum q_k + \sum \mathcal{I}_k, \label{eq:reduced_five:en} \\
\dudx{\alpha_2}{t} + \vc{u} \dpr \nabla \alpha_2  =  R_{p_2} + R_{q_2} + R_{Q_2}+ R_{I_2}, \label{eq:reduced_five:vol}
\end{align}
\end{subequations}
{\color{black}where $\overline{\overline{I}}$ is the unit tensor, $\rho = \sum \alpha_k \rho_k$ and $\overline{\overline{\tau}} = \sum \alpha_k \overline{\overline{\tau}}_k$ are the mixture density and the mixture viscous stress, respectively.}

The right hand side terms of  \cref{eq:reduced_five:vol} are 
\begin{subequations}
\begin{align}
&R_{p_2} = {{\alpha}_{2}}\frac{A-{{A}_{2}}}{{{A}_{2}}}\nabla\dpr{\vc{u}}, \quad
R_{q_2} =  A \frac{\Gamma_2 q_2 \alpha_1 - \Gamma_1 q_1 \alpha_2}{A_1 A_2}, \nonumber\\
&R_{Q_2} =  A \frac{\Gamma_2 \mathcal{Q}_2 \alpha_1 - \Gamma_1 \mathcal{Q}_1 \alpha_2}{A_1 A_2},\quad
R_{I_2} =  A \frac{\Gamma_2 \mathcal{I}_2 \alpha_1 - \Gamma_1 \mathcal{I}_1 \alpha_2}{A_1 A_2}. \nonumber
\end{align}
\end{subequations}
 where $\Gamma_k$ is the phase Gruneisen coefficient, $\Gamma_k=V_k \left( {\frac{\partial p_k}{\partial e_k}}\right)_{{V_k}}$, $V_k=1/\rho_k$, $a_k^2 = \overline{\gamma}_k p_k V_k$ is the phase speed of sound,  
 $\overline{\gamma}_k = -\frac{V_k}{p_k}\left(\frac{\partial p_k}{\partial V_k}\right)_{s_k}$ is the phase adiabatic exponent, $A_k=\rho_k a_k^2$, and
{\color{black}$1/A = \sum \left( \alpha_k / A_k \right)$}.

{\color{black}
In the case of the SG EOS (\ref{eq:eos_stiffened}), these parameters are
\begin{align}
\Gamma_k = \gamma_k -1, \\
\overline{\gamma_k} = \gamma_k \frac{p_k + p_{\infty,k}}{p_k} > \gamma_k > \Gamma_k, \label{eq:gamgam} \\
A_k  = \gamma_k (p_k + p_{\infty, k}) = \frac{\gamma_k \left( \gamma_k - 1 \right) C_{v,k} T_k}{V_k}. \label{eq:AK}
\end{align}
}
The first term on the right hand side (r.-h.s.) of \cref{eq:reduced_five:vol} $R_{p_k}$ comes from pressure relaxation. In fact, $W_{p_k} = -p R_{p_k}$ represents the rate of work performed on material interfaces to maintain pressure equilibrium under compression or expansion \cite{kreeft2010new}. The significance of this term for  spherical bubble dynamics and multiphase flows has been demonstrated in \cite{SCHMIDMAYER2020109080} and  \cite{murrone2005five}, respectively.

In the limit of the sharp material interface, i.e. $\alpha_k = 1, \; \alpha_l = 0 \; (l \neq k)$,  the first r.-h.s. term  of \cref{eq:reduced_five:vol} $R_{p_k}$ vanishes. {\color{black}The second term $R_{q_k}$ and the third term $R_{I_k}$  also vanish in accordance of the definitions \cref{eq:q} and \cref{eq:I}. However, the term $R_{Q_k}$  due to temperature relaxation still remains.} This means that for compressible multicomponent problem with heat conduction,  the thermal relaxation can not be neglected even for interface-tracking methods where the diffused zone is absent. Therefore, vanishing the r.-h.s. of \cref{eq:reduced_five:vol} and using just  the pure advection equation for volume fraction may lead to errors that come from physical defects instead of numerical ones.

If we define the mixture entropy as
\begin{equation}\label{eq:mixture_entropy}
s = Y_1 s_1 + Y_2 s_2,
\end{equation}
from  \cref{eq:phase_entropy_eqn} one can deduce 
\begin{equation}\label{eq:entropy_inequality}
\dudx{ \rho s}{t} + \nabla \dpr \left( \rho \boldsymbol{u} s \right) = \frac{\alpha_1 \overline{\overline{\tau}}_{1} : \overline{\overline{D}}_1}{T_1} + \frac{\alpha_2 \overline{\overline{\tau}}_{2} : \overline{\overline{D}}_2}{T_2} + \frac{\mathcal{Q}_{1}+q_{1} + \mathcal{I}_1}{T_1}  + \frac{\mathcal{Q}_{2}+q_{2} + \mathcal{I}_2}{T_2}.
\end{equation}
{\color{black}
\begin{proposition}
In the absence of heat flows through the external boundaries of the control volume, \cref{eq:entropy_inequality} is non-negative.
\end{proposition}
}
\begin{proof}\label{eq:proof5eqn_entropy}
The first two terms on the r.-h.s. of \cref{eq:entropy_inequality} can be proven to be non-negative, $\alpha_k \overline{\overline{\tau}}_{k} : \overline{\overline{D}}_k \geq 0$, by using the definition \cref{eq:newton_vis} and simple tensor manipulations.

Also, due to \cref{eq:relaxations}
\begin{equation}
\frac{\mathcal{Q}_{1}}{T_1} + \frac{\mathcal{Q}_{2}}{T_2} = \varsigma \frac{\left(T_2 - T_1 \right)^2}{T_1 T_2} \geq 0.
\end{equation}
 
The heat conduction term can be recast as
\begin{equation}
\frac{q_k}{T_k} = - \frac{\nabla \cdot \mathbf{q}_k}{T_k} =
 - \nabla \cdot (\frac{\mathbf{q}_k}{T_k}) + \mathbf{q}_k \cdot \nabla(\frac{1}{T_k})
\end{equation}

The first term is of the divergence type and represents external heat inflow to the phase material particle. The second term is positive due to  the Fourier's law of heat conduction $\vc{q}_k =  - \alpha_k \lambda_k \nabla T_k$, $\alpha_k \geq 0 , \; \lambda_k \geq 0$ and $\mathcal{I}_k \geq 0$. Therefore, except the heat inflow terms, the r.-h.s of the mixture entropy equation (2.11) is non-negative. This means that the mixture entropy respects the second law of thermodynamics.
\end{proof}

The model \cref{eq:reduced_five} assumes two temperatures and only one energy equation. {\color{black}In this paper, we prefer to consider} the energy exchanges and thermal conduction with the six-equation model that consists of two energy equations and physically consistent with the five-equation model \cref{eq:reduced_five}. 

\subsection{The reduced six-equation model}
\label{subsec:six_eqn}
We first separate the physical process into three stages: the mechanical stage, the thermal relaxation stage and the heat conduction stage, {\color{black}and then build thermodynamical consistency for each stage.}

{\color{black}
In the mechanical stage the pressure equilibrium is reached with the instantaneous pressure relaxation. Then thermal relaxation drives the phase temperatures to equilibrium. The heat conduction proceeds maintaining the obtained pressure equilibrium and temperature equilibrium.
}

\paragraph{Mechanical stage}
For the mechanical stage, we temporarily omit thermal relaxation and conduction. In the limit of instantaneous velocity relaxation, one can obtain the following six-equation model with one velocity from the Baer-Nunziato model \cref{eq:bn}
\begin{subequations}\label{eq:bn_six}
\begin{align}
\dudx{\alpha_k\rho_k}{t} + \nabla\dpr(\alpha_k\rho_k\vc{u}) = 0, \label{eq:bn_six:mass} \\
\dudx{\rho \vc{u}}{t} + \nabla\dpr\left(\rho \vc{u}\tpr\vc{u}\right) + \nabla \left( \alpha_1 p_1 +  \alpha_2 p_2 \right) = \nabla \dpr \overline{\overline{\tau}}, \label{eq:bn_six:mom} \\
\dudx{\alpha_k \rho_k e_k}{t} + \nabla\dpr\left( \alpha_k \rho_k e_k  \vc{u}  \right)  + \alpha_k p_k \nabla\dpr \vc{u}
  =  -  p_I \mathcal{F}_k + \alpha_k \overline{\overline{\tau}}_k : \overline{\overline{D}}, \label{eq:bn_six:en} \\
\dudx{\alpha_2 }{t} + \vc{u} \dpr \nabla\alpha_2 = \mathcal{F}_2. \label{eq:bn_six:vol}
\end{align}
\end{subequations}

The corresponding balance equation for mixture entropy is 
\begin{align}\label{eq:entropy_inequality_sixeqn}
\dudx{ \rho s}{t} + \nabla \dpr \left( \rho \boldsymbol{u} s \right) = \frac{\alpha_1 \overline{\overline{\tau}}_{1} : \overline{\overline{D}}_1}{T_1} + \frac{\alpha_2 \overline{\overline{\tau}}_{2} : \overline{\overline{D}}_2}{T_2} + {\color{black}\frac{\left(  {{p}}_1 - {{p}}_{I}\right)\mathcal{F}_{1}}{T_1}  + \frac{\left(  {{p}}_2 - {{p}}_{I}\right)\mathcal{F}_{2}}{T_2}}
\end{align}

{\color{black}As long as the interface pressure ${p}_{I}$ is assumed to be a convex combination of ${p}_1$ and ${p}_2$, i.e., 
\begin{equation}\label{eq:int_pres}
{p}_{I} = Z_1 {p}_1 + Z_2 {p}_2 \quad (Z_1, Z_2 \in [0,\,1], \; Z_1 + Z_2  = 1),
\end{equation} 
the term ${\left(  {{p}}_k - {{p}}_{I}\right)\mathcal{F}_{k}}$ remains non-negative and the second law of thermodynamics is respected. 
%
%
}


{\color{black}When solving internal energy equations (\ref{eq:bn_six:en}), the total energy equation (\ref{eq:reduced_five:en}) will have to be supplemented to keep
the energy conservation as in   \cite{saurel2009simple}.

This stage consists of the hydrodynamic, the viscous and the pressure relaxation processes. The last relaxation process drives the phase pressures into an equilibrium pressure $p = \lim_{\eta \to \infty} p_1 = \lim_{\eta \to \infty} p_2$.}

\paragraph{Thermal relaxation stage}
The procedure is similar to that for deducing the model for phase transition in \cite{zein2010}. {\color{black}Having reached the pressure equilibrium after the mechanical stage, we continue to build our model for the thermal relaxation on the basis of the following physical assumptions}:
\begin{enumerate}
\item[•] The mechanical relaxation happens much faster than the thermal relaxation, which means that the temperature relaxation goes in the state of pressure equilibrium.
\end{enumerate}

This assumption is reasonable for many applications, such as detonations and deflagration. The estimation analysis performed in  \cite{kapila2001two} demonstrates that the time scale for thermal relaxations are much larger than that for mechanical relaxations.

The thermal relaxation  process is assumed to be governed by the following equations:
\begin{subequations}
  \label{eq:thermal_relax}
\begin{align}
\dudx{\alpha_k\rho_k}{t} = 0, \label{eq:thermal_relax:mass} \\
\dudx{\rho \vc{u}}{t} = 0, \label{eq:thermal_relax:mom} \\
\dudx{\alpha_1 \rho_1 e_1}{t} =  \mathcal{Q}_1^{\prime}, \label{eq:thermal_relax:en1} \\
\dudx{\alpha_2 \rho_2 e_2}{t} =  \mathcal{Q}_2^{\prime}, \label{eq:thermal_relax:en2} \\
\dudx{\alpha_2 }{t} = r_{0} \frac{\mathcal{Q}_2^{\prime}}{p} . \label{eq:thermal_relax:vol}
\end{align}
\end{subequations}
{\color{black}where $Q_k^{\prime}$ is the thermal relaxation term defined in \cref{eq:relaxations}, which results in the variation of the phase temperature and the volume fraction.}  The term ${\mathcal{Q}_2^{\prime}}/{p}$ represents the volume fraction change rate if no phase temperature variation is considered. 

{\color{black}The parameter $r_0$ is a dimensionless coefficient, balancing the phase temperature change and volume fraction change. It is determined in such a way that the pressure equilibrium condition is maintained, i.e.,}
\begin{equation}\label{eq:pressure_equilibrium}
\frac{\partial p_1}{\partial t} = \frac{\partial p_2}{\partial t}.
\end{equation}

Thus one can obtain
\begin{equation} \label{eq:r0}
r_0 = \frac{\Gamma_1 / \alpha_1 + \Gamma_2 / \alpha_2}{ \left( A_1/\alpha_1 + A_2/\alpha_2  \right) / p -  \left( \Gamma_1 / \alpha_1 + \Gamma_2 / \alpha_2 \right) },
\end{equation}
or
\begin{equation} \label{eq:r0_}
r_0 = \frac{\Gamma_1 / \alpha_1 + \Gamma_2 / \alpha_2}{  \left( \overline{\gamma}_1 -  \Gamma_1\right)/ \alpha_1 + \left( \overline{\gamma}_2 - \Gamma_2\right)/ \alpha_2  },
\end{equation}

{\color{black}
According to \cref{eq:r0,eq:thermal_relax:mass,eq:AK}, $r_0$ is a function of $T_1$, $T_2$ and $\alpha_2$, and from \cref{eq:gamgam}, it satisfies
\begin{equation}\label{eq:r0_fun}
r_0 = r_0 \left( T_1, T_2, \alpha_2 \right) > 0.
\end{equation}}

By using \cref{eq:gibbs,eq:mixture_entropy,eq:r0_fun,eq:relaxations} one can deduce
\begin{equation}
\dudx{\rho s}{t} + \nabla \dpr \left( {\rho s \vc{u}} \right) = \left( 1 + r_0 \right) \left(\frac{1}{T_2} - \frac{1}{T_1} \right) \mathcal{Q}_2^{^{\prime}}  \geq 0 
\end{equation}

{\color{black}This means that the model for the thermal relaxation is consistent with the second law of thermodynamics.}



\paragraph{Heat conduction stage}
We build our model for thermal conduction under the following physical assumptions:
\begin{enumerate}
\item[•] The process of heat conduction goes under the condition of equilibrium in both pressure and temperature. 
\end{enumerate}

This assumption means that the heat conduction time scale is larger enough than the heat transfer scale so that temperature equilibrium always holds. In fact, this is a default assumption for models based on single temperature formulation, for example, the four-equation conservative model in  \cite{Lemartelot2014}. 

The heat conduction process (including the external heat source)  is modeled by the following system of equations:
\begin{subequations}
  \label{eq:ht}
\begin{align}
\dudx{\alpha_k\rho_k}{t} = 0, \label{eq:ht:mass} \\
\dudx{\rho \vc{u}}{t} = 0, \label{eq:ht:mom} \\
\dudx{\alpha_1 \rho_1 e_1}{t} =  \delta q_1 + q_1 + \mathcal{I}_1, \label{eq:ht:en1} \\
\dudx{\alpha_2 \rho_2 e_2}{t} =  \delta q_2 + q_2 + \mathcal{I}_2, \label{eq:ht:en2} \\
\dudx{\alpha_2 }{t} = \frac{r_{10}}{p} \delta q_2 + \frac{r_{1}}{p} \left( q_1 + \mathcal{I}_1 \right) + \frac{r_{2}}{p} \left( q_2  + \mathcal{I}_2\right). \label{eq:ht:vol}
\end{align}
\end{subequations}
{\color{black}Here, the term $\delta q_k$ represents an interphase heat conduction, and $\delta q_1 + \delta q_2 = 0$ for energy conservation. It plays a vital role  in maintaining thermodynamical consistency. In fact, without this interphase heat conduction the entropy inequality {\color{black}does not hold}.}

The interphase heat conduction is defined as a linear combination of $q_k$ and $\mathcal{I}_k$,
\begin{equation}\label{eq:linear_dq}
\delta q_2 = \widehat{r}_1 (q_1 + \mathcal{I}_1) + \widehat{r}_2 (q_2 + \mathcal{I}_2).
\end{equation}

{\color{black}We then define the parameters $r_{10}, \; r_1, \; r_2, \; \widehat{r}_1, \; \widehat{r}_2$ in this model in the following manner.}

\begin{enumerate}
\item[•] Defining the parameters $r_{10}, \; r_1, \; r_2$ 

{\color{black}The pressure  equilibrium condition \cref{eq:pressure_equilibrium} 
should be maintained, thus, one can obtain} 
\begin{subequations}
\begin{align}
r_{10} = r_{0}, \\
r_{1} = \frac{ - \Gamma_1 / \alpha_1}{  \left( \overline{\gamma}_1 -  {\Gamma}_1\right)/ \alpha_1 + \left( \overline{\gamma}_2 - \Gamma_2\right)/ \alpha_2  }, \\
r_{2} = \frac{ \Gamma_2 / \alpha_2}{   \left( \overline{\gamma}_1 -  {\Gamma}_1\right)/ \alpha_1 + \left( \overline{\gamma}_2 - \Gamma_2\right)/ \alpha_2 }.
\end{align}
\end{subequations}

\item[•] Defining the parameters $\widehat{r}_1, \; \widehat{r}_2$ 

{\color{black}The temperature equilibrium condition should be satisfied, i.e.,}
\begin{equation}\label{eq:temperature_equilibrium}
\frac{\partial T_1}{\partial t} = \frac{\partial T_2}{\partial t},
\end{equation}
which yields the coefficients in \cref{eq:linear_dq} as
\begin{subequations}
\begin{align}
\widehat{r}_1 = - \frac{r_1 \mathcal{Y} - m_2 C_{v,2} }{r_0 \mathcal{Y} - m_1 C_{v,1} - m_2 C_{v,2}}, \\
\widehat{r}_2 = - \frac{r_2 \mathcal{Y} + m_1 C_{v,1} }{r_0 \mathcal{Y} - m_1 C_{v,1} - m_2 C_{v,2}}, \\
\mathcal{Y} =  m_1 C_{v,1} G_{2} + m_2 C_{v,2} G_{1},
\end{align}
\end{subequations}
{\color{black}where $m_k = \alpha_k \rho_k$,  
$G_k = 1+\Gamma_k/c_{v,k}$, with $c_{v,k}$ being the dimensionless specific heat, $c_{v,k} = pV_k/(C_{v,k} T_k)$.}
\end{enumerate}


By using \cref{eq:ht} and the Gibbs relation \cref{eq:gibbs}, one can deduce
\begin{equation}\label{eq:entropy_ht}
\dudx{\rho s}{t} + \nabla \dpr \left( {\rho s \vc{u}} \right) = \frac{q_1}{T_1} + \frac{q_2}{T_2} + \frac{\mathcal{I}_1}{T_1} + \frac{\mathcal{I}_2}{T_2} + \left( \frac{1}{T_2} - \frac{1}{T_1} \right) \delta q_2 + \left ( \frac{p_1}{T_1} \frac{\text{D} \alpha_1}{\text{D} t} +  \frac{p_2}{T_2} \frac{\text{D} \alpha_2}{\text{D} t} \right).
\end{equation}

Since we maintain the pressure equilibrium \cref{eq:pressure_equilibrium} and temperature equilibrium \cref{eq:temperature_equilibrium} in this process, the last two terms in \cref{eq:entropy_ht} vanish. 
As for the first two terms, according to \cref{eq:proof5eqn_entropy}, we have
\begin{equation}
\int_{V} \frac{q_{k}}{T_{k}} \mathrm{~d} \tau \geq 0,
\end{equation}
if the net heat flux across the surface of the volume $V$ vanishes. 

We further deduce
\begin{equation}
\int_{V}\left[\frac{\partial \rho s}{\partial t}+\nabla \cdot(\rho s \boldsymbol{u})\right] \mathrm{d} \tau \geq 0,
\end{equation}
and as $V$ is an arbitrary closed domain, one can write 
\begin{equation}\label{eq:ent_neq1}
 \dudx{\rho s}{t} + \nabla \dpr \left( {\rho s \vc{u}} \right)  \geq 0.
\end{equation}

From \cref{eq:entropy_ht} we see that without temperature equilibrium, the entropy inequality \cref{eq:ent_neq1} may be violated.

\subsection{The final model}

We summarize the finial model for compressible two-phase flows with viscosity and heat conduction as follows:
\begin{subequations}
  \label{eq:bn_sixf}
\begin{align}
\dudx{\alpha_k\rho_k}{t} + \nabla\dpr(\alpha_k\rho_k\vc{u}) = 0, \label{eq:bn_sixf:mass}\\
\dudx{\rho \vc{u}}{t} + \nabla\dpr\left(\rho \vc{u}\tpr\vc{u}\right) + \nabla \left( \alpha_1 p_1 +  \alpha_2 p_2 \right) = \nabla \dpr \overline{\overline{\tau}}, \label{eq:bn_sixf:mom} \\
  \dudx{\alpha_k \rho_k e_k}{t} + 
  \nabla\dpr\left(
\alpha_k \rho_k e_k  \vc{u} 
  \right)  + \alpha_k p_k \nabla\dpr \vc{u}
  =  -  p_I \mathcal{F}_k + \alpha_k \overline{\overline{\tau}}_k : \overline{\overline{D}} \nonumber \\ + \mathcal{Q}^{{\prime}}_k + \delta q_k  + q_k + \mathcal{I}_k,  \label{eq:bn_sixf:en} \\
\dudx{\alpha_2 }{t} + \vc{u} \dpr \nabla\alpha_2 = \mathcal{F}_2 + \frac{r_0}{p} \left( \mathcal{Q}_2^{\prime} +  \delta q_2\right) \nonumber \\ + \frac{r_1}{p} \left( q_1 + \mathcal{I}_1 \right) + \frac{r_2}{p} \left( q_2 + \mathcal{I}_2 \right). \label{eq:bn_sixf:vol}
\end{align}
\end{subequations}

The mechanical stage  can violate the temperature equilibrium state of the phases that is reached and maintained through following temperature relaxations. One can see that the temperature relaxation $\mathcal{Q}_k$ in the considered model \cref{eq:bn} consists of two parts: the thermal relaxation $\mathcal{Q}_k^{\prime}$ and the phase heat conduction $\delta q_k$, with the former being much faster than the latter. The former ensures the initial temperature equilibrium before the heat conduction  progresses, the latter maintains this temperature equilibrium while the heat conduction  in and between the phases. Thus, temperature equilibrium is still maintained after the heat conduction. 

For each stage, the entropy inequality remains valid. Thus, after implementing the fractional step method corresponding to the three relaxation stages, the solution  obtained should not be contrary to the second law of thermodynamics.

Since the final model is non-conservative for the mixture total energy, we supplement it with the mixture total energy equation of the five-equation model \cref{eq:reduced_five:en} in order to correct the non-conservativeness. Similar idea is adopted in \cite{saurel2009simple} in the absence of the diffusion processes.

\section{Numerical methods}
\label{sec:numer_met}

In this section we describe the numerical methods for solving the above proposed model.
The numerical method is based on the operator splitting technique that consists of three stages: the mechanical (solving the hyperbolic, viscous, and the pressure relaxation parts of the equations), the temperature relaxation, and the heat conduction stage. We address these steps separately. 

\subsection{Hyperbolic part}
The homogeneous hyperbolic part of the governing equations \cref{eq:bn_sixf} to be solved first is as follows:
\begin{subequations}
  \label{eq:bn_six1}
\begin{align}
\dudx{\alpha_k\rho_k}{t} + \nabla\dpr(\alpha_k\rho_k\vc{u}) = 0, \label{eq:bn_six1:mass} \\
\dudx{\rho \vc{u}}{t} + \nabla\dpr\left(\rho \vc{u}\tpr\vc{u}\right) + \nabla \left( \alpha_1 p_1 +  \alpha_2 p_2 \right) = 0, \label{eq:bn_six1:mom} \\ 
\dudx{\alpha_k \rho_k e_k}{t} +  \nabla\dpr\left(\alpha_k \rho_k e_k  \vc{u}_k  \right)  + \alpha_k p_k \nabla\dpr \vc{u} =  0, \label{eq:bn_six1:en} \\
\dudx{\ \rho E}{t} + \nabla\dpr\left[ \left( \rho E  + \alpha_1 p_1 + \alpha_2 p_2 \right) \vc{u}  \right] = 0, \label{eq:bn_six1:ent} \\
\dudx{\alpha_2 }{t} + \vc{u} \dpr \nabla\alpha_2 = 0. \label{eq:bn_six1:vol}
\end{align}
\end{subequations}

As mentioned above, we adopt the idea similar to that  of \cite{saurel2009simple}, i.e.,  using a redundant equation for the mixture total energy \cref{eq:bn_six1:ent} to correct the solution of the non-conservative equations for phase internal energies  \cref{eq:bn_six1:en}.  One can rewrite  \cref{eq:bn_six1} (without the redundant equation) into the following system with respect to the primitive variable $\vc{Z} = \left[\rho_1 \;\; \rho_2 \;\; u \;\; v \;\; p_1 \;\; p_2 \;\; \alpha_2 \right]^{\text{T}}$
\begin{equation}
\dudx{\vc{Z}}{t} + \vc{A}\dudx{\vc{Z}}{x} = 0.
\end{equation}

It can be shown that the matrix $\vc{A}$ has 7 real eigenvalues (i.e. $u\pm c$ and $u$ of multiplicity 5) and  the corresponding set of six linearly independent right eigenvectors. Thus, the system is hyperbolic.

The mixture speed of sound for this model is
\begin{equation}
c^2  = Y_1 c_1^2 + Y_2 c_2^2.
\end{equation}

This ensures monotonic variation of the characteristic velocity across the interface zone and therefore more robust compared with the five-equation model \cref{eq:reduced_five} where the mixture speed of sound is given by non-monotonic Wood's formulae.  

We recast \cref{eq:bn_six1} in the vector compact form (in 1D) as:
\begin{equation}\label{eq:vectorForm}
\dudx{\vc{U}}{t} + \dudx{\vc{F}\left(\vc{U}\right)}{x} + \vc{R} \left( \vc{U} \right) \dudx{{u}}{x} = 0,
\end{equation}
{\color{black}where
\[
\vc{U} = \left[ \alpha_1 \rho_1 \;\; \alpha_2 \rho_2 \;\;  \rho u \;\; \rho v \;\; \alpha_1 \rho_1 e_1 \;\; \alpha_2 \rho_2 e_2 \;\; \rho E \;\; \alpha_2  \right]^{\text{T}},
\]
\[
\vc{F}\left(\vc{U}\right) = u \vc{U} + (\alpha_1 p_1 + \alpha_2 p_2) \vc{D},
\]
\[
\vc{D}\left(\vc{U}\right) = \left[ 0 \;\; 0 \;\; 1 \;\; 0 \;\; 0 \;\; 0 \;\; u \;\; 0 \right]^{\text{T}},
\]
\[
\vc{R}\left(\vc{U}\right) = \left[ 0 \;\; 0 \;\; 0 \;\;  0 \;\; \alpha_1 p_1 \;\; \alpha_2 p_2 \;\; 0 \;\; -\alpha_2  \right]^{\text{T}}.
\]}

We use the Godunov method coupled with the approximate Riemann solver HLLC to solve \cref{eq:vectorForm}:
\begin{equation}
\vc{U}_{i}^{n+1} = \vc{U}_{i}^{n} - \frac{\Delta t}{\Delta x} \left[  \vc{F}\left(\vc{U}_{i+1/2}^{*}\right) -\vc{F}\left(\vc{U}_{i-1/2}^{*}\right)  \right] - \frac{\Delta t}{\Delta x} \vc{R}\left(\vc{U}_{i}^{n}\right) \left( u_{i+1/2}^{*} -  u_{i-1/2}^{*} \right),
 \end{equation} 
where $\vc{U}_{i+1/2}^{*} = \vc{U}_{i+1/2}^{*} \left( \vc{U}_{i}, \vc{U}_{i+1} \right)$ is the Riemann solution at the cell face $i+1/2$. Here we use the the three-wave approximate Riemann solver HLLC  \cite{ShenZ2014a,ShenZ2014b,Shen2016A,Toro2009Riemann}.
{\color{black}The dimensional spitting method is used for extension to multiple dimensions.}

\subsection{Viscous part}
Viscous terms have no impact on the mass  balance equations and affect only the momentum and energy equations. The corresponding splitted equations are read as
\begin{equation}
\label{eq:viscous}
\begin{split}
\dudx{\alpha_k \rho_k}{t} = 0,\quad
\dudx{\alpha_k}{t} = 0, \quad
  \dudx{\rho\vc{u}}{t}  = \nabla \dpr \overline{\overline{\tau}}, \\
  \dudx{\alpha_k \rho_k e_k}{t}  = \alpha_k \overline{\overline{\tau_k}} : \overline{\overline{D}},\quad
 \dudx{\rho E}{t}  = \nabla \dpr \left( \overline{\overline{\tau}} \dpr \vc{u} \right).
 \end{split}
\end{equation}

To solve the parabolic PDE for velocity, we use an efficient method of local iterations  based on Chebyshev parameters \cite{Zhukov2010,SHVEDOV1998}. A brief introduction on this method is given below.

Consider the following 1D parabolic PDE 
\begin{equation}
\dudx{v}{t} = Lv + f\left( x, t \right), \quad x\in G \subset\mathbb{R} 
\end{equation}
where $L$ is a linear elliptic self-adjoint positive-definite operator. 

Given a grid $\Omega_h = \cup [{x_{j-1/2}, x_{j+1/2}}]$ with a space step $h$, consider also a discrete operator $L_h{v_j}$ that approximates the operator $L$ with $\mathcal{O}(h^2)$ on smooth solutions. For example, it can be the 1D reduction of  the 7-point (in 3D) symmetric discretization of $L$ obtained with the finite volume method used in the present paper (see below). The operator $L_h$ is self-adjoint and has real positive eigenvalues within an interval $[\lambda_{min},\lambda_{max}]$.

%
%
%
%
%
%

The method of local iterations \cite{Zhukov2010} is realized as $2P-1$ explicit iterations, where $P = \left \lceil \pi / 4 \sqrt{\tau \lambda_{max} + 1} \right \rceil$, with $\tau$ being the time step and $\left \lceil x \right \rceil$ denoting the maximal integer to be greater than or equal to $x$. These explicit iterations are written as follows (for details see \cite{Zhukov2010}):
\begin{equation}\label{eq:LIM}
v^{(m)} = \frac{1}{1 + \tau b_m} \left( v^{n} + \tau b_m v^{(m-1)} - \tau L_h v^{(m-1)} + \tau f^{(n)}\right), \quad m=1,2,\cdots, 2P-1,
\end{equation}
where $v^{(m)}$ is the solution after $m$-th iteration, $b_m$ is a set of iteration parameters, \[(b_1, b_2,\cdots, b_{2P-1})=(a_P,a_{P-1},\cdots,a_2,a_P, a_{P-1},\cdots,a_1),\]  
Here, 
\begin{equation}
a_m=\frac{\lambda_{max}}{1+\beta_1}(\beta_1-\beta_m), \quad  m=1,\cdots,P,
\end{equation}
and the sequence $(\beta_1,\cdots,\beta_P)$ represents the roots of the Chebyshev polynomial $T_P (x)$  :
${{\text{cos}} \frac{(2j-1)\pi}{2P}, \; j=1,\cdots,P},$
arranged in the increasing order.

Since $b_{2P-1} = 0$, the last iteration becomes
\begin{equation}\label{eq:LIM_last}
v^{(2P-1)} = v^{n} + \tau L_h v^{(2P-2)} + \tau f^{n},
\end{equation}
which is the pure explicit step and $v^{(2P-2)}$ can be viewed as a predicted solution.

This scheme ensures the monotonicity of the solution \cite{Zhukov2010}. 
Each explicit iteration of \cref{eq:LIM} is  a conventional explicit step, making its parallel realization quite straightforward.

According to \cref{eq:viscous}, the mixture density $\rho$ does not vary with time at this stage. Therefore, the momentum equation takes in 1D the following form:
\begin{equation}\label{eq:viscous1D}
 \rho\dudx{u}{t} = \frac{\partial}{\partial x} \left( \frac{4}{3} \mu \dudx{u}{x} \right).
\end{equation}

The above method of local iterations is applied to \cref{eq:viscous1D}. The operator $L_h$ that approximates the r.-h.s. is given by central differences as
\begin{equation}
 L_h = \frac{1}{\Delta x} \left( F^{vis}_{i+1/2} - F^{vis}_{i-1/2}\right),
 \end{equation} 
where \[F^{vis}_{i+1/2} = \frac{4}{3} \mu_{i+1/2} \dudx{{u}}{x} \big{|} _{i+1/2} \] represents the viscous flux across the cell face $i+1/2$. 

The last iteration step is given in the conservative form,
\[ \frac{(\rho u)_{i}^{n+1}-(\rho u)_{i}^{n}}{\Delta t} = L_h(\widehat{u}_i),\] with $\widehat{u}$ being the predicted velocity after the first $(2P-2)$ iterations.

Once the velocity is calculated, the total energy is then updated as follows:
\begin{equation}
 \frac{(\rho E)_{i}^{n+1} - (\rho E)_{i}^{n}}{\Delta t} = \frac{1}{\Delta x} \left( {\widehat{u}}_{i+1/2} F^{vis}_{i+1/2} - {\widehat{u}}_{i-1/2} F^{vis}_{i-1/2}\right).
\end{equation} 
where $F_{i+1/2}^{vis}$ is determined by the velocity $\widehat{u}$ calculated in the first $2P-2$ local iteration.

Note that $ \alpha_k \overline{\overline{\tau}} : \overline{\overline{D}} = \nabla \dpr \left( \alpha_k  \overline{\overline{\tau}}  \dpr \mathbf{u} \right) - \left[ \nabla \dpr \left( \alpha_k  \overline{\overline{\tau}} \right) \right] \dpr \mathbf{u} $, then one can update the internal energies as follows:
\begin{align}
 \frac{(\alpha_k \rho_k e_k)_{i}^{n+1} - (\alpha_k \rho_k e_k)_{i}^{n}}{\Delta t}& = \frac{1}{\Delta x} \left( {\alpha}_{k,i+1/2} {\widehat{u}}_{i+1/2} F^{vis}_{i+1/2} - {\alpha}_{k,i-1/2} {\widehat{u}}_{i-1/2} F^{vis}_{i-1/2}\right) \nonumber \\
&- \frac{1}{\Delta x} u_{i} \left( {\alpha}_{k,i+1/2} F^{vis}_{i+1/2} - {\alpha}_{k,i-1/2} F^{vis}_{i-1/2}\right).
\end{align}

Extensions of the above algorithm to multiple dimensions can be done straightforwardly  in the directional splitting manner.

\subsection{Pressure relaxation part}
Next step is to drive phase pressures into an equilibrium state by performing instantaneous pressure relaxation  procedures when $\tau = 1 / \eta \to 0$. The process can be described with the following equations:

\begin{equation}
  \label{eq:bn_prelax}
  \dudx{\alpha_k\rho_k}{t} = 0, \quad
  \dudx{\rho \vc{u}}{t} = 0, \quad
  \dudx{\alpha_k \rho_k e_k}{t} =  -  p_I \mathcal{F}_k,\quad
  \dudx{\alpha_k }{t} = \mathcal{F}_k,
\end{equation}
where $\mathcal{F}_k$ is defined in \cref{eq:relaxations}.

{\color{black}Here we use the relaxation algorithm proposed in \cite{saurel2009simple}. This algorithm consists of the following basic steps:
\begin{enumerate}
\item[(1)] Combining \cref{eq:bn_prelax} and \cref{eq:eos_stiffened}, one can obtain the relaxed volume fraction as a function of the equilibrium pressure ${p}^{(1)}$, i.e., $\alpha_k^{(1)} = \alpha_k({p}^{(1)})$. By using the saturation constraint $\sum \alpha_k({p}^{(1)}) = 1$, we can find ${p}^{(1)}$ and $\alpha_k^{(1)}$.
\item[(2)] Having $\alpha_k^{(1)}$, we then re-evaluate the pressure by using the mixture total energy $\rho E$ (solved from the mechanical part of \cref{eq:reduced_five:en}) to ensure the conservativeness of energy and obtain the final pressure as $p^{(2)} = p(\alpha_k^{(1)}, \rho e)$, where $\rho e = \rho E - \rho \mathbf{u} \cdot \mathbf{u} / 2$.
\item[(3)] The phase internal energies are recalculated according to $e_k = e_k(p^{(2)}, \alpha_k^{(1)})$.
\end{enumerate}}

It is reported that this solution algorithm turns to be only about 5\% more expensive than that of the five-equation model  \cite{SCHMIDMAYER2020109080}.

\subsection{Temperature relaxation and heat conduction parts}
\label{subsec:temp_relax_ht}
%


The system of equations for the temperature relaxation read:
\begin{subequations}\label{eq:Trelax}
\begin{align}
\dudx{\alpha_k\rho_k}{t} = 0, \label{eq:Trelax:mass}\\
\dudx{\rho \vc{u}}{t} = 0, \label{eq:Trelax:mom}\\
\dudx{\alpha_k \rho_k e_k}{t} =  \mathcal{Q}_k^{\prime}, \label{eq:Trelax:en}\\
\dudx{\alpha_2 }{t} + \vc{u} \dpr \nabla\alpha_2 = \frac{r_0}{p} \mathcal{Q}_2^{\prime}. \label{eq:Trelax:vol}
\end{align}
\end{subequations}

And the heat conduction process is described by 
\begin{subequations}\label{eq:ht1}
\begin{align}
\dudx{\alpha_k\rho_k}{t} = 0, \label{eq:ht1:mass} \\
\dudx{\rho \vc{u}}{t} = 0, \label{eq:ht1:mom} \\
\dudx{\alpha_k \rho_k e_k}{t} =  {\delta q}_k + q_k + \mathcal{I}_k , \label{eq:ht1:en}\\
\dudx{\alpha_2 }{t} = \frac{r_0}{p} {\delta q}_2 + \frac{r_1}{p} \left( q_1 + \mathcal{I}_1\right) + \frac{r_2}{p} \left( q_2 + \mathcal{I}_2 \right). \label{eq:ht1:vol}
\end{align}
\end{subequations}

We see that formally  \cref{eq:Trelax} is a particular case of  \cref{eq:ht1} when $q_k + \mathcal{I}_k = 0, \; \delta q_k = \mathcal{Q}_k^{\prime}$. Therefore, we first deal with numerical solutions of \cref{eq:ht1} and then extend to \cref{eq:Trelax}.

Considering $e_k = e_k(T_k, \rho_k)$  and eliminating  $\delta q_2$ from \cref{eq:ht1}, one can deduce the following relation between phase temperatures and volume fraction $\alpha_2$:
%
%
%
%
%
%
%
%
%
\begin{equation}\label{eq:ht_eq2}
\mathcal{C}_2 \dudx{T_2}{t} + \mathcal{B}_2\dudx{\alpha_2}{t} = \mathcal{R}_1 \left( q_1 + \mathcal{I}_1 \right) + \mathcal{R}_2 \left( q_2 + \mathcal{I}_2\right).
\end{equation}
%
\begin{align}\label{eq:ht_eq3}
&\mathcal{B}_2 \mathcal{C}_1\dudx{T_1}{t} + \left( \mathcal{B}_1 + \mathcal{B}_2 \right) \mathcal{C}_2 \dudx{T_2}{t}  = \left( \mathcal{B}_2 + \mathcal{B}_1 \mathcal{R}_1 \right) \left( q_1 + \mathcal{I}_1 \right) + \left( \mathcal{B}_2 + \mathcal{B}_1 \mathcal{R}_2 \right) \left( q_2 + \mathcal{I}_2\right).
\end{align}
{\color{black}where 
\[\mathcal{C}_k =m_k C_{v,k}, \] 
\[\mathcal{B}_1 \left(T_k, \alpha_2 \right) = p\left(T_k, \alpha_2 \right) (G_2-G_1), \]
\[ \mathcal{B}_2\left(T_k, \alpha_2 \right)=p\left(T_k, \alpha_2 \right) \frac{1-r_0\left(T_k, \alpha_2 \right) G_2}{r_0\left(T_k, \alpha_2 \right)}, \]
\[ \mathcal{R}_1 \left( \alpha_2 \right) = \frac{- \Gamma_1 / \alpha_1}{\Gamma_1 / \alpha_1 + \Gamma_2 / \alpha_2}, \] 
\[\mathcal{R}_2 \left( \alpha_2 \right) = 1+\frac{ \Gamma_2 / \alpha_2}{\Gamma_1 / \alpha_1 + \Gamma_2 / \alpha_2}.\]
}

In the case of the SG EOS, we have:
\begin{align}
\mathcal{B}_1 = p_{\infty,1} - p_{\infty,2},\\
\mathcal{B}_2 = \frac{p}{r_{0}} + p_{\infty,2}. \label{eq:B2}
\end{align}

It can be seen that $\mathcal{C}_1, \;\mathcal{C}_2, \;\mathcal{B}_1$ are all constants  in this case, while $\mathcal{B}_2$ is a function of $\alpha_2, \; T_1, \; T_2$ due to \cref{eq:r0_fun}) and $p = p(T_k, \rho_k) = p(T_k, \frac{m_k}{\alpha_k})$. Here, $m_k$ is constant as a result of \cref{eq:Trelax:mass}.

{\color{black}
Using \cref{eq:r0,eq:AK,eq:B2}, $\mathcal{B}_2(T_1, T_2, \alpha_2)$ can be explicitly written as
\begin{equation}\label{eq:B2_2}
\mathcal{B}_2(T_1, T_2, \alpha_2) = \frac{\Gamma_1 \mathcal{C}_1 T_1 / \alpha_1^2 + \Gamma_1 p_{\infty,1}/\alpha_1 + \Gamma_2 \mathcal{C}_2 T_2 / \alpha_2^2 + \Gamma_2 p_{\infty,2}/\alpha_2}{\Gamma_1 / \alpha_1 + \Gamma_2 / \alpha_2} + p_{\infty,2}.
\end{equation}
}

\subsubsection{Temperature relaxation}\label{subsubsec:TemperatureRelaxation}
The temperature relaxation process is assumed to be much faster than phase heat conduction so that we take $q_1 = q_2 =  0$ and $\mathcal{I}_1 = \mathcal{I}_2 =  0$. In this case, \cref{eq:ht_eq2,eq:ht_eq3} are reduced to the following: 
\begin{subequations}
\begin{align}
\mathcal{C}_2 \dudx{T_2}{t} + \mathcal{B}_2 \dudx{\alpha_2}{t} = 0, \label{eq:ht_eq2_TR} \\
\mathcal{B}_2 \mathcal{C}_1\dudx{T_1}{t} + \left( \mathcal{B}_1 + \mathcal{B}_2 \right) \mathcal{C}_2 \dudx{T_2}{t} = 0. \label{eq:ht_eq3_TR}
\end{align}
\end{subequations}

In the model considered, we neglect a finite temperature relaxation time and assume the temperature equilibrium to occur within the time step. Using the superscript ``0'' and ``$\prime$'' to denote parameters before and after the temperature relaxation stage, an implicit discretization of \cref{eq:ht_eq2_TR,eq:ht_eq3_TR} can be written as 
{
\begin{subequations}
\begin{align}
\mathcal{C}_2 \left(  T^{\prime} - T_2^{0} \right) + \mathcal{B}_2(T_1^{av}, T_2^{av}, {\alpha}_2^{av}) \left(  {\alpha_2}^{\prime} - {\alpha_2}^{0} \right) = 0, \label{eq:ht_eq2_TR_dis}\\
\mathcal{B}_2(T_1^{av}, T_2^{av}, {\alpha}_2^{av}) \mathcal{C}_1 \left(  T^{\prime} - T_1^{0} \right) + \left( \mathcal{B}_1  + \mathcal{B}_2(T_1^{av}, T_2^{av}, {\alpha}_2^{av}) \right) \mathcal{C}_2 \left(  T^{\prime} - T_2^{0} \right)  = 0. \label{eq:ht_eq3_TR_dis}
\end{align}
\end{subequations}
}
Here the parameters $\mathcal{C}_1, \; \mathcal{C}_2$, and $\mathcal{B}_1$ are all constants, while $\mathcal{B}_2$ is a function of the phase temperatures and the volume fraction, $\mathcal{B}_2=\mathcal{B}_2(T_1, T_2,  \alpha_2)$ that is approximated by the average values $T_k^{av} = {(T_k^0+T^{\prime})}/{2}$ and ${\alpha}_2^{av} = {(\alpha_2^0+\alpha_2^{\prime})}/{2}$, i.e., $\mathcal{B}_2=\mathcal{B}_2(T_1^{av}, T_2^{av}, {\alpha}_2^{av})$.
This system is solved with the Newton method or the simple iterative method. In the present work we use the latter.

{\color{black}
\begin{remark}
If we look at \cref{eq:mixeos2} from the perspective of the temperature relaxation, the relaxed temperature defined by the one-temperature five-equation  model can be viewed as an averaged temperature:
\begin{equation}
\rho e = \sum \left( \alpha_k \rho_k C_{v,k} T_k^{0} +  \alpha_k p_{\infty,k} + \alpha_k \rho_k q_k\right) = \rho C_v T^{\prime} + p_{\infty} + \rho q,
\end{equation}
with
\begin{equation}\label{eq:Taverage}
T^{\prime} = \frac{\mathcal{C}_1 T_1^{0} + \mathcal{C}_2 T_2^{0}}{\mathcal{C}_1  + \mathcal{C}_2 }.
\end{equation}
In fact, in the case when $p_{\infty,1} = p_{\infty,2}$ the solution of \cref{eq:ht_eq3_TR} coincides with \cref{eq:Taverage}. Otherwise, we obtain a solution different from \cref{eq:Taverage}. Moreover, no corresponding volume fraction variation is considered in the  one-temperature five-equation  model.
\end{remark}
}

\subsubsection{Heat conduction}
The heat conduction process  goes under the temperature equilibrium condition $T_1 = T_2 = T$,  so that \cref{eq:ht_eq2,eq:ht_eq3} describe the change in time of temperature and volume fraction:  
\begin{subequations}
\begin{align}\label{eq:ht_eq2_ht}
\dudx{T}{t} = \mathcal{V}_1 \left( q_1 + \mathcal{I}_1\right) + \mathcal{V}_2 \left( q_2 + \mathcal{I}_2\right),\\ \label{eq:ht_eq3_ht}
\dudx{\alpha_2}{t} = \mathcal{U}_1 \left( q_1 + \mathcal{I}_1\right) + \mathcal{U}_2 \left( q_2 + \mathcal{I}_2\right),
\end{align}
\end{subequations}
where
\[
\mathcal{V}_1 = \frac{\mathcal{B}_2 + \mathcal{B}_1 \mathcal{R}_1}{\left( \mathcal{B}_1 + \mathcal{B}_2 \right) \mathcal{C}_2 + \mathcal{C}_1 \mathcal{B}_2}, \quad
\mathcal{V}_2 = \frac{\mathcal{B}_2 + \mathcal{B}_1 \mathcal{R}_2 }{\left( \mathcal{B}_1 + \mathcal{B}_2 \right) \mathcal{C}_2 + \mathcal{C}_1 \mathcal{B}_2},
\]
and
\[
\mathcal{U}_1 = \frac{\mathcal{R}_1-\mathcal{C}_2 \mathcal{V}_1}{\mathcal{B}_2}, \quad
\mathcal{U}_2 = \frac{\mathcal{R}_2-\mathcal{C}_2 \mathcal{V}_2}{\mathcal{B}_2}.
\]

Initial data for this system of ODE are $T^\prime$ and $\alpha^\prime$ obtained as the result of solving the temperature relaxation step (see \cref{subsubsec:TemperatureRelaxation}). 


Note that the coefficients $\mathcal{V}_1, \;\mathcal{V}_2, \;\mathcal{U}_1, \; \mathcal{U}_2$ are functions of $T$ and $\alpha_2$. The heat conduction coefficients commonly depend on temperature, i.e., $\lambda_k = \lambda_k (T_k)$. For example, for the thermal conductivity in completely ionized gas  $\lambda_k = \mathcal{O} (T_k^{\frac{5}{2}})$ \cite{Spitzer1953}. Therefore, \cref{eq:ht_eq2_ht,eq:ht_eq3_ht} represent a system of non-linear PDEs, with the spatial differential operator being applied  only to $T$.


To solve this system of parabolic equations we implement the method of local iterations described above (\cref{eq:LIM}). 
The term due to heat conduction $q_k$ is approximated with the central difference scheme. For example, assuming the 1D case on a uniform grid, $q_k$ is discretized as
\begin{equation}\label{eq:qk_dis}
\left( q_k \right)_{i} = \frac{\Lambda_{k,i+1/2}  T_{i+1} - \left( \Lambda_{k,i+1/2} + \Lambda_{k,i-1/2} \right)  T_{i} + \Lambda_{k,i-1/2}  T_{i-1} }{\Delta x ^2}
\end{equation}
where $\Lambda_{k} = \alpha_k \lambda_k$.

The method of local iterations is applied to solve \cref{eq:ht_eq2_ht} for temperature with iterative recalculation of volume fraction in \cref{eq:ht_eq3_ht}. The computational algorithm is formulated in \Cref{alg:explicit}.

\begin{algorithm} 
\caption{The iterative algorithm for solving \cref{eq:ht_eq2_ht,eq:ht_eq3_ht}}
\label{alg:explicit}
\begin{algorithmic}
\STATE{Define the discretized solution \\ $\mathbb{T}:=\{T_{1}, T_{2}, \ldots, T_{N}\}, \quad \mathbb{A} := \{{\alpha}_{2,1}, {\alpha}_{2,2}, \ldots, {\alpha}_{2,N} \}, \quad \mathbb{T}^{(1)}:=\mathbb{T}^{\prime}, \quad \mathbb{A}^{(1)}:=\mathbb{A}^{\prime}$ }
\STATE{Define $it:=1, \quad Conv := -1, \quad tol$ }
\WHILE{$Conv < 0$}
\STATE{Calculate the parameters $\mathcal{V}_k, \;\mathcal{U}_k, \;\Lambda_k$ by using $\mathbb{T}^{(it)}, \; \mathbb{A}^{(it)}$}
\STATE{Solve \cref{eq:ht_eq2_ht} with respect to $T$ by using the method  of local iterations (or a conventional implicit scheme) to obtain $\mathbb{T}^{\prime}$ }
\STATE{Solve \cref{eq:ht_eq3_ht} with respect to $\alpha_2$ with $\mathbb{T}^{\prime}$ and $\mathbb{A}^{(it)}$ to obtain $\mathbb{A}^{\prime}$}
\STATE{Set $\mathbb{T}^{(it+1)} = \mathbb{T}^{\prime}$, \quad $\mathbb{A}^{(it+1)} = \mathbb{A}^{\prime}$}
\STATE{Calculate $err = \Vert \mathbb{T}^{(it+1)} - \mathbb{T}^{(it)}  \Vert$}
\IF{$err < tol$}
\STATE{$Conv$  = 1}
\ENDIF
\STATE{Update $it := it + 1$}
\ENDWHILE
\RETURN $\mathbb{T}^{(it+1)},\quad \mathbb{A}^{(it+1)}$
\end{algorithmic}
\end{algorithm}

\subsection{Evolution of constant pressure and temperature profiles}
\label{subsec:Evolution_PVT}
For the interface-capturing schemes, an important property is the preservation of constant velocity and pressure profiles, which is referred to as the PV property in literature and given by the following definition.

\begin{definition}
Say that an interface-capturing numerical scheme has the PV property if it ensures  \[u_{i}^{n+1} = u = \text{const}, \; p_{i}^{n+1} = p = \text{const}\] providing that 
\[u_{i}^{n} = u = \text{const}, \; p_{i}^{n} = p = \text{const}.\]
\end{definition}

The numerical methods/models with this property have been studied, for example, in  \cite{abgrall1996prevent,Abgrall2001,shyue1998efficient,SHYUE199943,Menshov2018,Zhangchao2020,allaire2002five}. However, as pointed out in  \cite{Johnsen2012Preventing,Alahyari2015}, the methods with the PV property may result in erroneous temperature spikes in the vicinity of the material interfaces. This phenomenon is not problematic when dissipative processes are not considered. However, when heat conduction is involved, the numerical errors in temperature may affect the pressure through the energy equation. Therefore, for compressible multi-fluid problems, instead of the above PV property we require the following PVT property

\begin{definition}
An interface-capturing numerical scheme has the PVT property if it ensures  \[u_{i}^{n+1} = u = \text{const}, \; p_{i}^{n+1} = p = \text{const}, \; T_{i}^{n+1} = T = \text{const},\]
providing that
\[u_{i}^{n} = u = \text{const}, \; p_{i}^{n} = p = \text{const}, \; T_{i}^{n} = T = \text{const}\]
\end{definition}

Johnsen et al.  \cite{Johnsen2012Preventing,Alahyari2015} have proposed a methodology to get rid of the temperature spikes by introducing rules to define different mixture EOS for computing pressure and temperature. Their idea is similar to that of  \cite{abgrall1996prevent} for designing numerical methods to ensure the PV property. They developed their method based on the one-fluid formulation with single velocity, pressure and temperature. In this model the interfaces are represented by discontinuity in material properties. 

However, this method may result in multiple definitions of material properties, and thus ambiguity in interface locations. In fact, although they assume that the fluids are in temperature equilibrium, the resultant model formally allows two temperatures. Their definitions of the mixture EOS for computing temperature is equivalent to averaging the phasic temperatures according to \cref{eq:Taverage}. 

If we look at the problem from the perspective of the two-temperature model, the temperature averaging procedure (by defining the mixture EOS) should be interpreted as a physical process -- temperature relaxation. The impact of temperature relaxation process on volume fraction evolution is significant, as we demonstrate below. In the model of \cite{Johnsen2012Preventing,Alahyari2015} this impact is neglected and volume fraction is purely advected. As can be seen in our model formulation (see \Cref{sec:model}), the impact of temperature relaxations ($\mathcal{Q}_k^{\prime}$ in \cref{eq:thermal_relax} and $\delta q_2$ in \cref{eq:ht}) on volume fraction evolution has been included and numerically treated properly in \cref{subsec:temp_relax_ht}.

In the case of ideal gas EOS, we have
\[
p_{\infty,k} = 0, \; q_k = 0, \; \mathcal{B}_1 = 0.
\]
Then the solution of our temperature-relaxations equations (\ref{eq:ht_eq3_TR}) reproduce \cref{eq:Taverage}. In fact, as long as $p_{\infty,1} = p_{\infty,2}$,  \cref{eq:Taverage} holds. 
If the phasic temperatures before thermal relaxation are in equilibrium, then the averaging procedure does not change the temperature, nor the volume fraction.

Next we demonstrate that the proposed model preserves the PVT property, and is free of the temperature spike problem. Let us consider the following Riemann problem with the initial discontinuity:
\begin{equation}\label{eq:RiemannInitial}
\begin{array}{l}
u^{L}=u^{R}=u > 0 \\
\rho_{k}^{L}=\rho_{k}^{R}=\rho_{k}, \quad k=1,2 \\
e_{k}^{L}=e_{k}^{R}=e_{k}, \quad k=1,2 \\
\alpha_2^{L} \neq \alpha_2^{R}, \\
T_1 = T_2 = T,\\
{\gamma}_1 \neq {\gamma}_2.
\end{array}
\end{equation}

This problem is similar to that in  \cite{allaire2002five,Zhangchao2020}, the difference consists in that we additionally require an initial temperature equilibrium and consider the thermal relaxation process. 

\begin{proposition}
The solution to our model equations with initial discontinuity (\ref{eq:RiemannInitial}) ensures that
\begin{equation}
\begin{array}{l}
u^{*}=u \\
\rho_{k}^{*}=\rho_{k}, \quad k=1,2 \\
e_{k}^{*}=e_{k}, \quad k=1,2 \\
T_{1}^{*} = T_{2}^{*} = T.
\end{array}
\end{equation}
where the superscript ``*'' denotes the solution in the cell downstream the discontinuity after one time step.
\end{proposition}
\begin{proof}\label{proofPVT}
We apply a Riemann solver that resolves isolated contact discontinuity exactly (for example HLLC  \cite{Toro2009Riemann,Shen2016A}). After one time step, we have
\begin{equation}
\vc{U}_{i}^{*} = \xi \vc{U}_{i-1} + (1-\xi)\vc{U}_{i},
\end{equation}
where $\vc{U}$ is the solution vector defined in \cref{eq:vectorForm} and $\xi = u \Delta t / \Delta x$.
After some algebraic manipulations, one can obtain that 
\begin{equation}
u^{*} = u, \quad p_1^{*} = p_2^{*} = p, \quad {e}_k^{*} = {e}_k, \quad {\rho}_k^{*} = {\rho}_k.
\end{equation}

By using the EOS of each phase, one can deduce
\begin{equation}
T_k^{*} = T_k(p_k^{*}, \rho_k^{*}) = T. 
\end{equation}

Next we prove that the temperature relaxation \cref{eq:ht_eq2_TR,eq:ht_eq3_TR} with $T_k^{0} = T_k^{*}$ allows only one physically admissible solution $T^{\prime} = T$. 

For the case $p_{\infty,1} = p_{\infty,2}$, this consequence immediately comes from \cref{eq:Taverage}. 

For the case $p_{\infty,1} \neq p_{\infty,2}$, the proof is not so straightforward. For this case, let us assume that there exists another solution that $T^{\prime\prime} \neq T$ and satisfies \cref{eq:ht_eq2_TR,eq:ht_eq3_TR}. 
{\color{black}
By using \cref{eq:ht_eq3_TR_dis} and having in mind $T^{\prime\prime} - T \neq 0$, one obtains
\begin{equation}\label{eq:B2_}
\mathcal{B}_2 \left( T^{av}, \alpha_2^{av} \right) = - \frac{\mathcal{C}_2 \mathcal{B}_1}{\mathcal{C}} = B_2 = \textrm{const},
\end{equation}
where \[\mathcal{C} = \mathcal{C}_1 + \mathcal{C}_2,\] \[T^{av} = T_1^{av} = T_2^{av} = (T^{\prime\prime} + T) / 2,\] \[\alpha_2^{av} = (\alpha_2^{\prime\prime} + \alpha_2) / 2.\]

By using \cref{eq:B2_,eq:B2_2}, one can obtain 
\begin{equation}\label{eq:Tav_eq1}
T^{av} = \frac{b_1 \Gamma_1/\alpha_1^{av} + b_2 \Gamma_2/\alpha_2^{av}}{ \Gamma_1 {\mathcal{C}}_1/{(\alpha_1^{av})}^2+\Gamma_2 {\mathcal{C}}_2/{(\alpha_2^{av})}^2},
\end{equation}
where $b_1 = {B}_2 - p_{\infty,2} - p_{\infty,1}$, $b_2 = {B}_2 - 2p_{\infty,2}$, $\alpha_1^{av} = 1 - \alpha_2^{av}$.

By using \cref{eq:ht_eq2_TR_dis,eq:B2_}, one obtains
\begin{equation}
T^{\prime\prime} = T - 2\frac{{B}_2}{\mathcal{C}_2}\left( \alpha_2^{av} - {\alpha_2}^{0} \right),
\end{equation}
or
\begin{equation}\label{eq:Tav_eq2}
T^{av} = T - \frac{{B}_2}{\mathcal{C}_2}\left( \alpha_2^{av} - {\alpha_2}^{0} \right),
\end{equation}

Combination of \cref{eq:Tav_eq1,eq:Tav_eq2} leads to 
$\gamma_1 = \gamma_2$, which contradicts the initial condition (\cref{eq:RiemannInitial}). 
Thus, the temperature relaxation procedure does not violate the temperature equilibrium.
Moreover, since velocity and temperature are spatially uniform, the diffusion processes (heat conduction and viscosity) does not have any impact on the solution.
}
\end{proof}

\subsection{Extension to high order and interface sharpening}
\label{subsec:high-order}
The scheme can be extended to higher orders with the MUSCL or WENO scheme. Moreover, to minimize numerical diffusion of material interfaces, we apply special interface-sharpening schemes  \cite{ZHANG2019124610,Shyue2014An,Chiapolino2018,Chiapolino2017,Ch2019Interface}. One principle for implementing these schemes is that the reconstruction schemes for volume fractions and phase densities should be consistent, otherwise, the PVT property is violated. We give a simple explanation on this issue below.

Observing the above proof of the PVT property, one can see that an important condition for proving the temperature equilibrium is $\rho_{k}^{*} = \rho_k$. The high-order extensions should also ensure this condition. This is deduced from
\begin{subequations}
\begin{align}
  \left( \alpha_k \rho_k \right)_{i}^{*} = \xi \left( \alpha_k \rho_k \right)_{i-1} + \left( 1 - \xi \right) \left( \alpha_k \rho_k \right)_{i} \label{eq:rhostar1} \\
  \left( \alpha_k \right)_{i}^{*} = \xi \left( \alpha_k \right)_{i-1} + \left( 1 - \xi \right) \left( \alpha_k \right)_{i}. \label{eq:rhostar2}
\end{align}
\end{subequations}

The corresponding high-order formulation is given as follows:
\begin{subequations}
\begin{align}
\left( \alpha_k \rho_k \right)_{i}^{*} = \xi \left( \alpha_k \rho_k \right)_{i-1,RF} + \left( 1 - \xi \right) \left( \alpha_k \rho_k \right)_{i,RF}, \label{eq:rhostar_high1} \\
\left( \alpha_k \right)_{i}^{*} = \xi \left( \alpha_k \right)_{i-1,RF} + \left( 1 - \xi \right) \left( \alpha_k \right)_{i,RF}, \label{eq:rhostar_high2}
\end{align}
\end{subequations}
where the subscript ``RF'' represents the reconstructed values on the right face of the current cell.

Assume that we use a reconstruction scheme that is a homogeneous function of degree 1 with respect to the reconstruction stencil, i.e., $\text{Rec}\left( \beta \mathcal{P} \right) = \beta \text{Rec}\left( \mathcal{P} \right), \; \beta = \text{const}>0$ . Note that the TVD schemes are such functions.  Then the reconstructed values are as follows:
\begin{subequations}
\begin{align}
\left( \alpha_k \right)_{i,RF} = \text{Rec} \left( \mathcal{P} \right), \label{eq:rec_scheme1} \\
\left( \alpha_k \rho_k \right)_{i,RF} = \overline{\text{Rec}} \left( \rho_k {\mathcal{P}} \right) = \rho_k \overline{\text{Rec}} \left(  {\mathcal{P}} \right), \label{eq:rec_scheme2}
\end{align}
\end{subequations}
where $\mathcal{P}$ is the reconstruction stencil, for example, for the MINMOD scheme, $\mathcal{P} = \{ \left( \alpha_k \right)_{i-1}, \; \left( \alpha_k \right)_{i}, \; \left( \alpha_k \right)_{i+1} \}$.  $\text{Rec}$ and $\overline{\text{Rec}}$ are the reconstruction scheme for $\alpha_k$ and $\alpha_k \rho_k$, respectively.

Combination of \cref{eq:rhostar_high1,eq:rhostar_high2,eq:rec_scheme1,eq:rec_scheme2} leads to the conclusion that $\rho_{k}^{*} = \rho_k$ holds only when $\text{Rec} \left( \mathcal{P} \right) = \overline{\text{Rec}} \left( {\mathcal{P}} \right)$. This means that the same  scheme should be used for reconstruction of $\left( \alpha_k \right)_{i,RF}$ and $\left( \alpha_k \rho_k \right)_{i,RF}$; Otherwise, the temperature equilibrium is violated. The numerical results presented in \Cref{sec:numer_res} also confirm this fact.

\section{Numerical results}
\label{sec:numer_res}

In this section we perform several numerical tests with the purpose to verify our model and numerical methods and also compare with some other methods presented in literature. In the laser ablation problem (\cref{subsec:ablation}), the variables are measured in the centimetre-gram-microsecond system of units, and in SI units for other tests.

\subsection{Preservation of the PVT property}
The purpose of this test is to check the capability of different models to keep the PVT property.
We consider the translation of material interface with initially uniform velocity $u = 1.00 \times 10^{3}$m/s, pressure $p = 1.00 \times 10^{5}$Pa and temperature $T = 300.00$K. Pressures and temperatures are all in equilibrium.  The computational domain is $[0\text{m},1\text{m}]$, the material interface is initially located at $x = 0.20\text{m}$. The EOS parameters for the left component $\gamma=4.40, \; C_v = 500.00{\text{J}/(\text{kg}\cdot \text{K})}, \; p_{\infty} = 6.00 \times 10^{8}{\text{Pa}}$, and those for the right component $\gamma=1.40, \; C_v = 200.00{\text{J}/(\text{kg}\cdot \text{K})}, \; p_{\infty} = 0.00{\text{Pa}}$. 


Here and in the following subsections, we evaluate four different schemes:
\begin{enumerate}
\item[(a)] The conservative four-equation model with one temperature (4-eqn model 1T.)  \cite{Lemartelot2014},
\item[(b)] The five-equation model with one temperature (5-eqn model 1T.) \cite{Alahyari2015},
\item[(c)] The six-equation model with two temperatures (6-eqn model 2T.)  \cite{saurel2009simple},
\item[(d)] The proposed six-equation model with two temperatures and thermal relaxation (6-eqn model 2T. relaxed).
\end{enumerate}

\begin{figure}[htb]
\centering
\subfloat[Density]{\label{fig:fig00a}\includegraphics[width=0.5\textwidth]{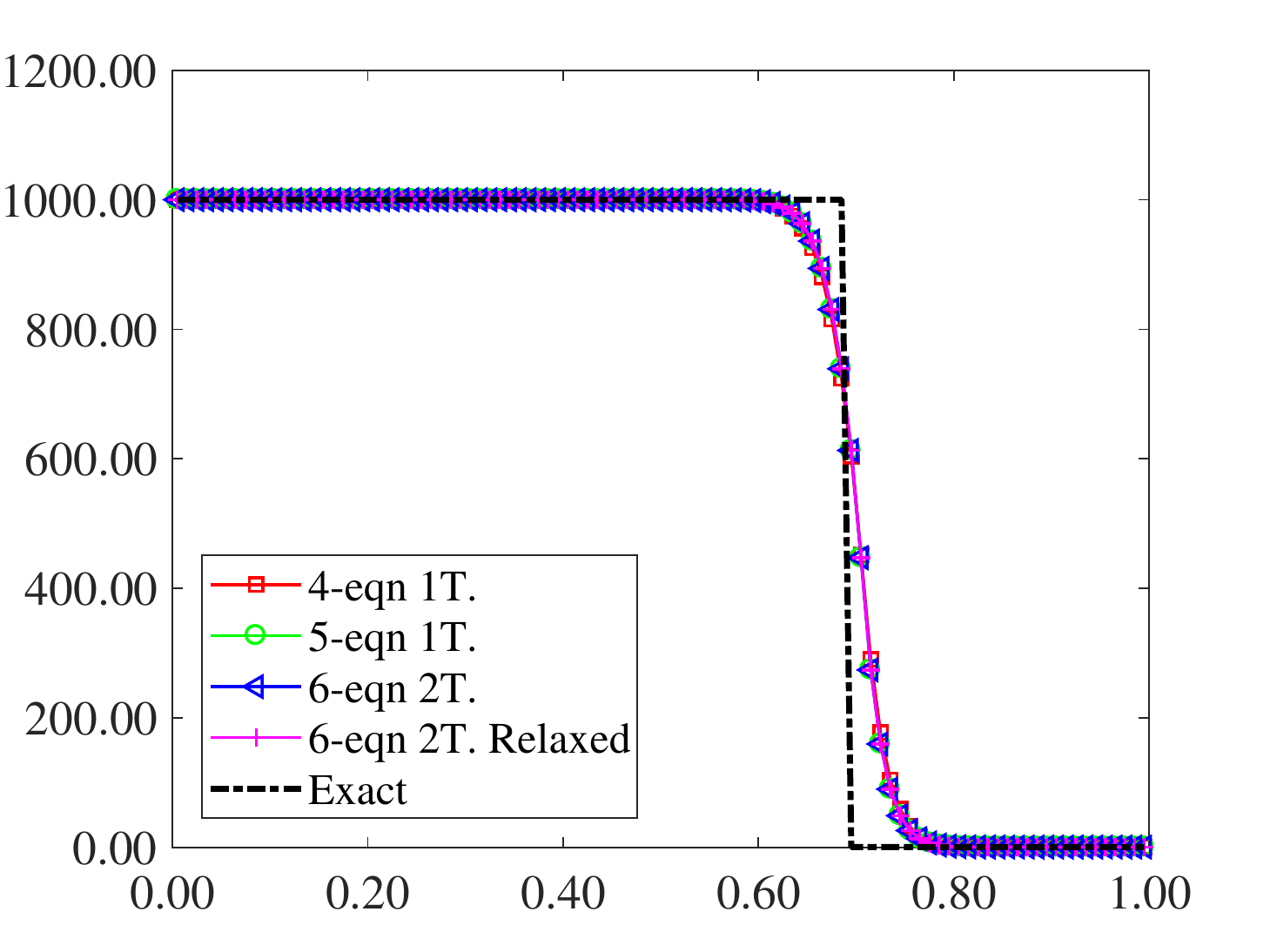}}
\subfloat[Volume/Mass fraction]{\label{fig:fig00b}\includegraphics[width=0.5\textwidth]{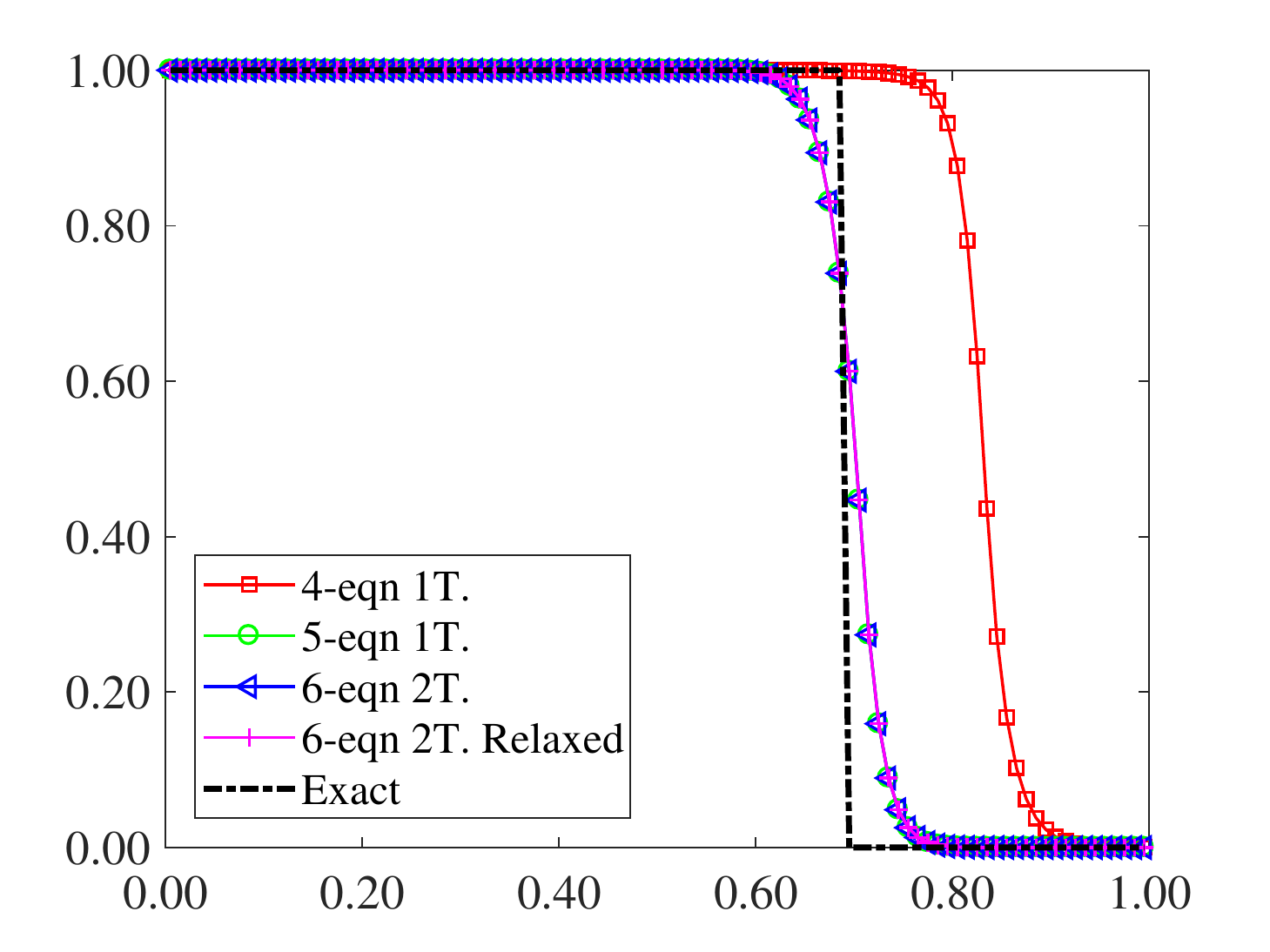}}\\
\subfloat[Velocity]{\label{fig:fig00c}\includegraphics[width=0.5\textwidth]{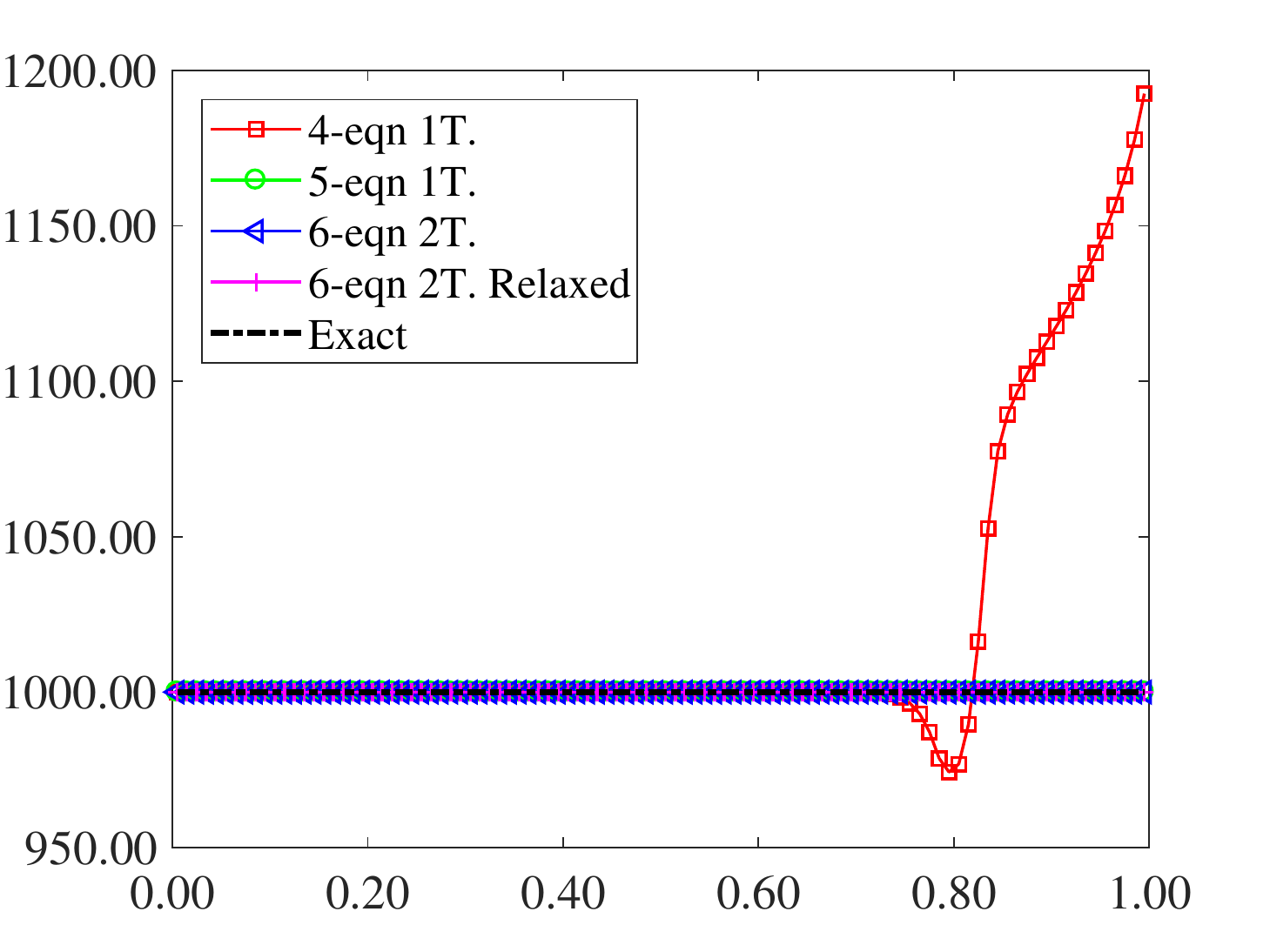}} 
\subfloat[Temperature]{\label{fig:fig00d}\includegraphics[width=0.5\textwidth]{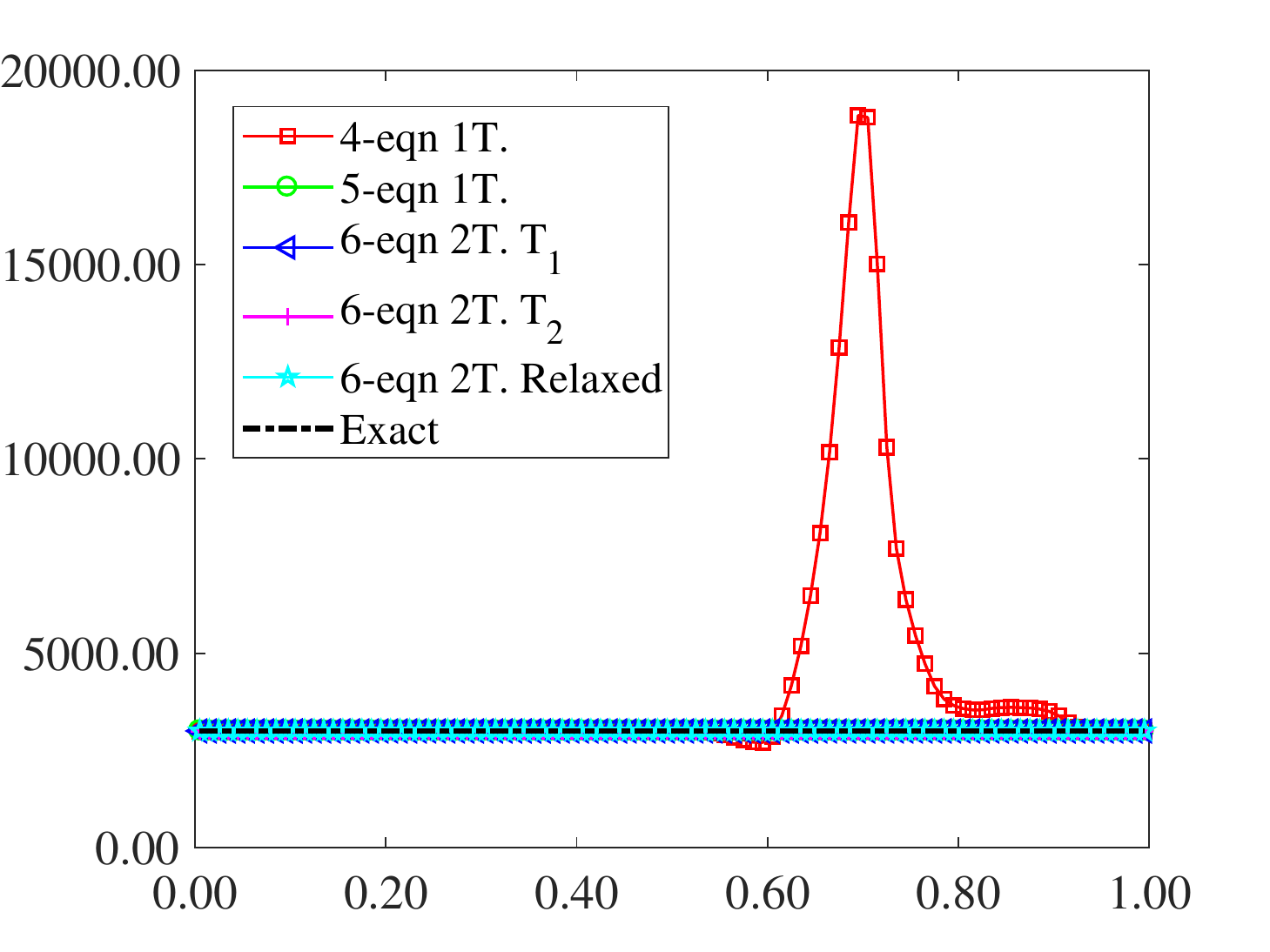}}\\
\subfloat[Pressure]{\label{fig:fig00e}\includegraphics[width=0.5\textwidth]{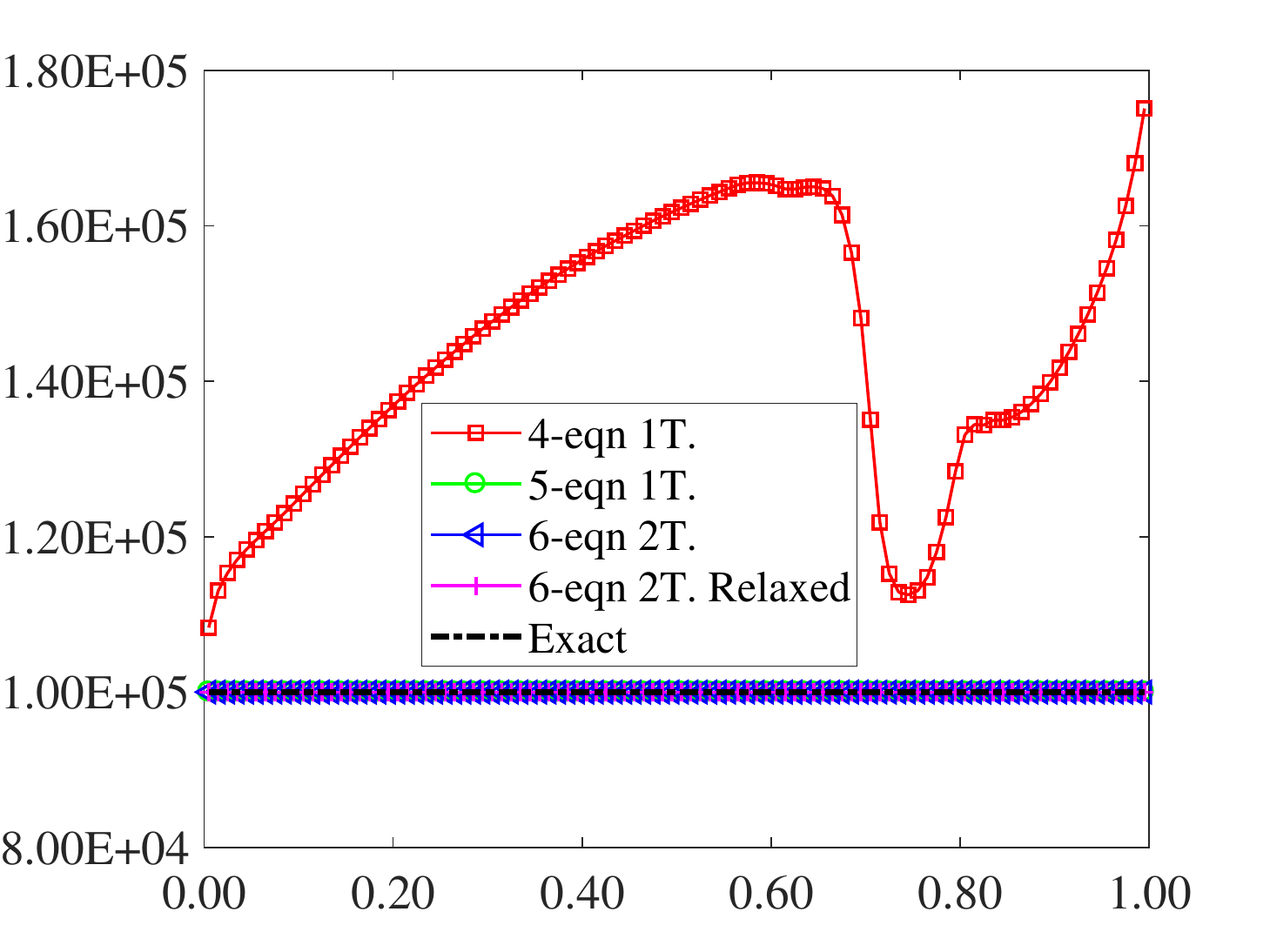}} 
\subfloat[Temperature]{\label{fig:fig00f}\includegraphics[width=0.5\textwidth]{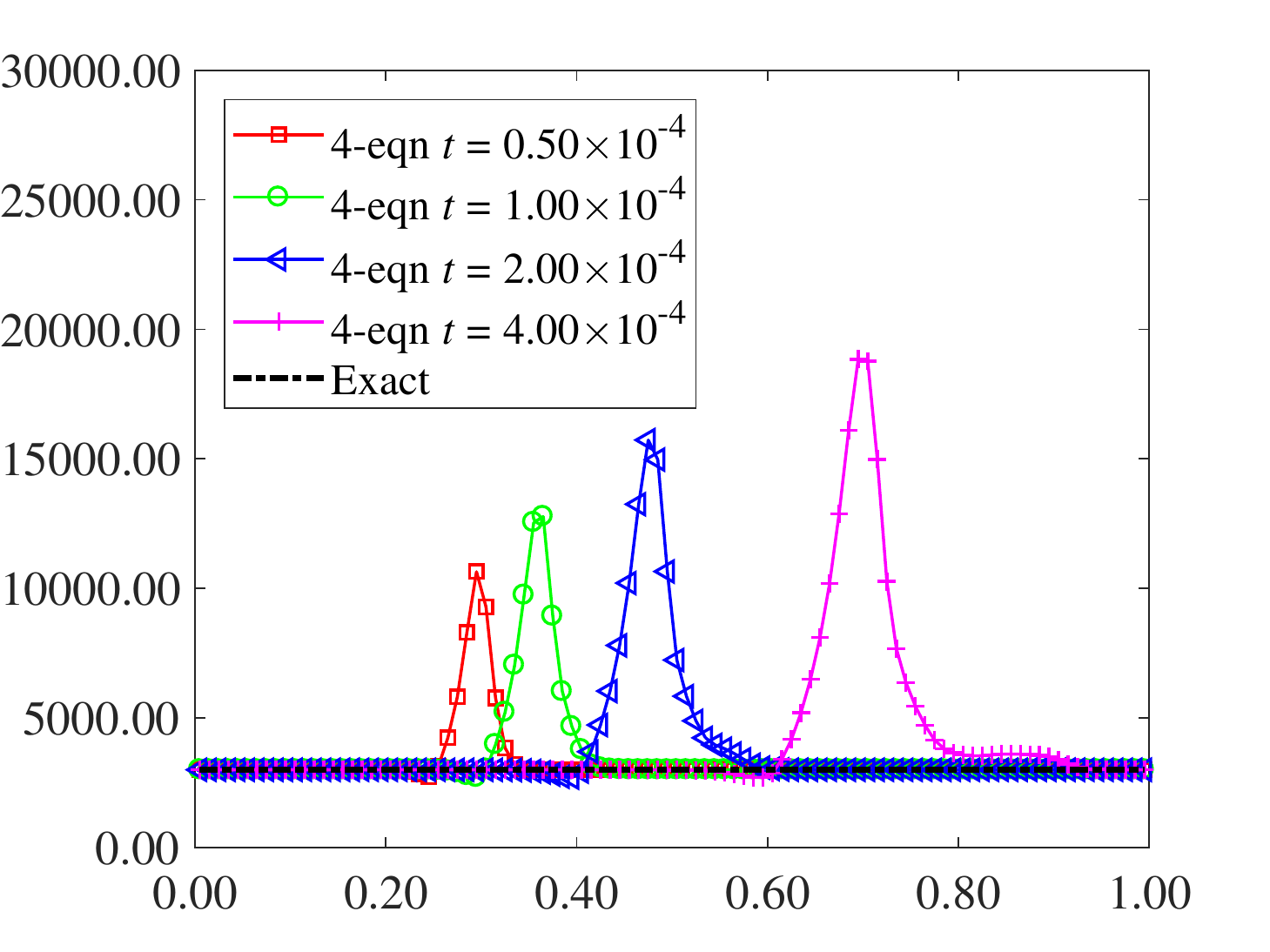}}
\caption{Pure translation of a two-fluid system: numerical solutions for different flow parameters.}
\label{fig:fig00} 
\end{figure}

\begin{figure}[htb]
\centering
\subfloat[Temperature]{\label{fig:fig0a}\includegraphics[width=0.5\textwidth]{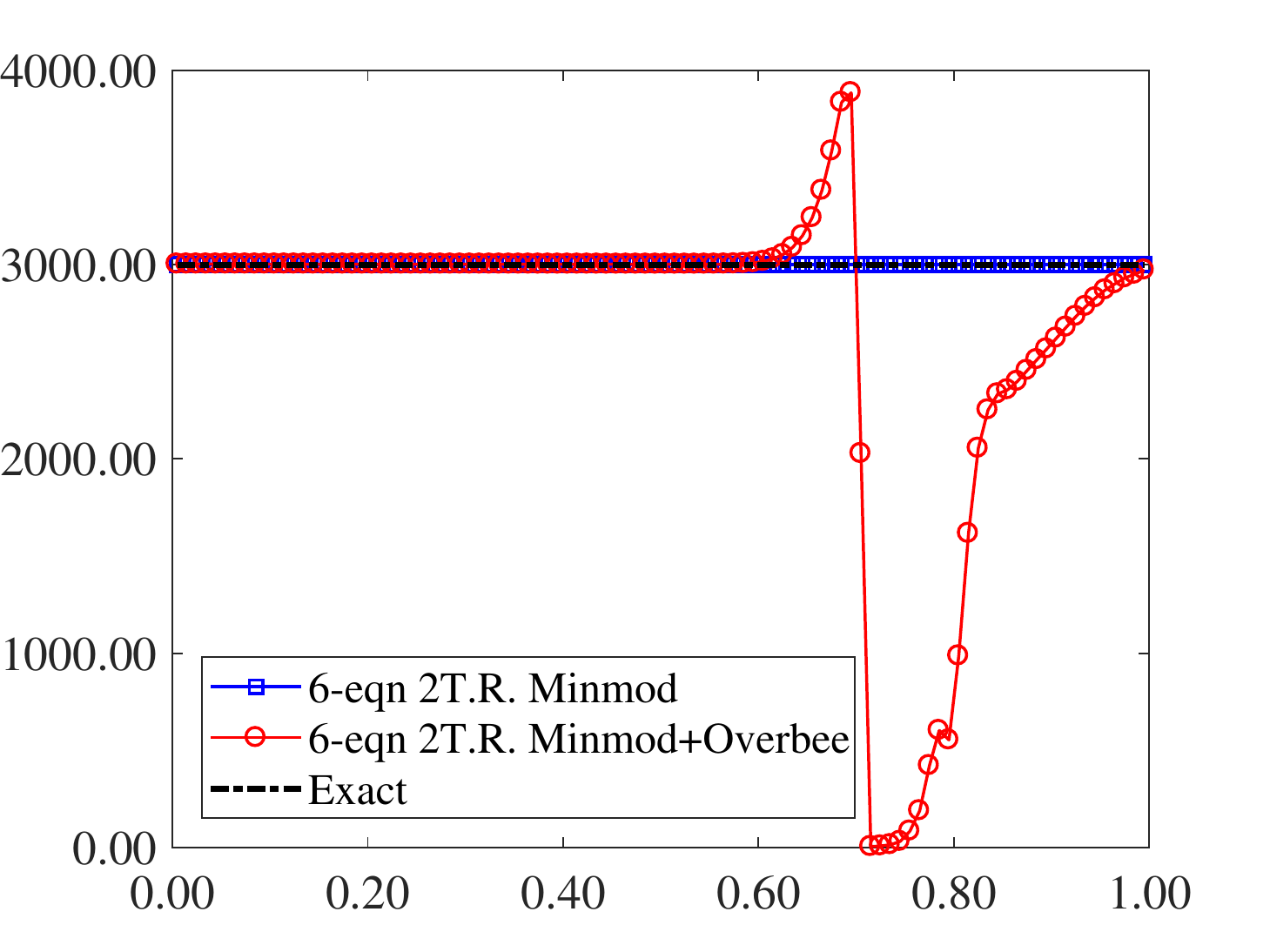}} 
\subfloat[Pressure]{\label{fig:fig0b}\includegraphics[width=0.5\textwidth]{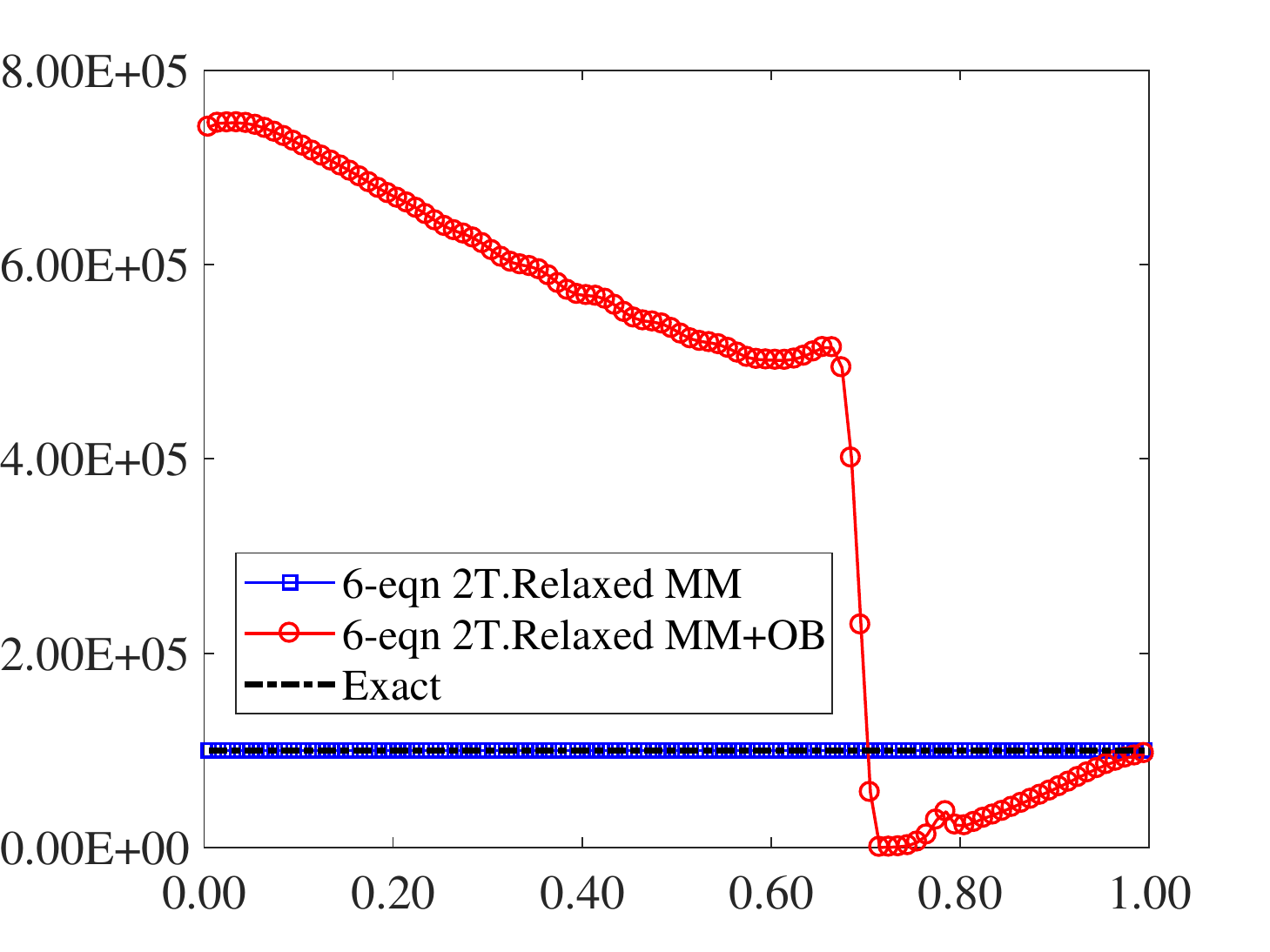}} 
\caption{Numerical results of temperature(left) and pressure(right) along the computational domain when different reconstruction schemes used.}
\label{fig:fig0} 
\end{figure}

We perform computations with the above models to the moment $t = 5.00 \times 10^{-4}$s. The numerical results are illustrated in \Cref{fig:fig00}.
The numerical results with different reconstruction schemes are displayed in \Cref{fig:fig0}.

It can be seen that the five-equation and the six-equation models maintain the PVT property very well, while the four-equation model triggers spurious oscillations in the velocity, pressure and temperature profiles (see \Cref{fig:fig00c,fig:fig00d,fig:fig00e}). From \Cref{fig:fig00f} we see that the non-physical temperature spikes in the numerical results obtained with the conservative four-equation model tend to increase with time. This error can infect the solution in the computational domain through the heat conduction.

Moreover, as mentioned in \cref{subsec:high-order} regarding high-order extensions, the reconstruction schemes for the volume fraction and phase densities should be consistent, otherwise non-physical spikes in temperature arise. If we use the interface-sharpening scheme Overbee \cite{Chiapolino2017} for reconstructing the volume fraction, and MINMOD for the phase densities, we obtain the results shown in \Cref{fig:fig0}. The defect  appeared in this figure arises because the inconsistent scheme combinations fail to maintain constant phase densities.

\subsection{Shock tube problem with heat conduction}
\label{subsec:shocktubetest}

In this section we consider a two-fluid shock tube problem with the purpose of evaluating different models. Two fluids are initially at rest and and separated by  the material interface  located at $x=0.7$m separating them. The fluid on the left has the EOS parameters as $\gamma=4.40, \; p_{\infty} = 6.00 \times 10^{8}{\text{Pa}}, \; C_v = 1606.00{{\text{J}/(\text{kg}\cdot \text{K})}}$, and and that on the right -- $\gamma=1.40, \; p_{\infty} = 0.00{\text{Pa}}, \; C_v = 714.00{{\text{J}/(\text{kg}\cdot \text{K})}}$. The initial pressure and temperature on both sides are given as follows:
\[ 0.00<x<0.70{\text{m}}:\quad  p=1.00 \times 10^{9}{\text{Pa}}, \; T = 293.02{\text{K}}, \]
\[ 0.70<x<1.00{\text{m}}:\quad  p=1.00 \times 10^{5}{\text{Pa}}, \; T = 7.02{\text{K}}. \]

The initial densities are determined from the corresponding EOS.

\paragraph{Test without heat conduction}

\begin{figure}[htb]
\centering
\subfloat[Density]{\label{fig:fig1a}\includegraphics[width=0.5\textwidth]{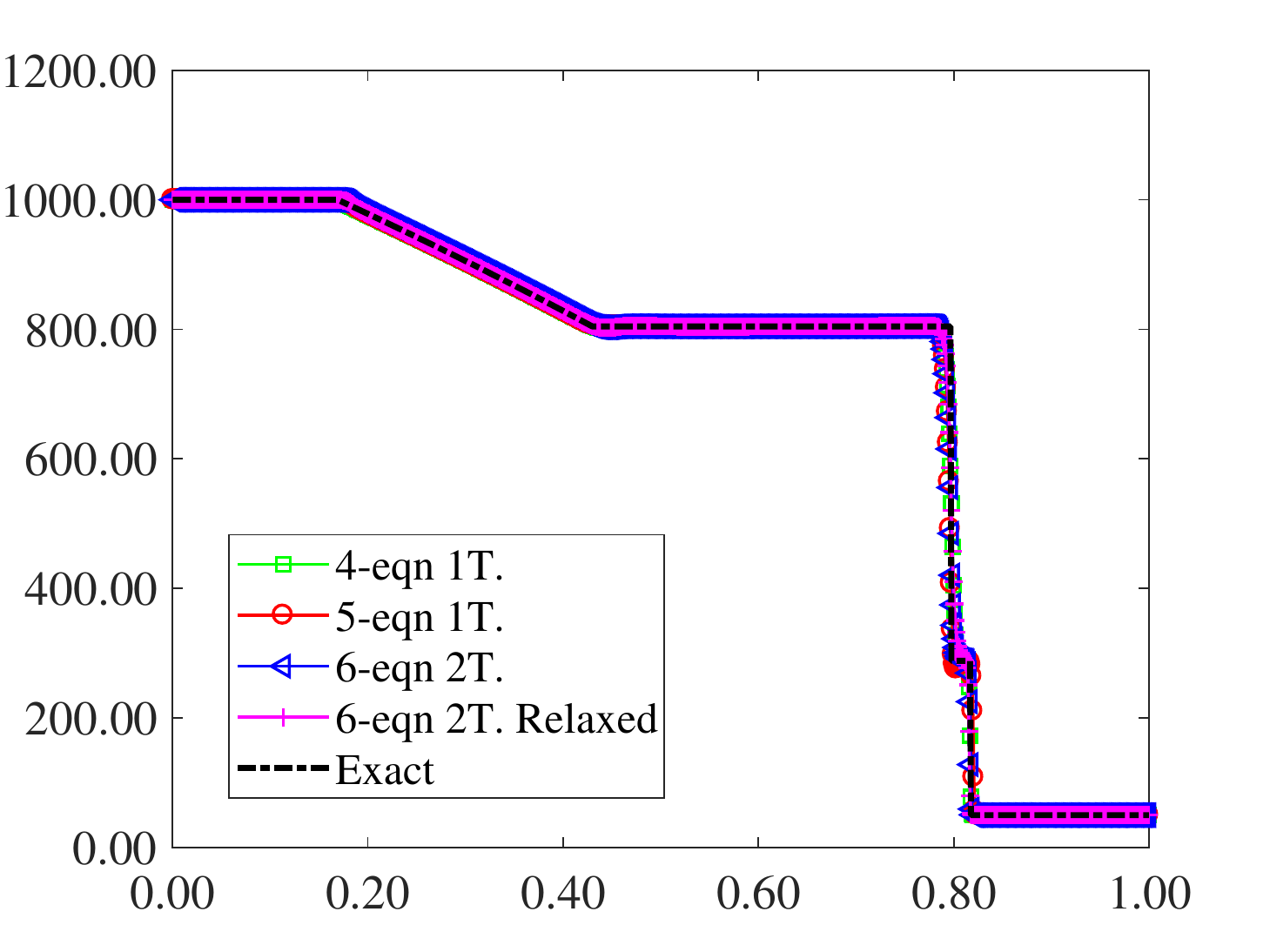}}
\subfloat[Density, locally enlarged]{\label{fig:fig1b}\includegraphics[width=0.5\textwidth]{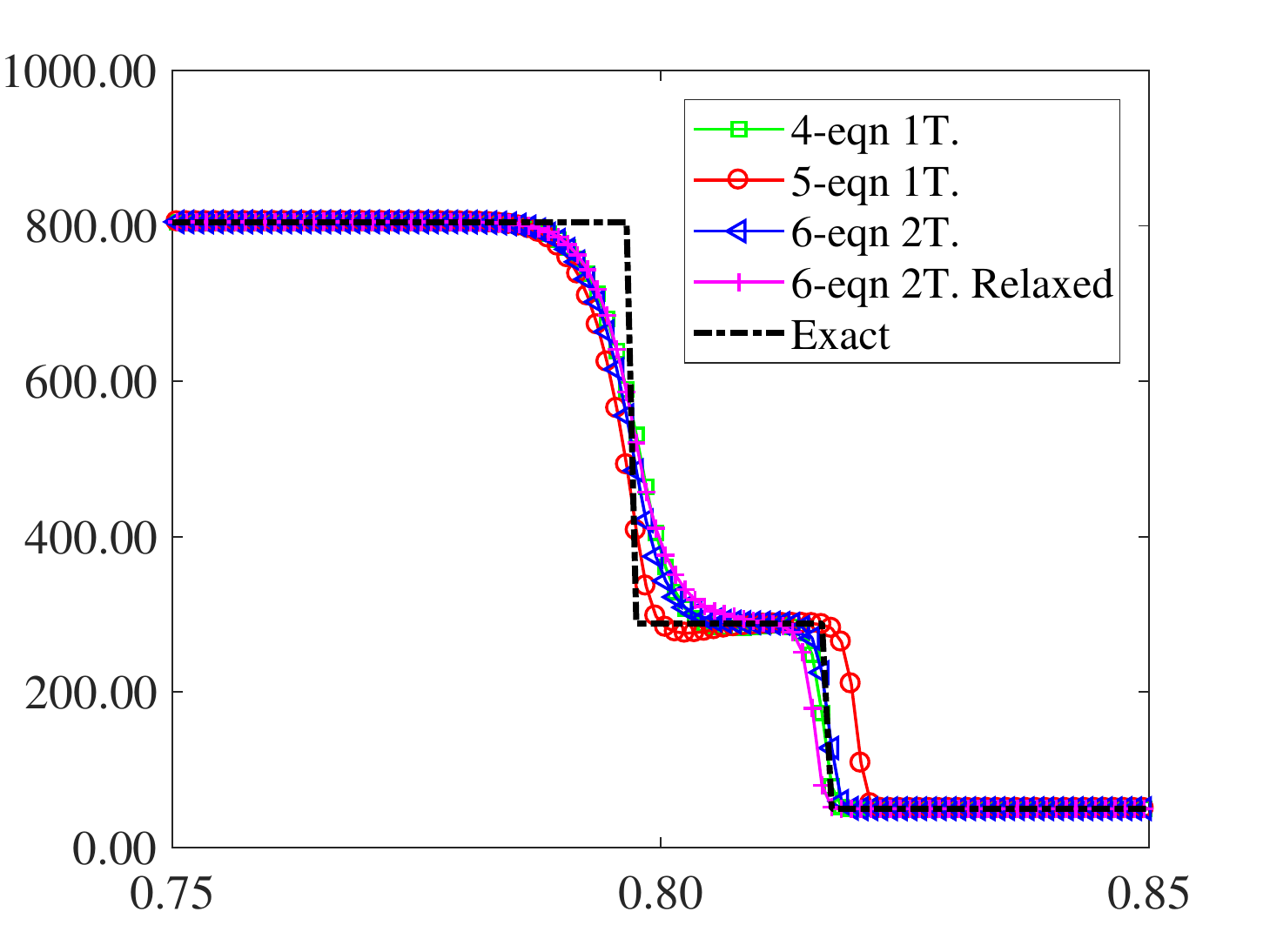}} \\
\subfloat[Temperature]{\label{fig:fig1c}\includegraphics[width=0.5\textwidth]{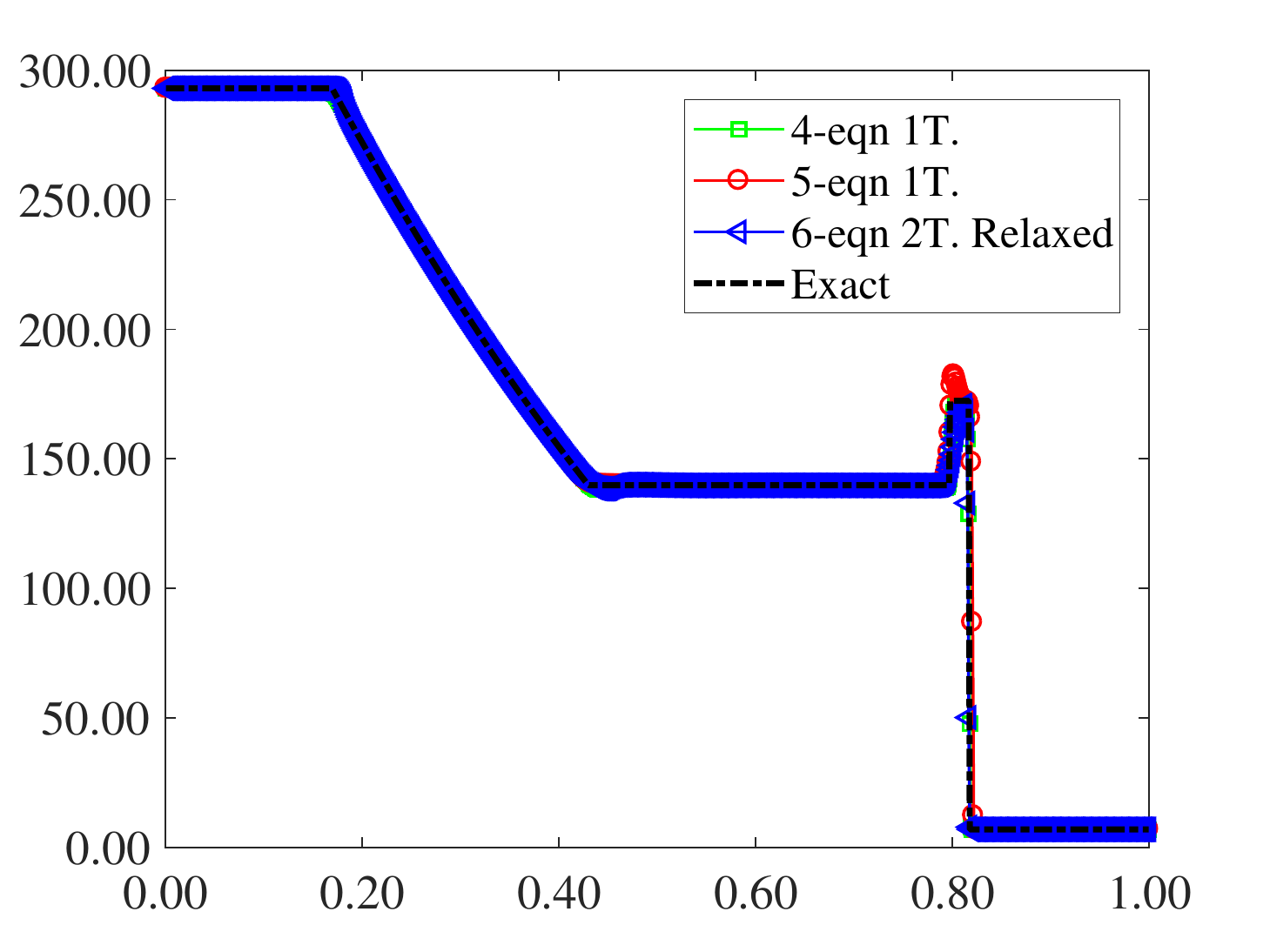}}
\subfloat[Temperature, locally enlarged]{\label{fig:fig1d}\includegraphics[width=0.5\textwidth]{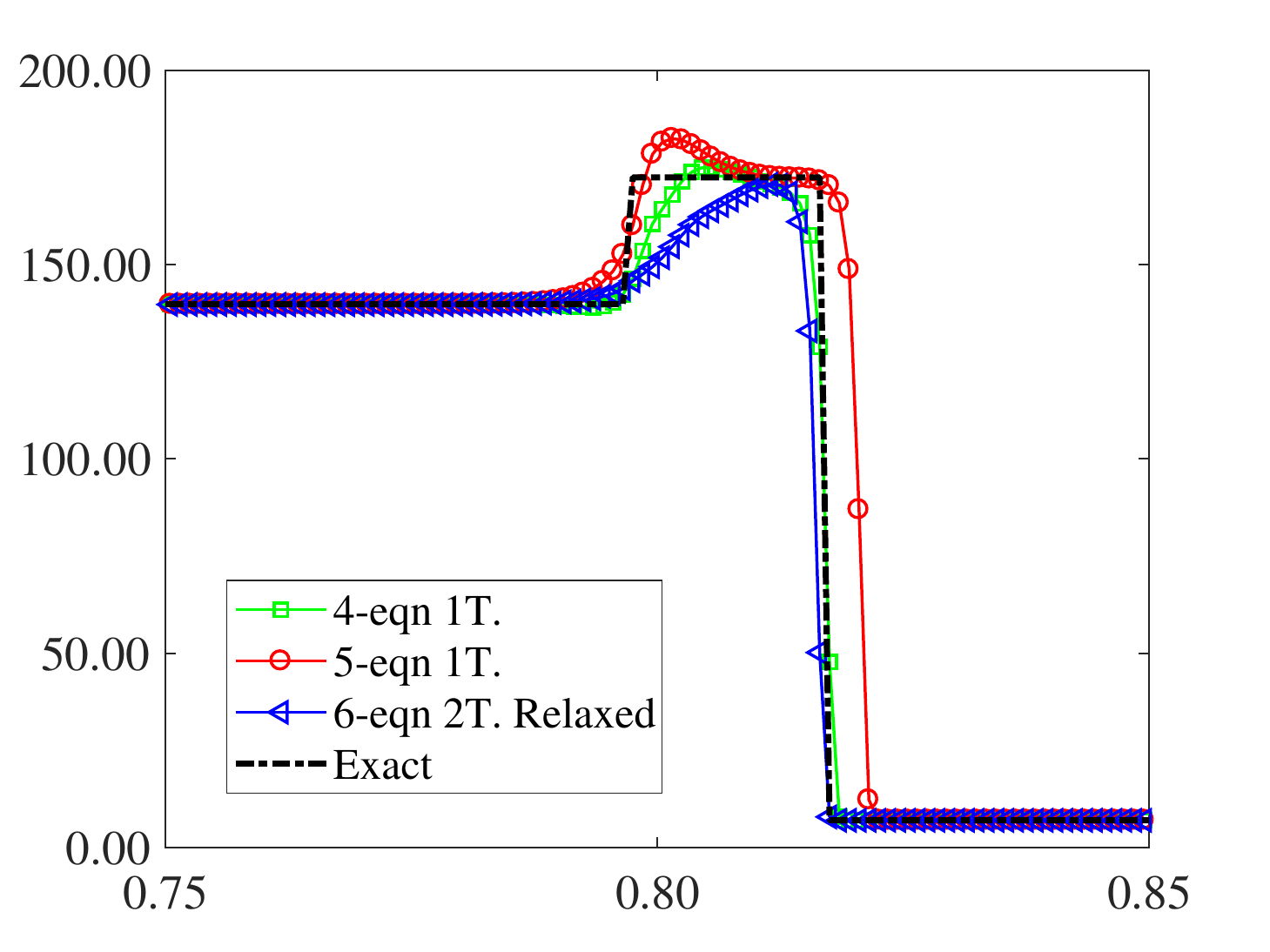}} \\
\subfloat[Temperatures obtained with 2T. model]{\label{fig:fig1e}\includegraphics[width=0.5\textwidth]{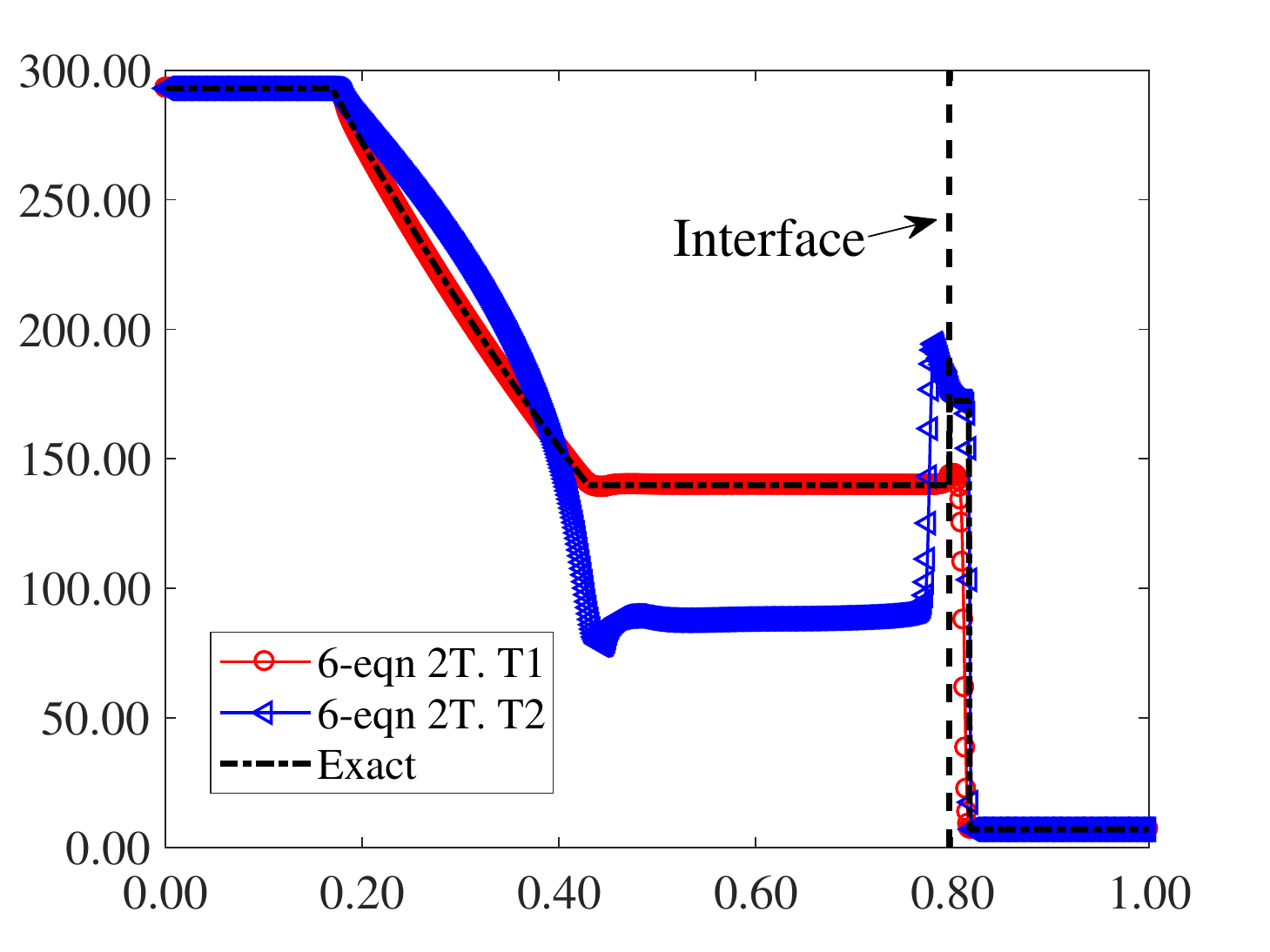}}
\caption{Numerical results for the two-fluid shock tube problem without heat conduction.}
\label{fig:fig1} 
\end{figure}

Computations are performed on a 1000-cell uniform grid. 
The obtained numerical results obtained at the time moment $t = 2.00 \times 10^{-4}$s  are compared to the exact Riemann solution in \Cref{fig:fig1}. 
The exact solution consists of a leftward rarefaction wave, a rightward contact wave (interface) and a rightward shock wave.

From the density profiles \Cref{fig:fig1b} and temperature profiles \Cref{fig:fig1d}, one can see that the shock wave velocity in the five-equation model with one temperature appears to be overestimated. This stems from the different estimation of mixture acoustic velocity inside the diffused zone. 

Note that as a solution to the Euler equations, the exact Riemann solution does not include any thermal relaxation. Therefore, the solution to the two-temperature six-equation model with no thermal relaxation is expected to better match the exact solution than that with thermal relaxation.
In \Cref{fig:fig1e} we plot the two temperatures of fluids calculated in the non-equilibrium model without temperature relaxation. As seen, the temperature of the first fluid quite well matches the exact solution on the left of the interface, while the temperature of the second fluid similarly does on the right. 
 Thermal relaxation drives the two temperatures into an equilibrium temperature -- the profile denoted as ``6-eqn model 2T. relaxed''   in \Cref{fig:fig1c,fig:fig1d}.

\begin{figure}[htb]
\centering
\subfloat[Density]{\label{fig:fig2a}\includegraphics[width=0.5\textwidth]{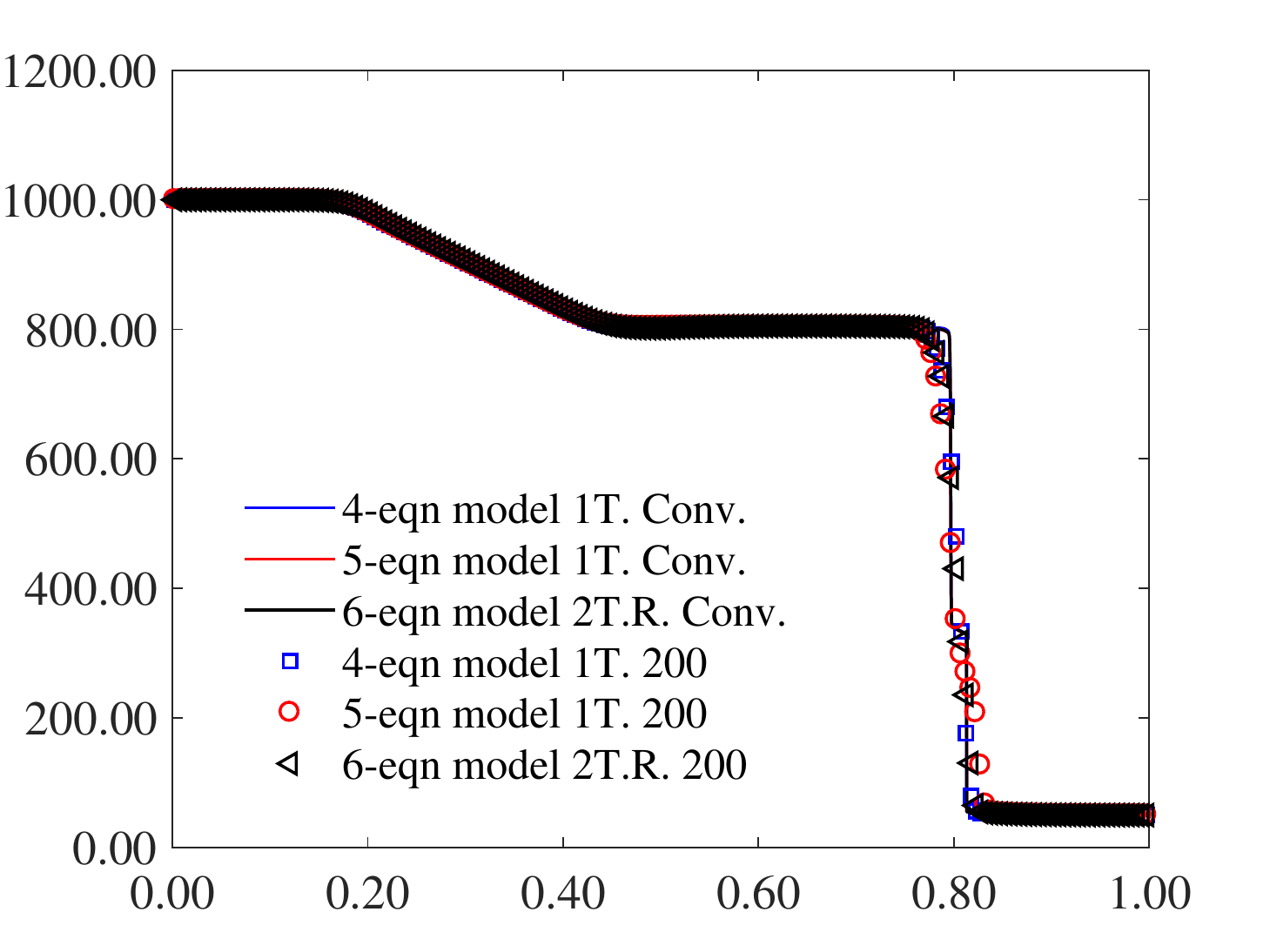}}
\subfloat[Density, locally enlarged]{\label{fig:fig2b}\includegraphics[width=0.5\textwidth]{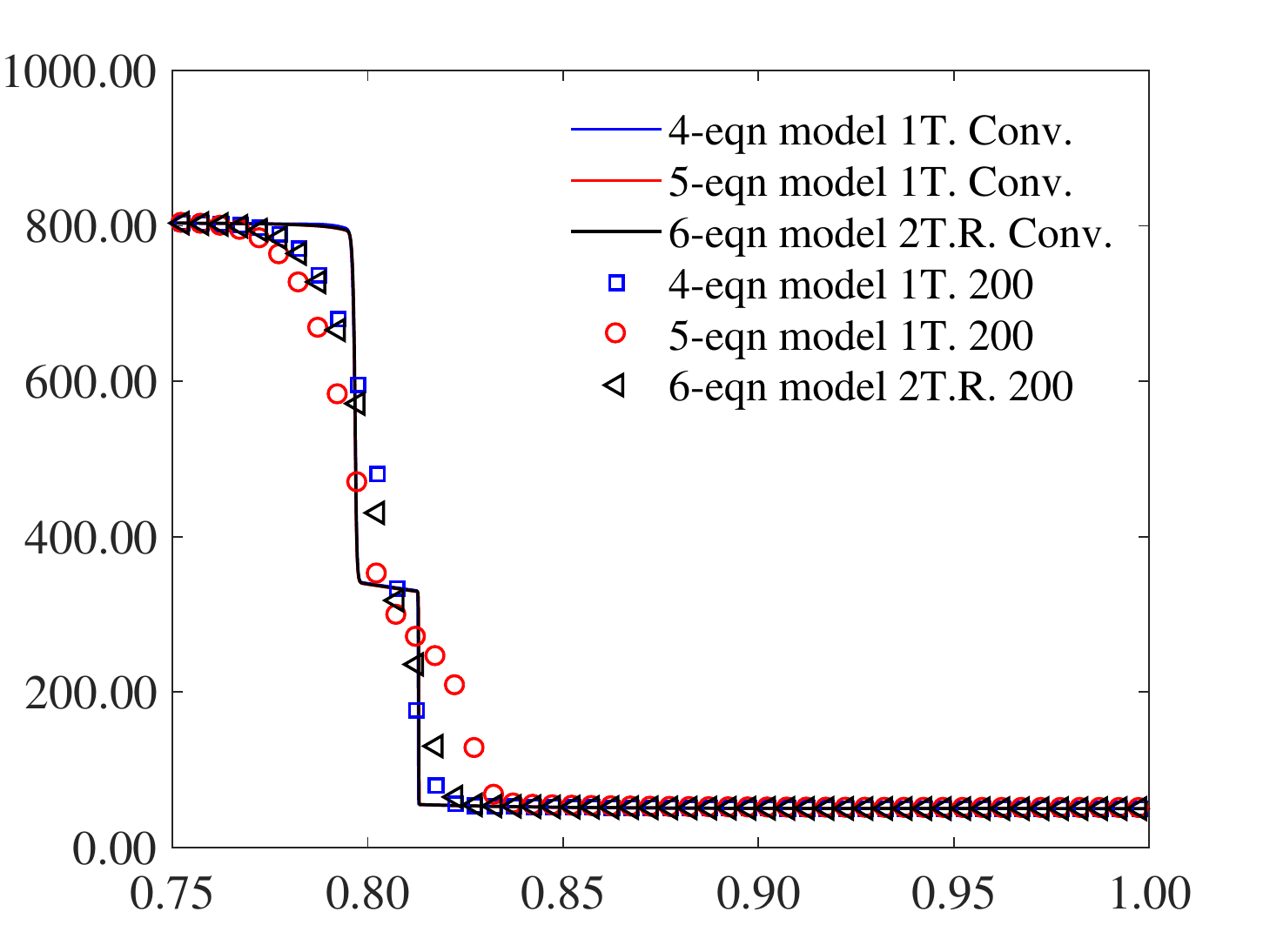}} \\
\subfloat[Temperature]{\label{fig:fig2c}\includegraphics[width=0.5\textwidth]{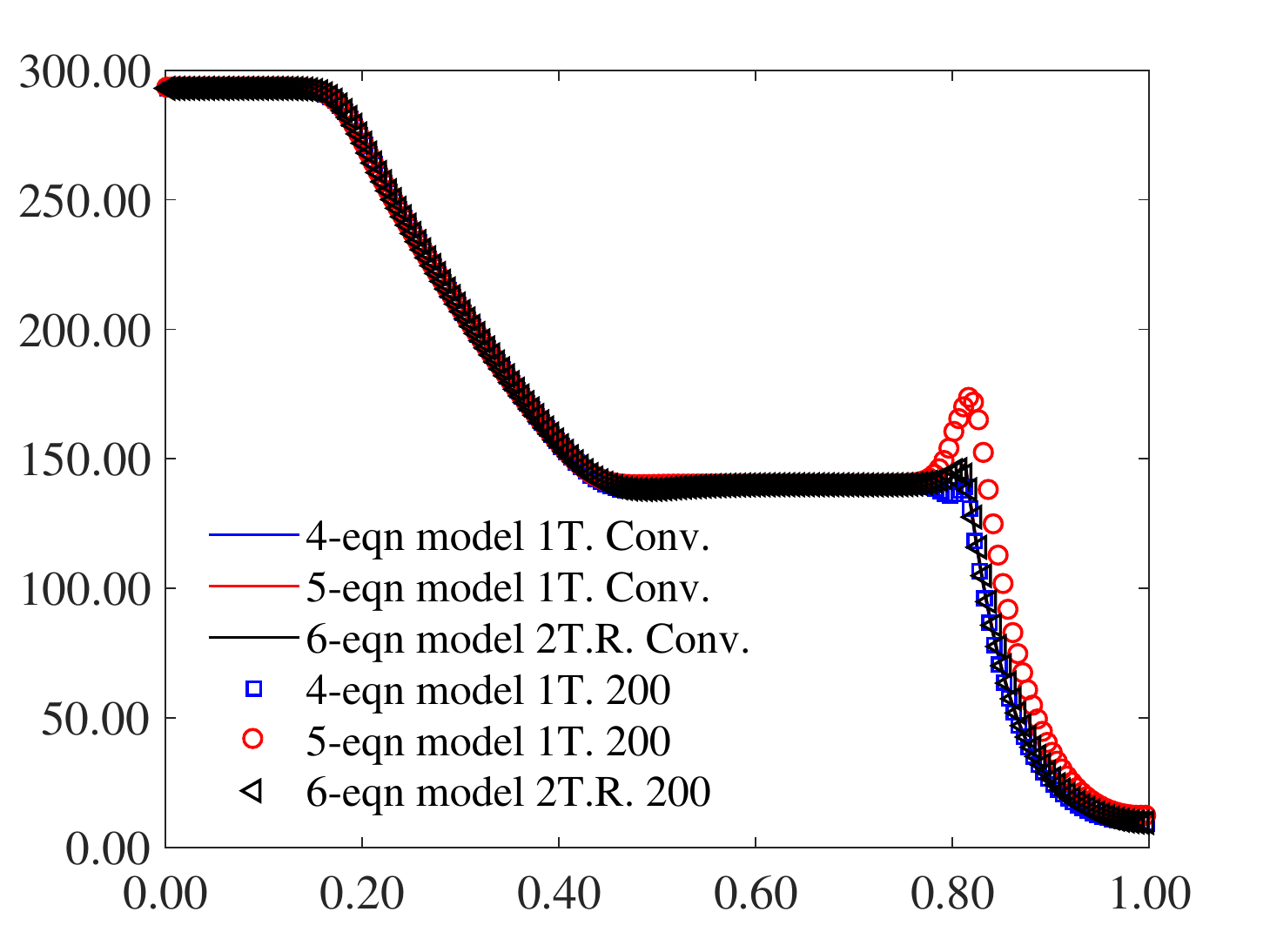}}
\subfloat[Temperature, locally enlarged]{\label{fig:fig2d}\includegraphics[width=0.5\textwidth]{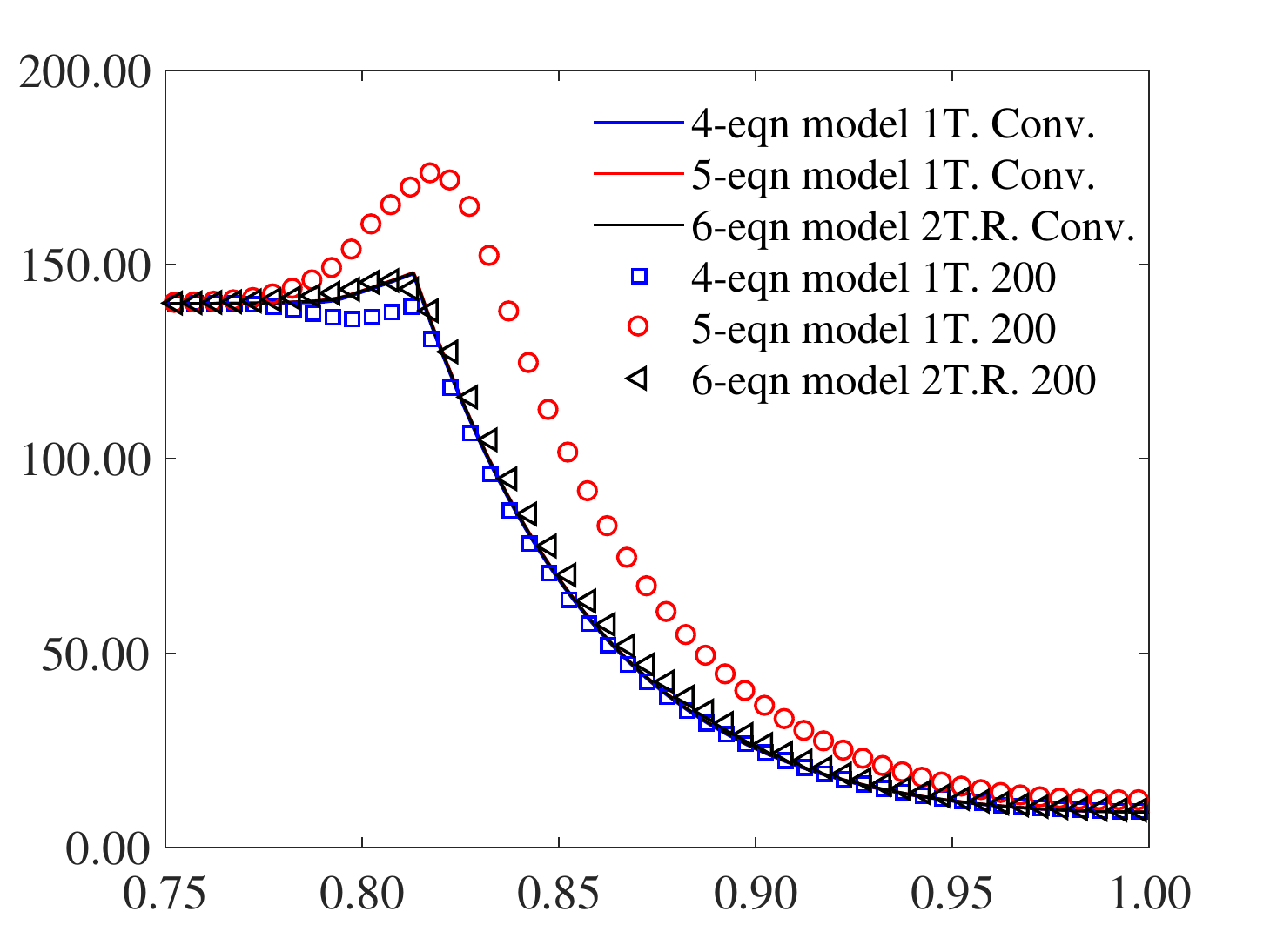}} \\
\subfloat[Velocity]{\label{fig:fig2e}\includegraphics[width=0.5\textwidth]{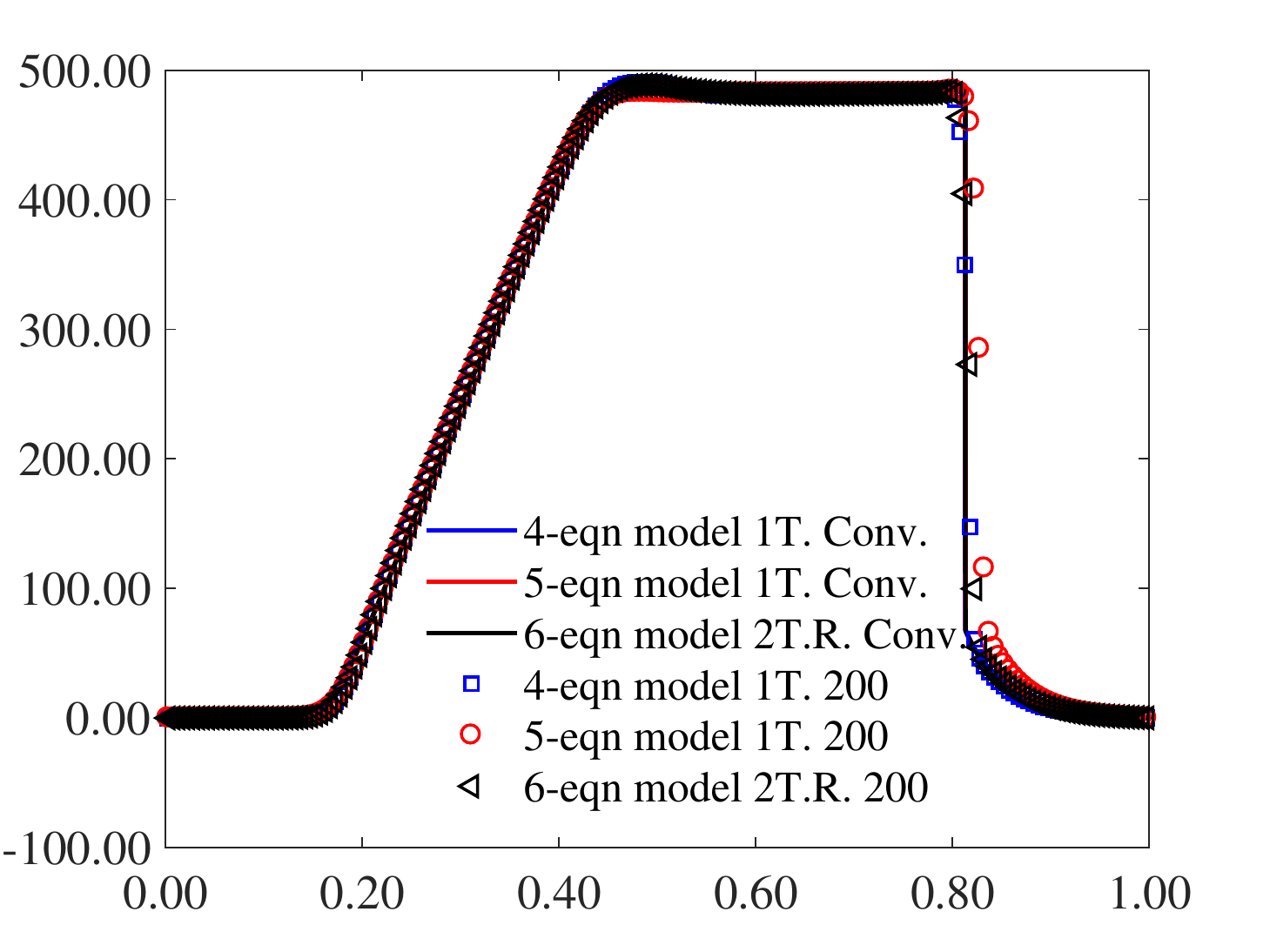}}
\subfloat[Velocity, locally enlarged]{\label{fig:fig2f}\includegraphics[width=0.5\textwidth]{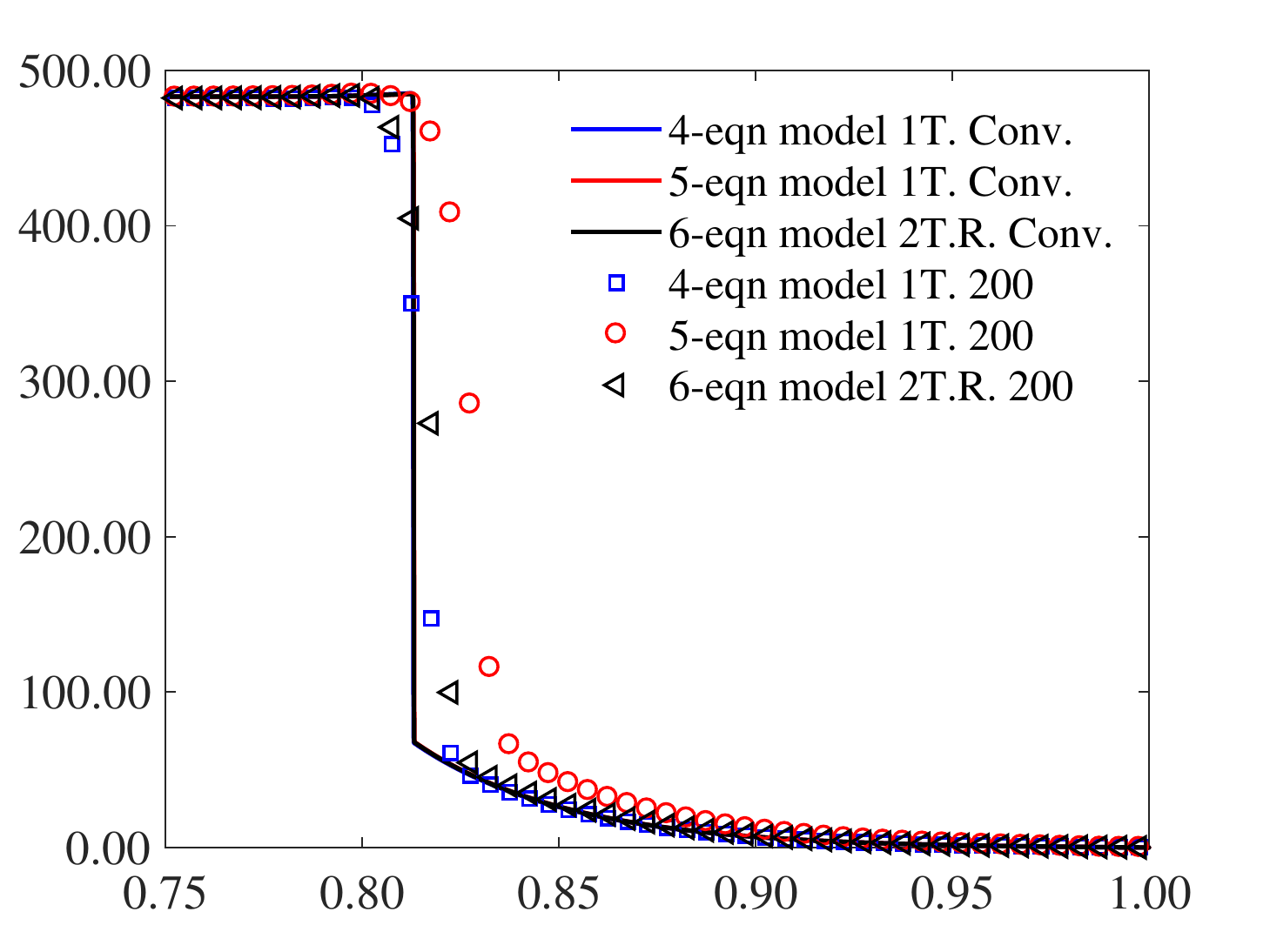}} 
\caption{Numerical results for the two-fluid shock tube problem with  equal phase  heat conductivity  $\lambda = 1.00 \times 10^6$.}
\label{fig:fig2} 
\end{figure}

\begin{figure}[htb]
\centering
\subfloat[Density]{\label{fig:fig3a}\includegraphics[width=0.5\textwidth]{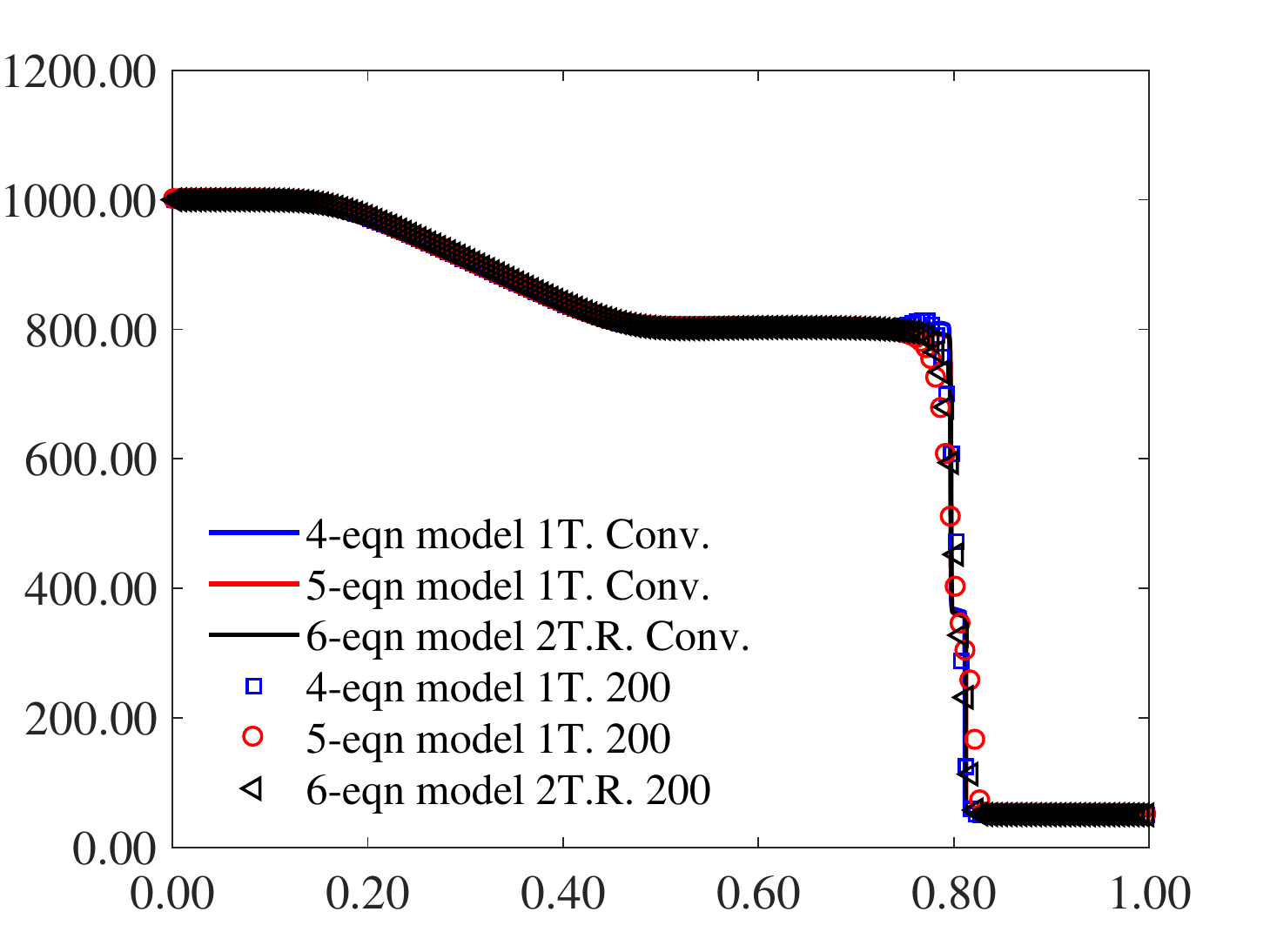}}
\subfloat[Density, locally enlarged]{\label{fig:fig3b}\includegraphics[width=0.5\textwidth]{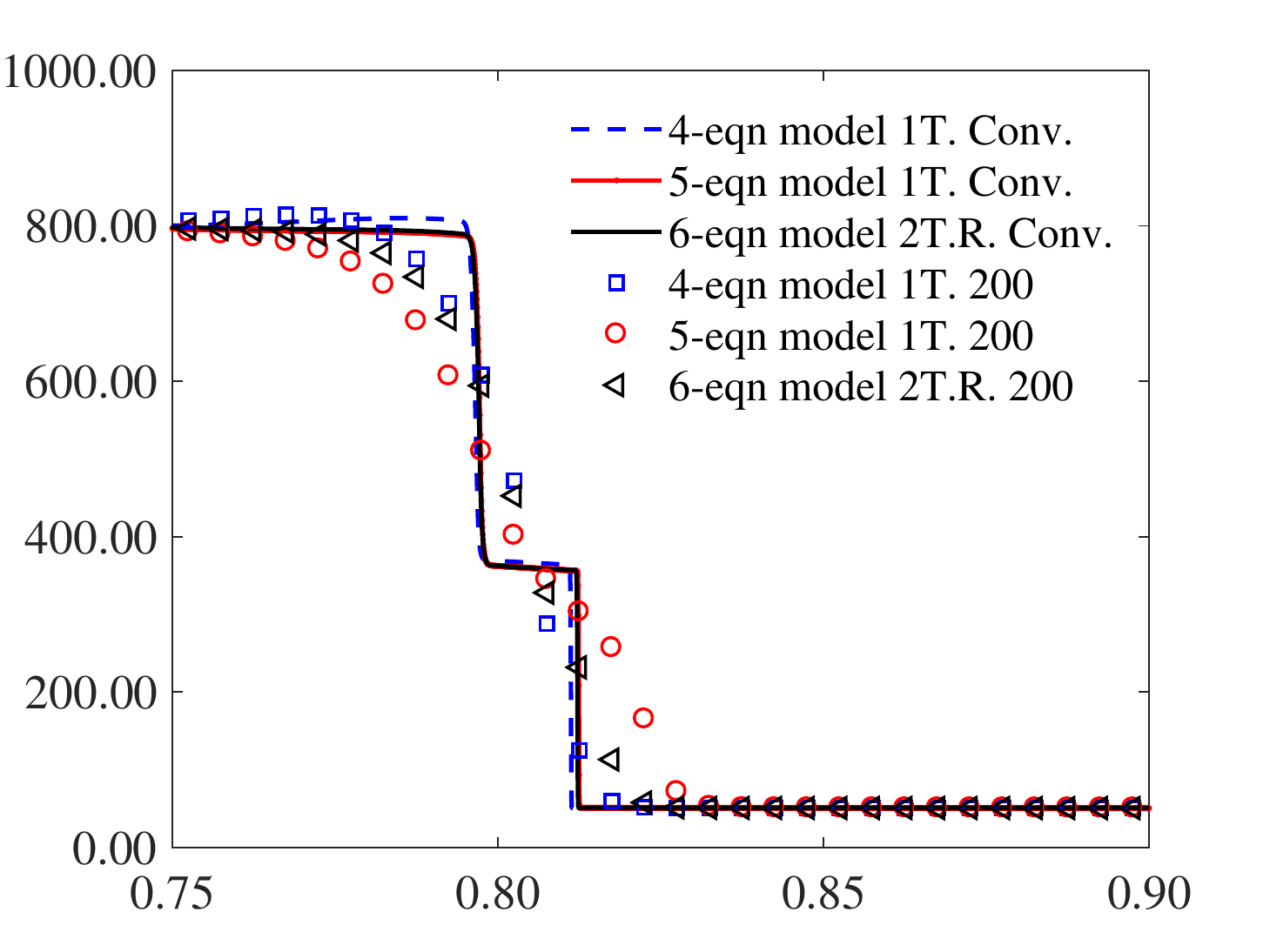}} \\
\subfloat[Temperature]{\label{fig:fig3c}\includegraphics[width=0.5\textwidth]{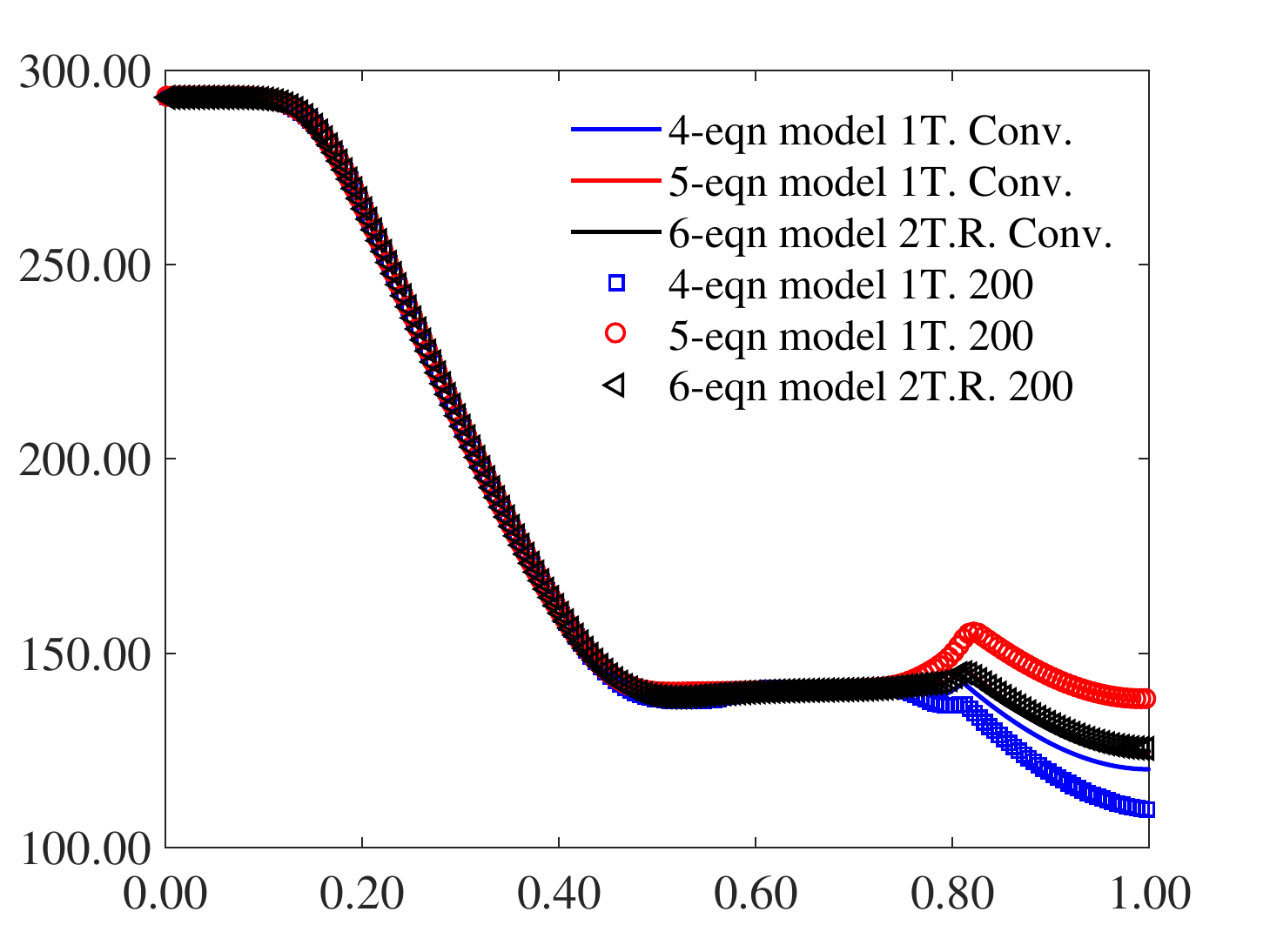}}
\subfloat[Temperature, locally enlarged]{\label{fig:fig3d}\includegraphics[width=0.5\textwidth]{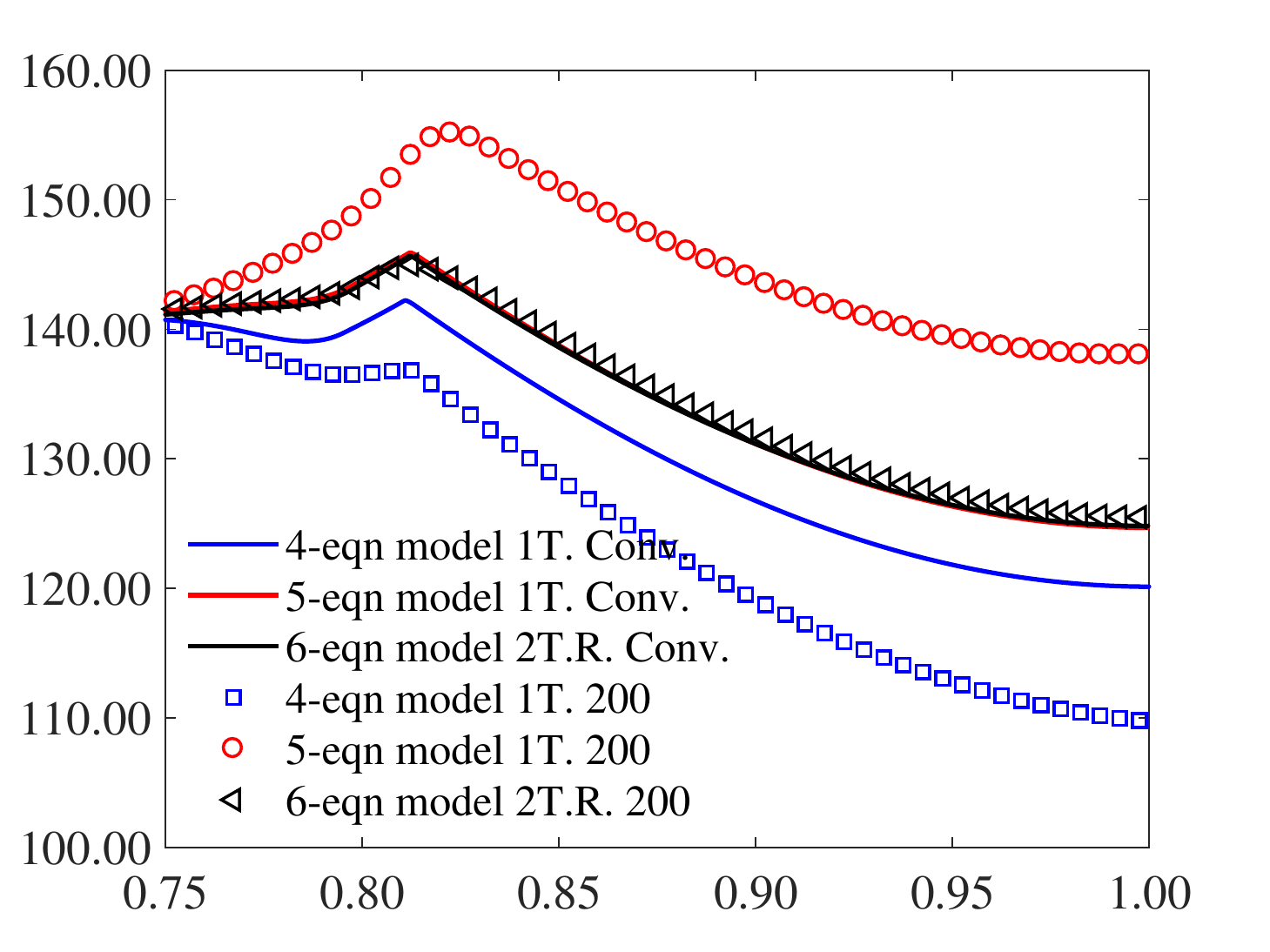}} \\
\subfloat[Velocity]{\label{fig:fig3e}\includegraphics[width=0.5\textwidth]{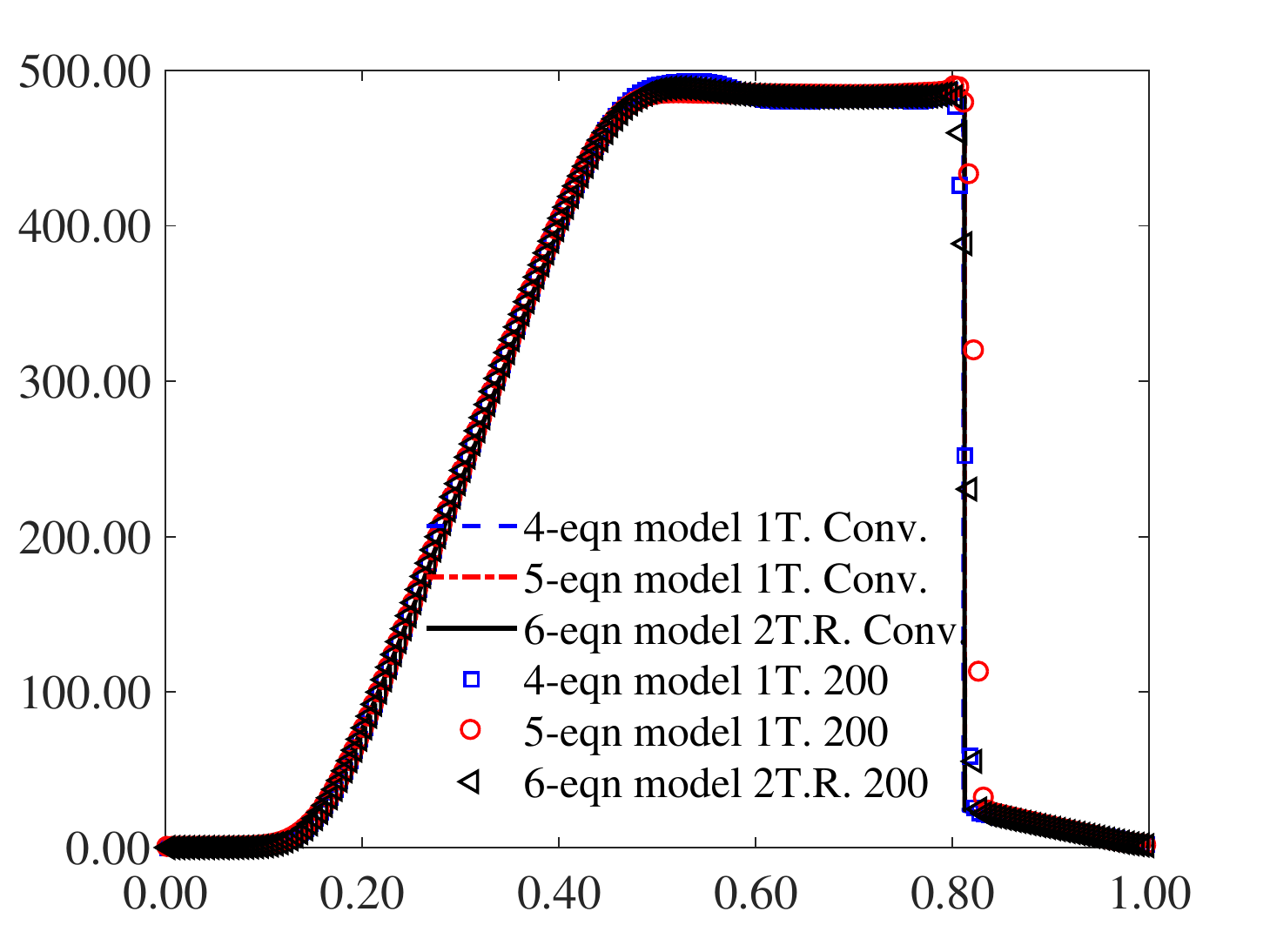}}
\subfloat[Velocity, locally enlarged]{\label{fig:fig3f}\includegraphics[width=0.5\textwidth]{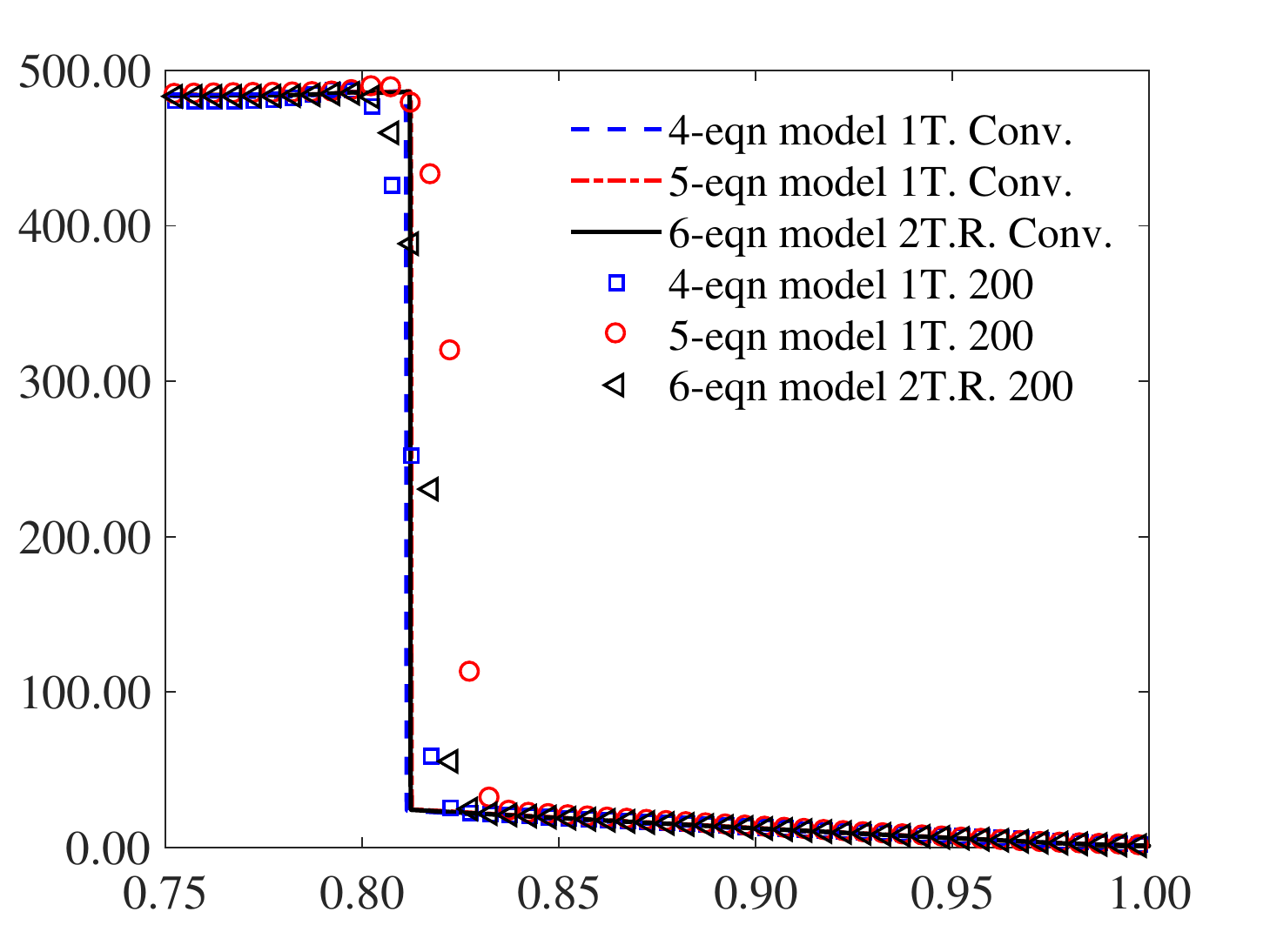}} 
\caption{Numerical results for the two-fluid shock tube problem with  equal phase  heat conductivity $\lambda = 1.00 \times 10^7$.}
\label{fig:fig3} 
\end{figure}

\paragraph{Test with equal phase thermal conductivity}
The above two-fluid shock tube problem is now considered with taking into account the phase heat conduction effect. The diffusion PDEs are solved with the explicit method of local iterations if not mentioned. The thermal conductivity is set to be a large number for comparison purpose. We first assign the same heat conduction coefficient for the two fluids $\lambda_1 \; = \lambda_2 = 1.00\times 10^6{\text{W}/(\text{m}\cdot\text{K})}$. The numerical results obtained with different models are compared in \Cref{fig:fig2}. The results marked as converged (``Conv'') are computed on a fine grid consisting of 20000 computational cells. {\color{black}The difference between the converged solutions of different models is indiscernible and they are taken as the reference solution for comparison.} To demonstrate the difference between the models, we show also the numerical solutions for a coarse grid of 200 cells. The results for the proposed model (6-eqn model 2T.R.) agree much better with the reference solution than the five-equation model (see \Cref{fig:fig2c,fig:fig2d}). 
The results of the four-equation model are also satisfactory since the heat conduction seems to be not strong enough to spread its erroneous temperature spikes. As the thermal conductivity is increased to $\lambda_1 \; = \lambda_2 = 1.00\times 10^7 {\text{W}/(\text{m}\cdot\text{K})}$, we find that these models do not converge to the same solution. This is demonstrated in \Cref{fig:fig3}, and is more clearly seen in the temperature profiles.   The results of the four-equation model on a 20000-cell grid diverge from those of the five-equation and six-equation models to the right of the material interface. This can be explained by the  numerical errors  in the {\color{black}diffused} zone, which then contaminate the results in the second fluid due to large thermal conductivity.

\begin{figure}[htb]
\centering
\subfloat[Density]{\label{fig:fig4a}\includegraphics[width=0.5\textwidth]{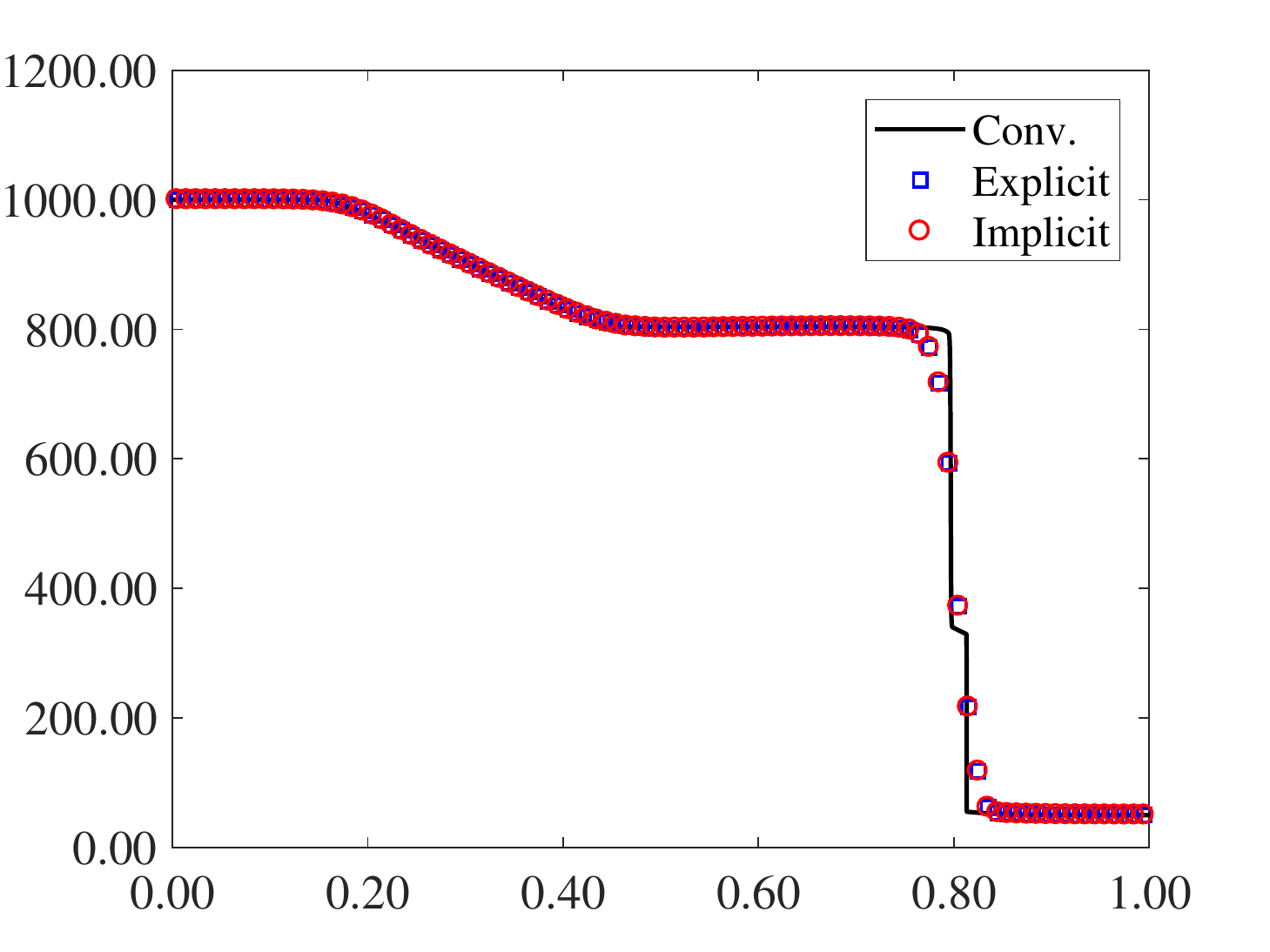}}
\subfloat[Temperature]{\label{fig:fig4b}\includegraphics[width=0.5\textwidth]{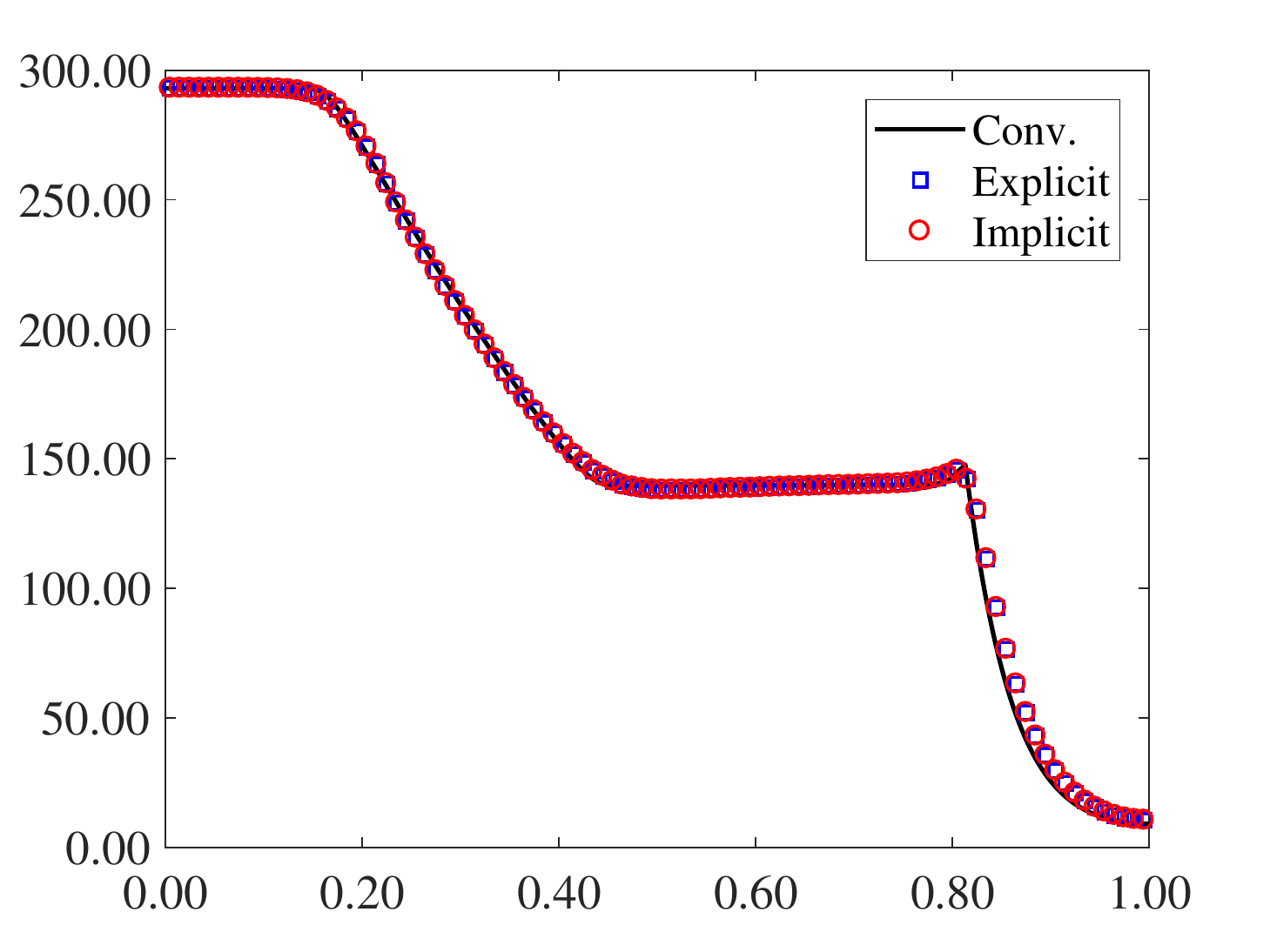}}
\caption{Numerical results for the two-fluid shock tube problem obtained with the explicit and implicit scheme on a 100-cell grid.}
\label{fig:fig4} 
\end{figure}

In \Cref{fig:fig4} we verify the explicit method of local iterations that is used to solve efficiently the parabolic part of the model (heat conduction). Here we compare this method with the implicit scheme solved by conventional Newtonian iterations. In the implicit scheme, the preconditioned conjugate gradient method is used for solving the system of algebraic equations. We see that the results obtained with both schemes on a 100-cell grid agree very well with the reference solution.

\paragraph{Test with non-uniform thermal conductivity and viscosity}

\begin{figure}[htb]
\centering
\subfloat[Temperature]{\label{fig:fig5a}\includegraphics[width=0.5\textwidth]{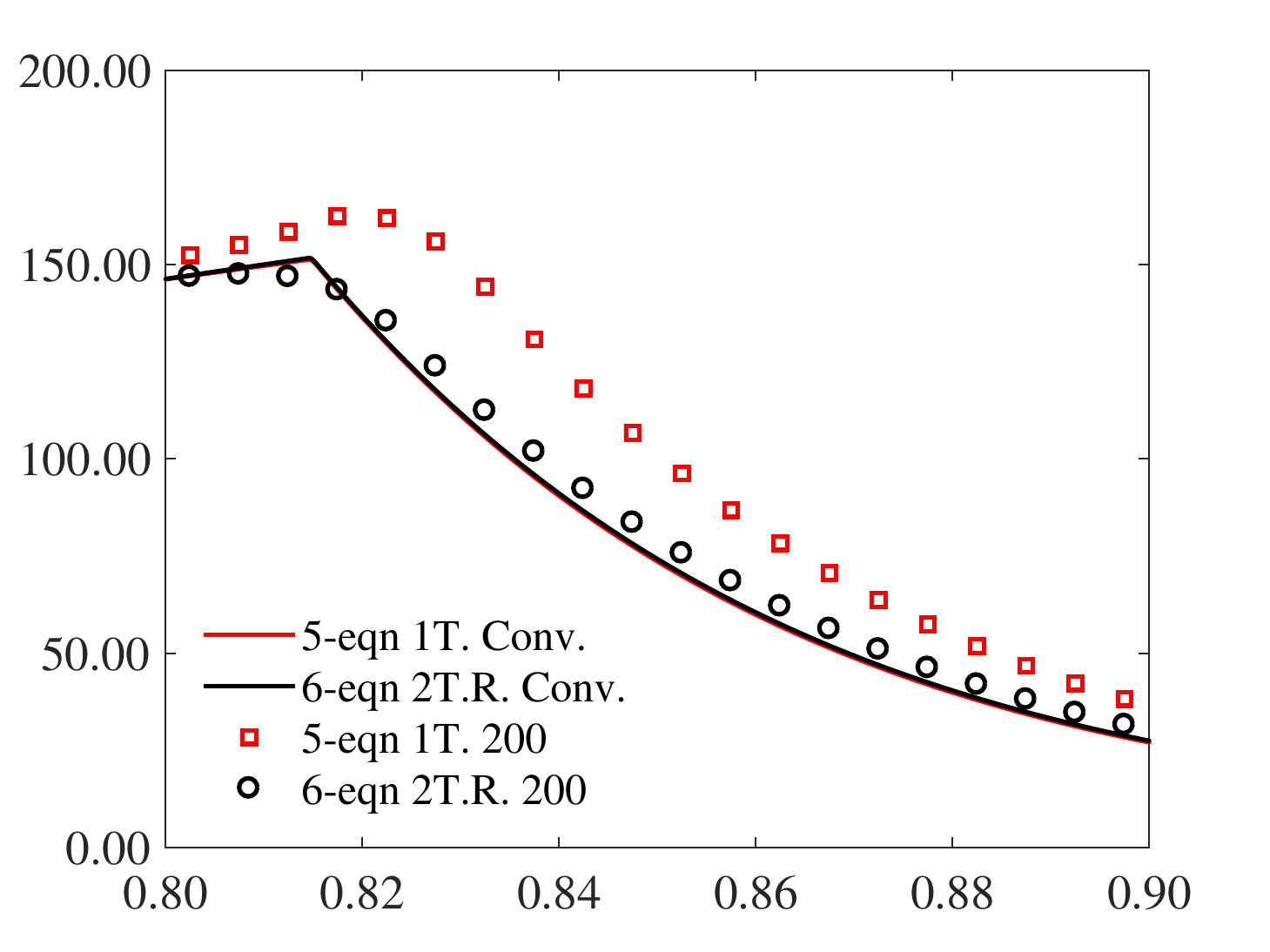}}
\subfloat[Velocity]{\label{fig:fig5b}\includegraphics[width=0.5\textwidth]{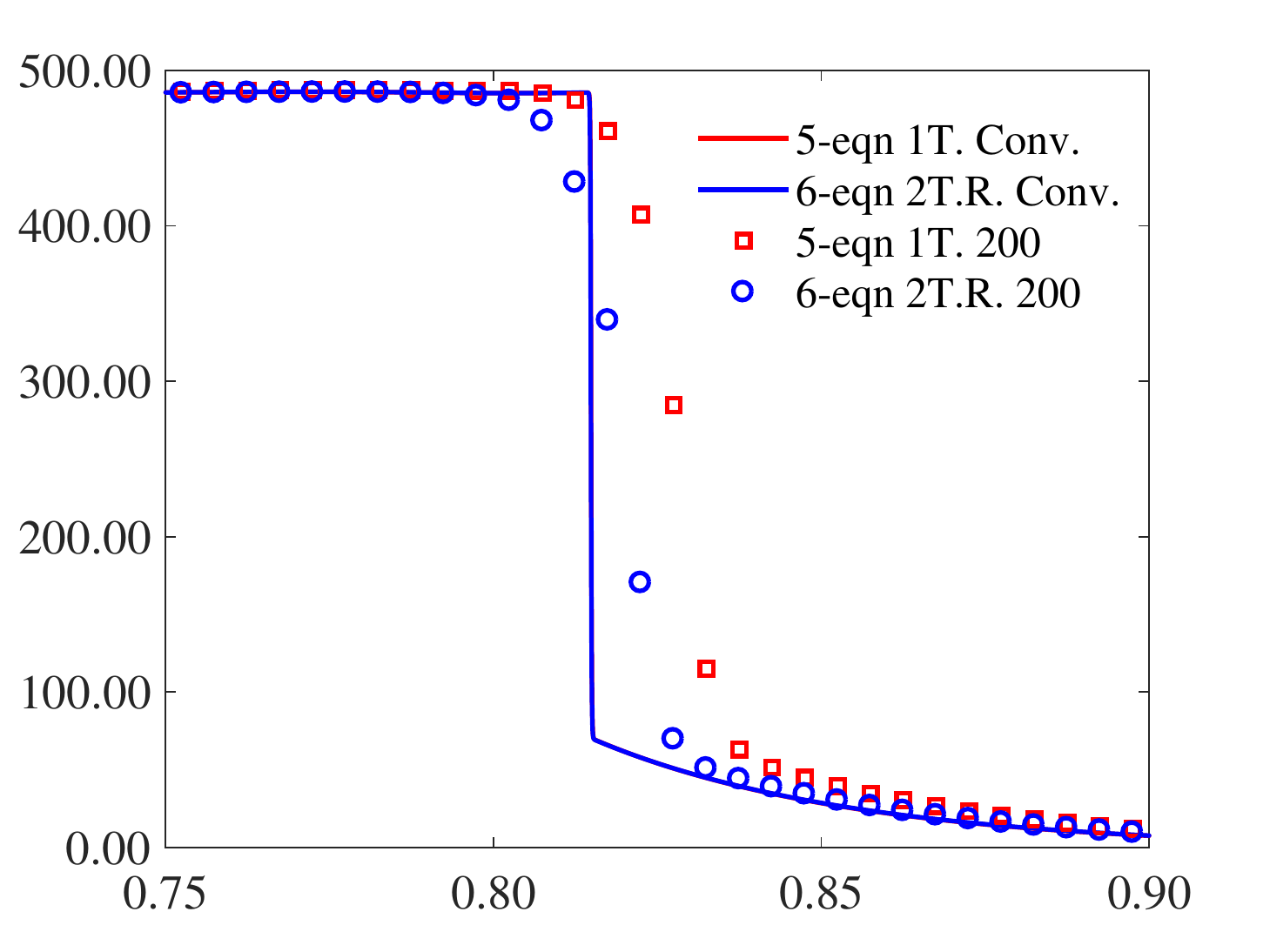}}
\caption{Numerical results for the two-fluid shock tube problem with for the case of different phase coefficients of viscosity and heat conductivity.}
\label{fig:fig5} 
\end{figure}

This test considers the shock tube problem for two fluids which have different thermal conductivities and viscosities. For the left fluid, the thermal conductivity and dynamic viscosity are assumed to be $\lambda_L = 1.00 \times 10^7{\text{W}/(\text{m}\cdot\text{K})}$ and $\mu_L = 5.00 \times 10^2{\text{Pa}\cdot\text{s}}$ and those for the right are $\lambda_R = 1.00 \times 10^6{\text{W}/(\text{m}\cdot\text{K})}$ and $\mu_R = 1.00{\text{Pa}\cdot\text{s}}$. The viscosity and thermal conductivity are averaged with volume fractions, i.e. $\lambda = \sum_{k}\alpha_k \lambda_k, \; \mu = \sum_{k}\alpha_k \mu_k$. The four-equation model does not provide solution of the volume fraction. Therefore, we test only the five-equation model and the six-equation model. From \Cref{fig:fig5} one can see that the convergence performance of the proposed six-equation model is still superior to that of the five-equation model.

\subsection{Two-phase problem}
In this section, numerical experiments are performed for two-phase flows where the phases are mixed and may occupy the same location in space.   

\paragraph{Shock propagation in a water-gas mixture}
The material properties  of the phases  are the same as in the previous  test. The volume fraction of each component is initially 0.50 in the whole computational domain.  Other initial data is given as follows:
\begin{align}\nonumber
x<0.5{\text{m}} : P = 1.00 \times 10^{9}{\text{Pa}}, \; T = 1000{\text{K}}; \\
\nonumber x>0.5{\text{m}} : P = 1.00 \times 10^{5}{\text{Pa}}, \; T = 300{\text{K}}.
\end{align}

\begin{figure}[htb]
\centering
\subfloat[Density]{\label{fig:MPWG_1}\includegraphics[width=0.5\textwidth]{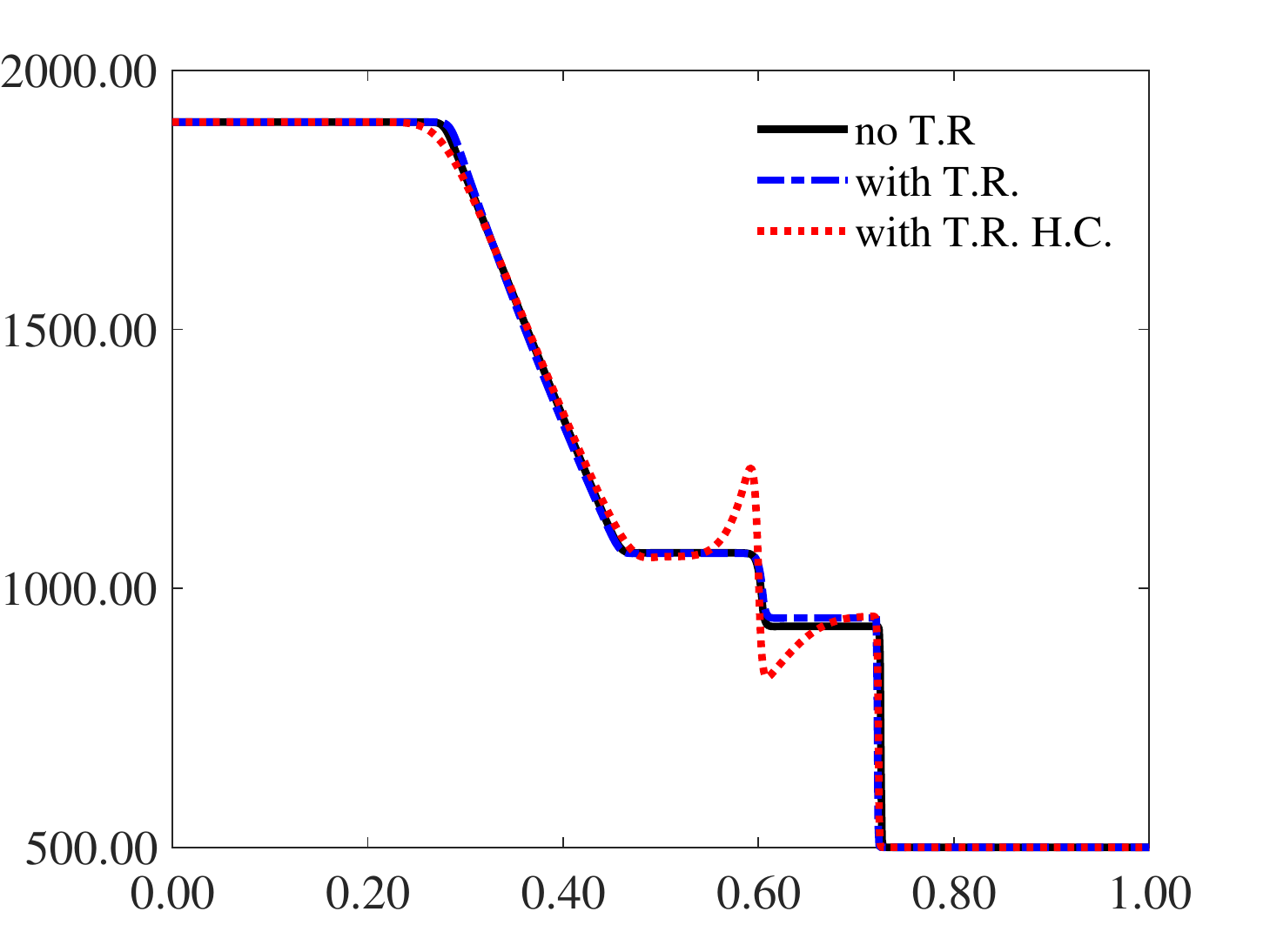}}
\subfloat[Temperature]{\label{fig:MPWG_2}\includegraphics[width=0.5\textwidth]{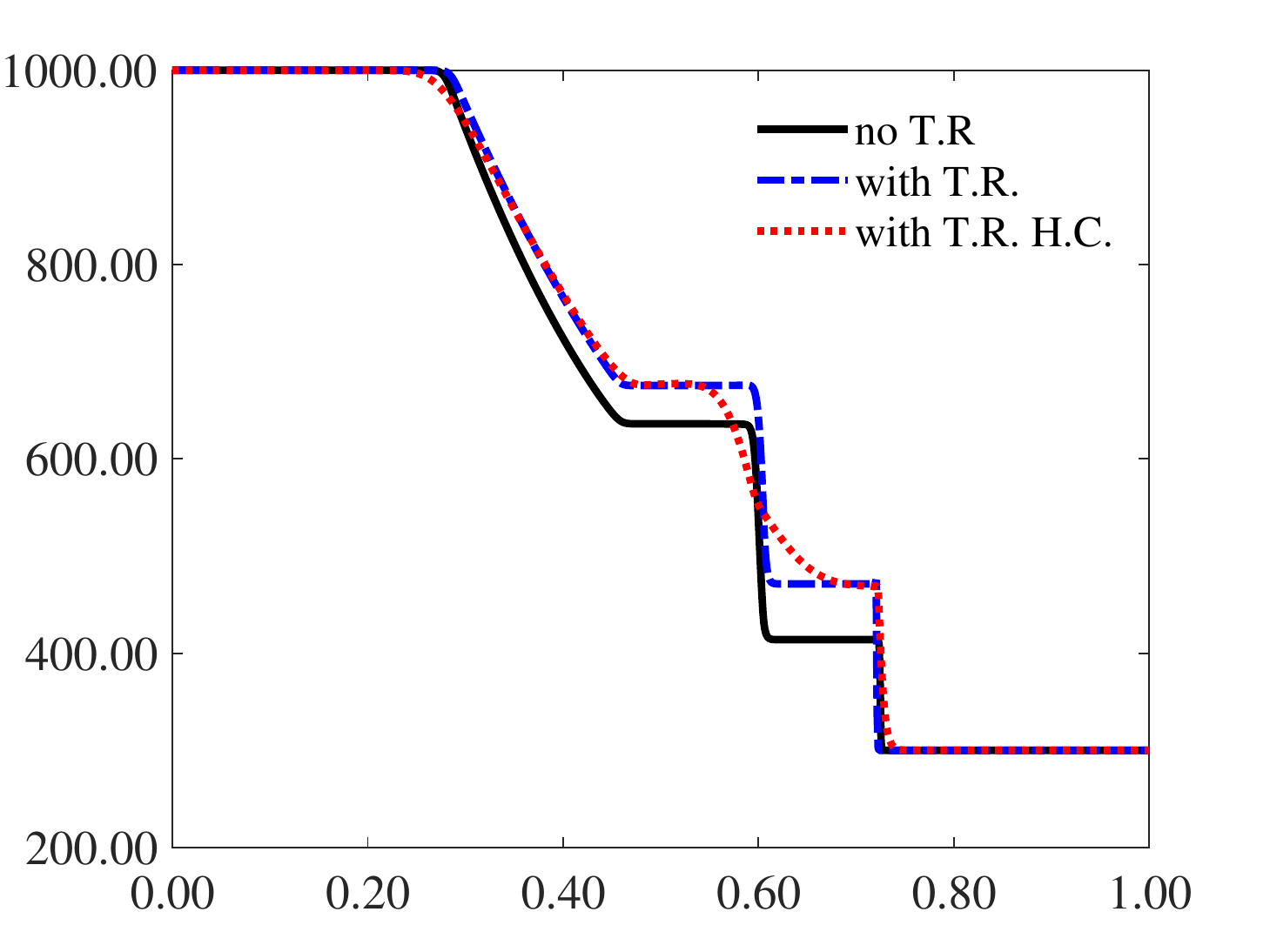}}\\
\subfloat[Volume fraction]{\label{fig:MPWG_3}\includegraphics[width=0.5\textwidth]{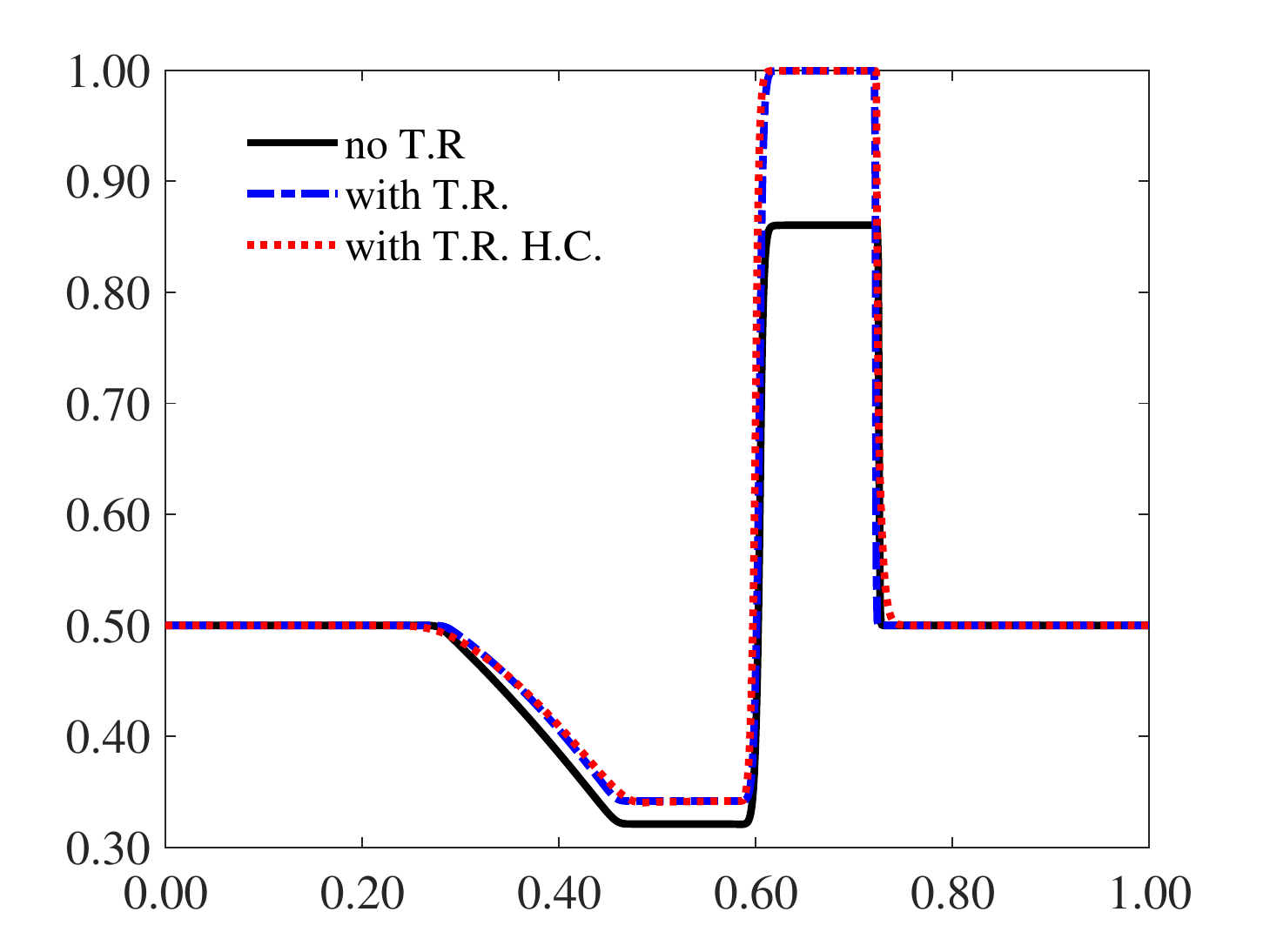}}
\subfloat[Pressure]{\label{fig:MPWG_4}\includegraphics[width=0.5\textwidth]{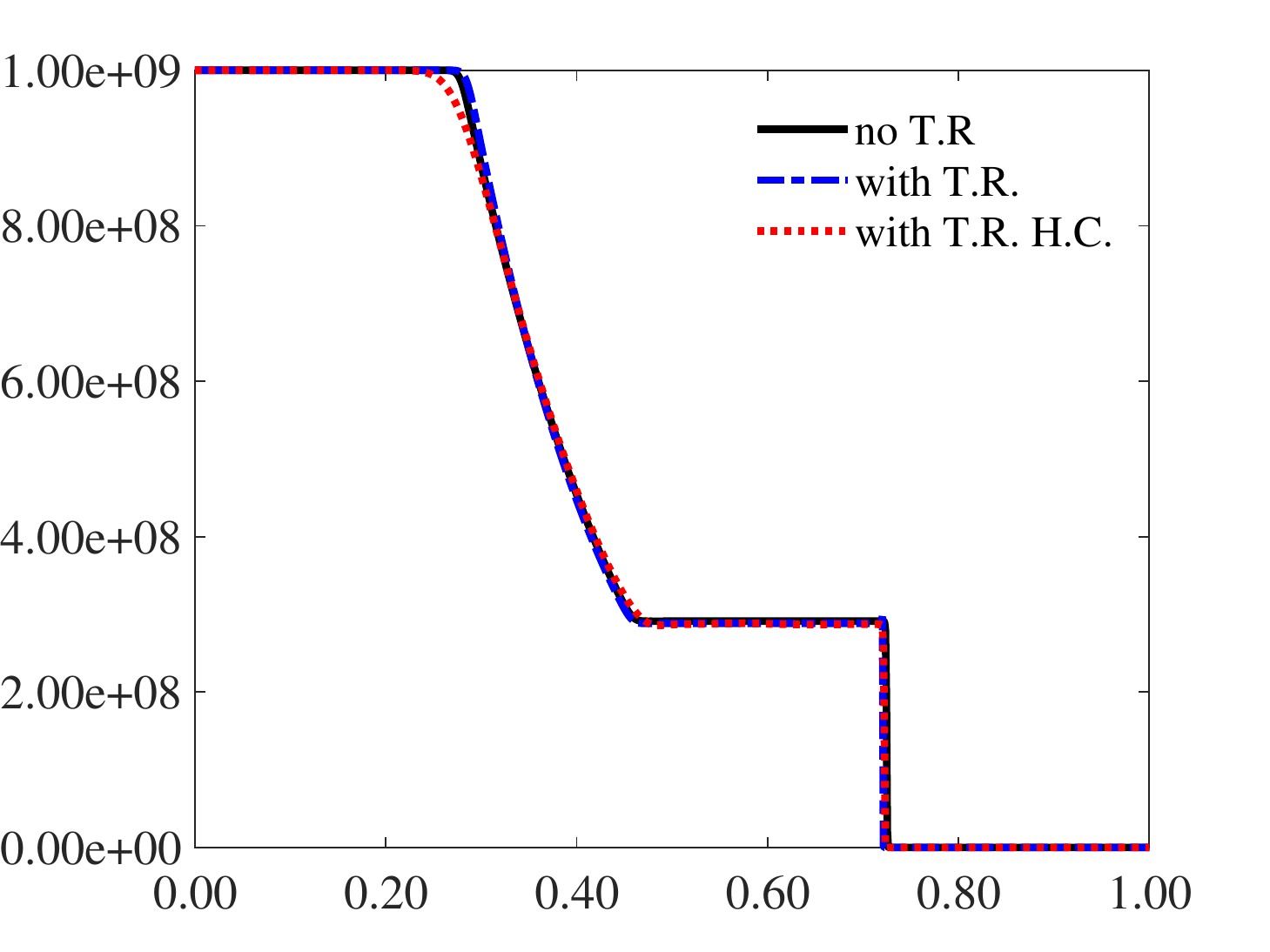}}
\caption{The numerical results for the water-gas multiphase shock tube problem. The lines `no T.R.',`with T.R.',`with T.R. H.C.',  represent the numerical results without temperature relaxation, with temperature relaxation, with temperature relaxation and heat conduction, respectively.}
\label{fig:MPWG} 
\end{figure}

Initial densities are determined by the corresponding EOS of each phase.
 For comparison purpose, the conductivities of gas and water are set by effective values of $1.00 \times 10^{7}{\text{W/m/K}}$ and $1.00 \times 10^{5}{\text{W/m/K}}$, respectively. Computations are performed to the moment $t=2.00 \times 10^{-4}$s on a 1000-cell uniform grid. The numerical results are shown in \Cref{fig:MPWG}. One can see that including temperature relaxation changes considerably the solution. The heat conduction process smears the temperature profile near the contact discontinuity, also resulting in corresponding changes in other variables.

\paragraph{Shock wave in solid alloys}
We further consider an alloy impact problem from  \cite{murrone2005five}. 
The alloy is composed of two components epoxy and spinel. The volume fractions of these two components are 0.595 and 0.415, respectively. We solve this problem as a two-phase one with the six-equation model. The materials are characterized by the following EOS parameters:
\begin{enumerate}
\item[•] Epoxy -- $\gamma = 2.94, \quad P_{\infty,1} = 3.20 \times 10^{9}{\text{Pa}}, \quad \rho = 1185.00{\text{kg}/\text{m}^3}$,
\item[•] Spinel -- $\gamma = 1.62, \quad P_{\infty,2} = 1.41 \times 10^{11}{\text{Pa}}, \quad \rho = 3622.00{\text{kg}/\text{m}^3}$.
\end{enumerate}

\begin{figure}[htb]
\centering
\includegraphics[width=0.65\textwidth]{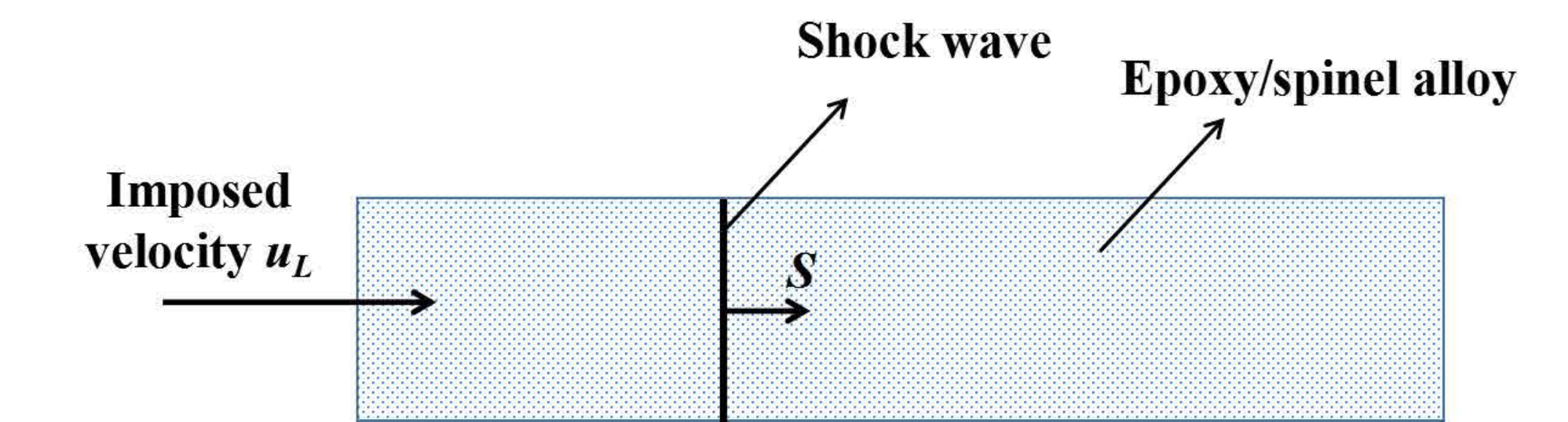}
\caption{Schematic of the alloy impact problem.}
\label{fig:alloy_impact_schematic} 
\end{figure}

The schematic of this problem is displayed in \Cref{fig:alloy_impact_schematic}. Calculations are carried out in the model without thermal relaxation, as the time scale of this problem is much smaller than the characteristic relaxation time. However, the mechanical relaxation is implemented.   

\begin{figure}[htb]
\centering
\includegraphics[width=0.55\textwidth]{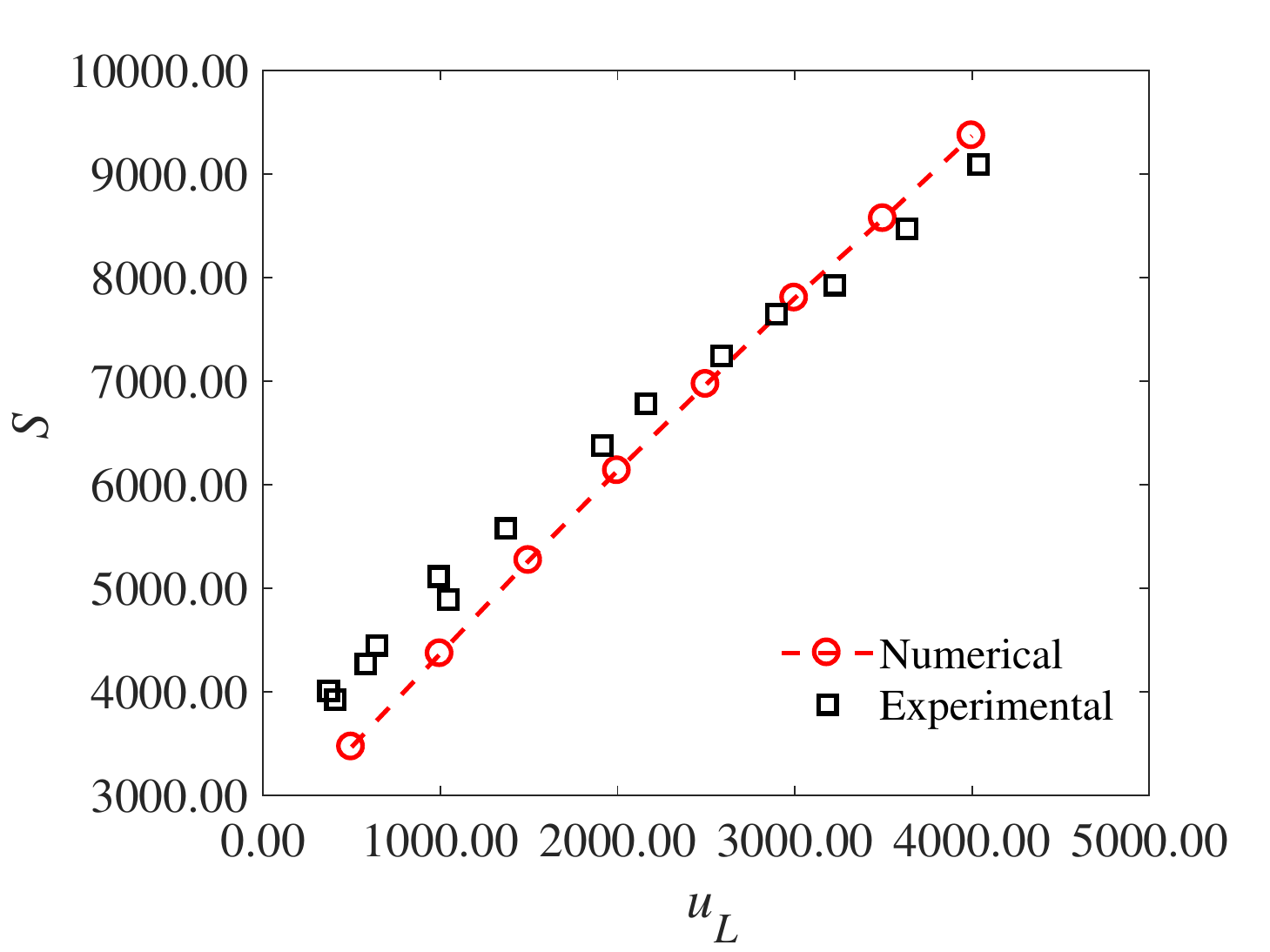}
\caption{Shock velocity diagramm for the epoxy/spinel alloy impact problem. The experimental data is from \cite{murrone2005five}.}
\label{fig:AlloyImpactSU} 
\end{figure}

For many metals, the shock velocity $S$ linearly depends on the impact velocity $u_L$. Calculations of the shock wave propagation are done for different velocities $u_L$. The results obtained are plotted in \Cref{fig:AlloyImpactSU}  and compared with the experimental data that is available from \cite{murrone2005five}. As can be seen, a linear profile of $S$ well agrees with the experimental data.

\subsection{Laser ablation problem}\label{subsec:ablation}

\begin{figure}[htb]
\centering
\subfloat[Problem setup]{\label{fig:pabsetupa}\includegraphics[width=0.5\textwidth]{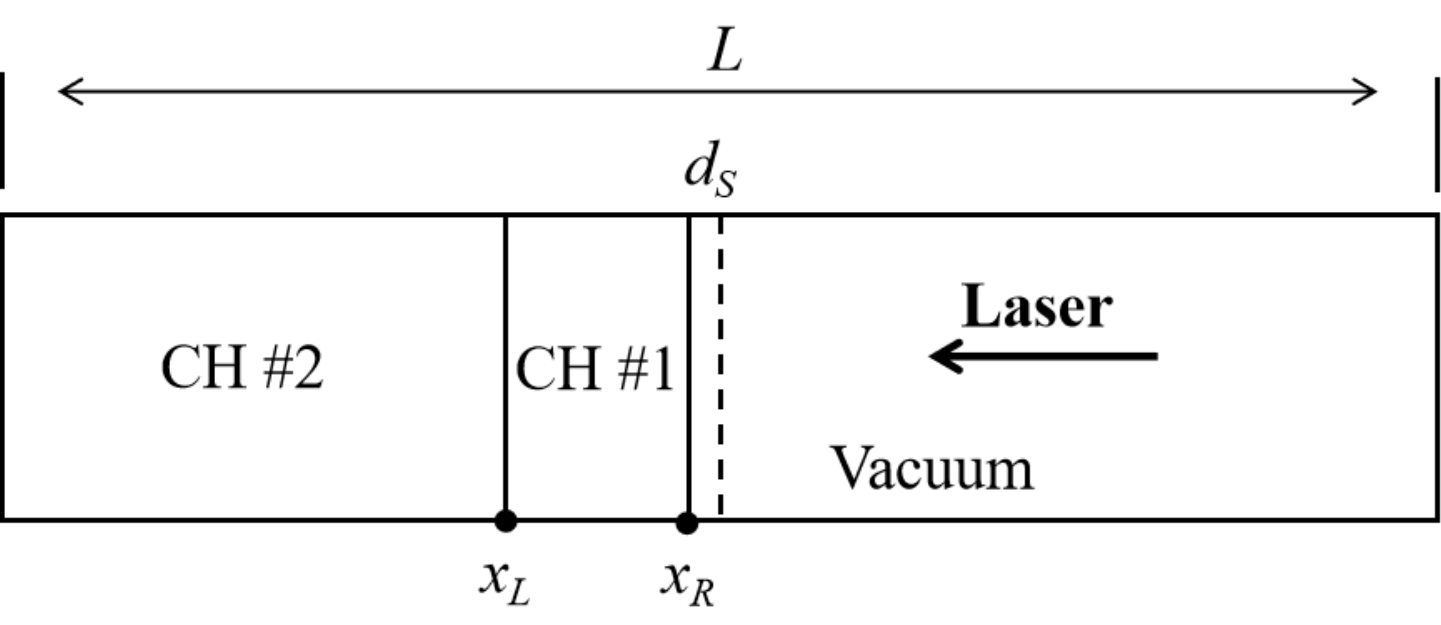}} \\
\subfloat[Critical density $\rho_{crt}$ and absorption distance $d_S$]{\label{fig:pabsetupb}\includegraphics[width=0.45\textwidth]{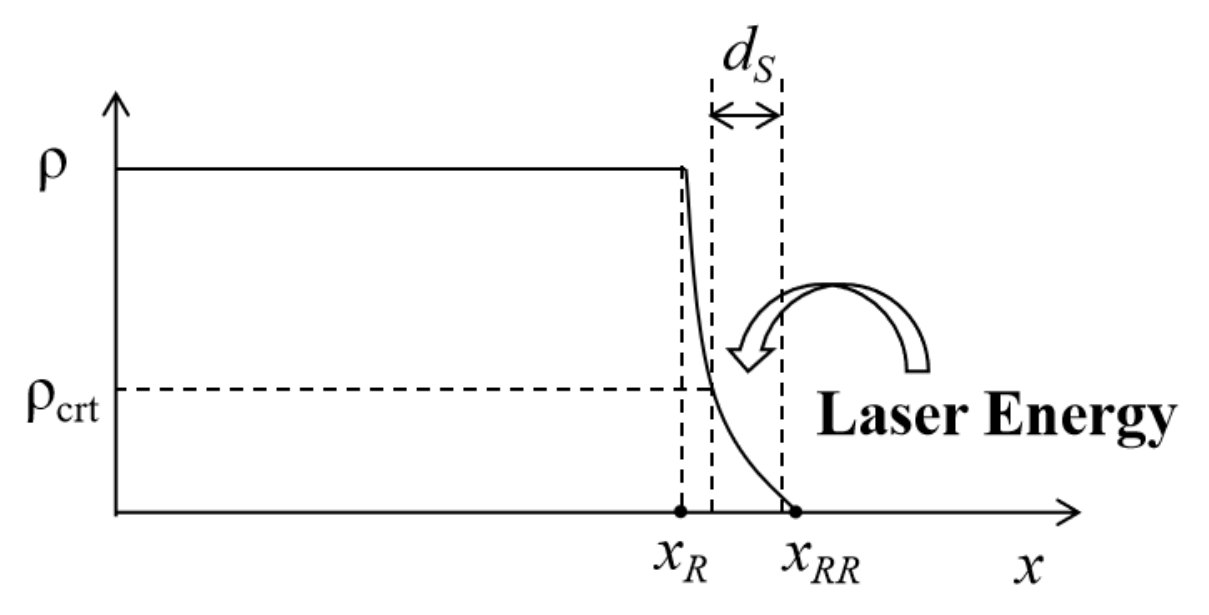}}
\caption{The laser ablation of a multicomponent planar target.}
\label{fig:pabsetup} 
\end{figure}


This section addresses an application problem related to the field of inertial confinement fusion (ICF) -- the laser ablation problem. In the direct-drive ICF capsule, the laser is used as an energy source to accelerate the plastic (CH, i.e., phenylethylene C$_8$H$_8$) target  creating high temperature and high pressure environment for inward implosion.

\paragraph{One-dimensional planar target}
First, the laser ablation problem is considered in the 1D approximation. It is assumed that the target is plane, and the laser emission is uniform and normal to the the target surface. The target is placed in vacuum that is approximated as a fluid with extremely low density. As shown in  \Cref{fig:pabsetup}, the laser radiation comes from the right and its energy is absorbed by the CH material that then turns to high temperature ablated plasma.  
The energy absorption occurs up to the critical density point (where the incident power energy equals the reflected one) and over a distance $d_{S}$ (absorption area) to the right of the critical density point.
We consider the following composite target consisting of two different CH materials separated by the material interface at a distance $x_L$
\begin{enumerate}
\item[•] CH \#1 \quad $\rho = 1.00, \; \gamma = 1.666, \; C_v = 86.27, $
\item[•] CH \#2 \quad $\rho = 0.80, \; \gamma = 1.220, \; C_v = 76.27, $
\end{enumerate}
hereinafter dimensions used are centimeter, gramm and microsecond.

The vacuum is approximated as the material CH \#2 with a density of $8.00\times 10^{-6}$. The critical density is $\rho_{crt} = 0.39$ that can be calculated according to the inverse bremsstrahlung absorption theory.

The electron, ion and photon in the plasma  are assumed to be in thermal equilibrium. 
The thermal conductivity of the plasma is approximated with the one-temperature Spitzer-Harm model \cite{Spitzer1953,Li2020Numerical} and  is a nonlinear function of density and temperature: 
\begin{equation}
\lambda_{{SH}}=9.44\left(\frac{2}{\pi}\right)^{3 / 2} \frac{\left(k_{{B}} T_{{e}}\right)^{5 / 2} k_{{B}} N_{\mathrm{e}}}{\sqrt{m_{{e}}} e^{4}} \frac{1}{N_{{i}} Z_{{e}}\left(Z_{{e}}+4\right) \ln \Lambda_{{ei}}},
\end{equation}
where $k_B$ is the Boltzmann constant, $T_e$ is the electronic temperature, $N_e$ is the electron density,  $e$ is the electronic charge, $m_e$ is the electronic mass, $N_i$ is the ion density, $Z_e$ is the degree of ionization.
For a certain plasma,
\begin{equation}
N_i = \frac{N_0}{A_c} \rho, \quad N_e = Z_e N_i,
\end{equation}
where $A_c$ is the average atomic number, $N_0$ is the Avogadro's number.

$\ln \Lambda_{{ei}}$ is the Coulomb logarithm of laser absorption and determined with
\begin{equation}
\ln \Lambda_{{ei}}=\left\{\begin{array}{ll}
\max \left(1, \ln \frac{l_{{D}}}{l_{{LD}}}\right), & \frac{Z_e^{2}}{3 k_{{B}} T_{{e}}} \geq l_{{dB}}, \\
\max \left(1, \ln \frac{l_{{D}}}{l_{{dB}}}\right), & \frac{Z_e^{2}}{3 k_{{B}} T_{{e}}}<l_{{dB}}.
\end{array}\right.
\end{equation}
where $l_D$ is Debye length, $l_{LD}$ is Landau length, $l_{dB}$ is De Broglie wavelength.

When each component obeys their own Spitzer-Harm relation, the four-equation model is not applicable since the conductivity is averaged by the volume fraction that is absent in this model. Therefore, for comparison purpose, we assume equal phase conductivity defined with the same Spitzer-Harm relation.

Within the absorption distance $d_S = 2.00 \times 10^{-3}$, the deposited laser power intensity is assumed to be constant, $\mathcal{I} = 1.00 \times 10^{3}$. In the vicinity of the right interface, the density is smoothed in the region from $x_R$ to $x_{RR}$ by an exponential function of the spatial coordinate. The geometry of the computational domain is specified as $L = 1.00 \times 10^{-1}$ and $x_L = 0.45L, \; x_R = 0.50L, \; x_{RR} = 0.51L$.
The initial temperature $T = 3.00 \times 10^{-4}$ in the whole computational domain. Pressure is calculated with the EOS of each material. 

Calculations are performed with three models (the four-equation model, the one-temperature five-equation model and the proposed six-equation model) and  two grids consisting of 1200 and 9600 equally distributed cells. To improve the material interface resolution, we implement the MUSCL scheme with the  Overbee limiter \cite{Chiapolino2017}. This scheme is applied to phase masses $\alpha_k \rho_k$ and the volume fraction $\alpha_1$ for the five- and six-equation model, and to mixture density $\rho$ and mass fraction $Y_1$ for the four-equation model.  

\begin{figure}[htb]
\centering
\subfloat[Density]{\label{fig:PABLA}\includegraphics[width=0.5\textwidth]{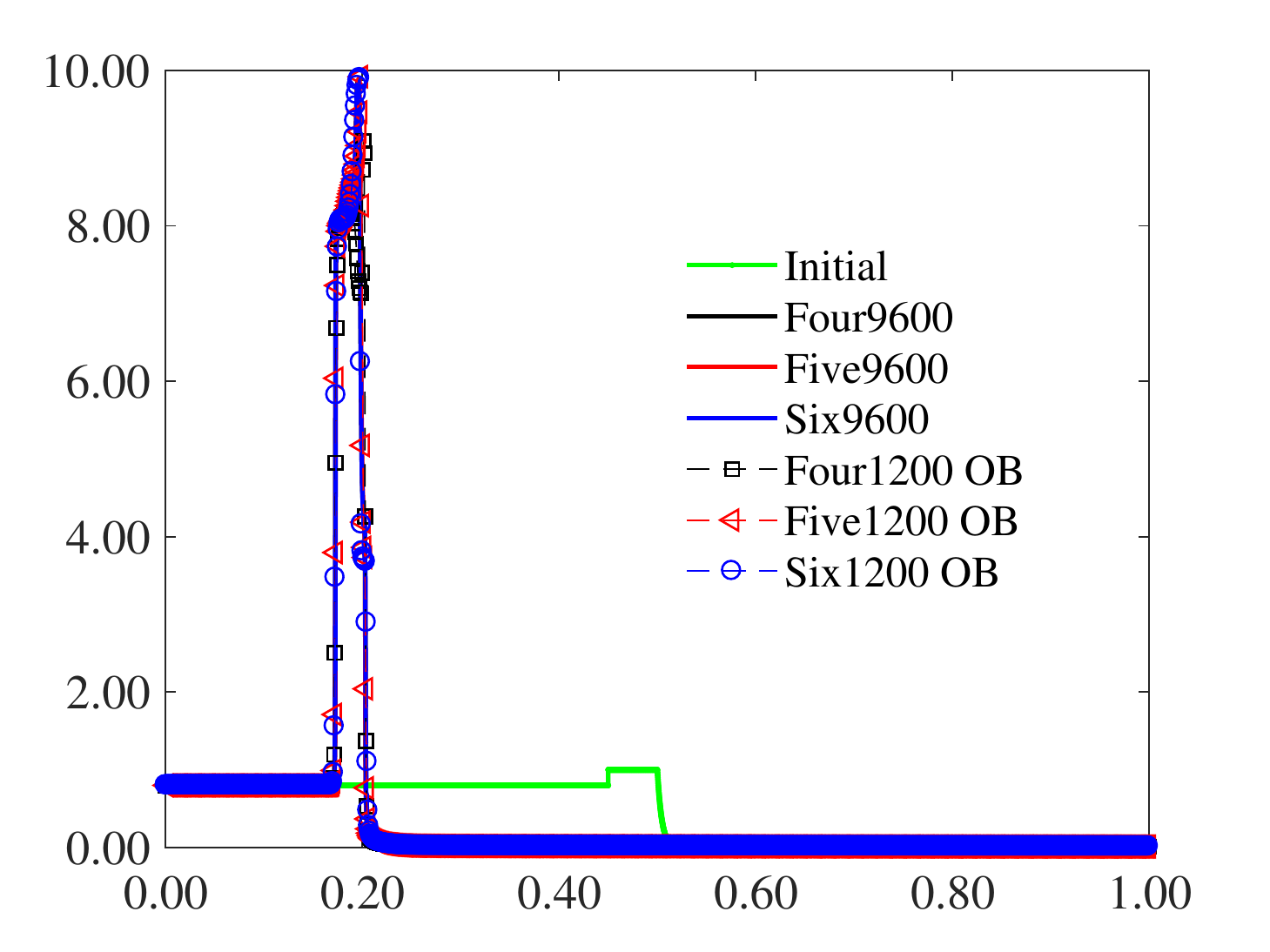}}
\subfloat[Density, locally enlarged]{\label{fig:PABLB}\includegraphics[width=0.5\textwidth]{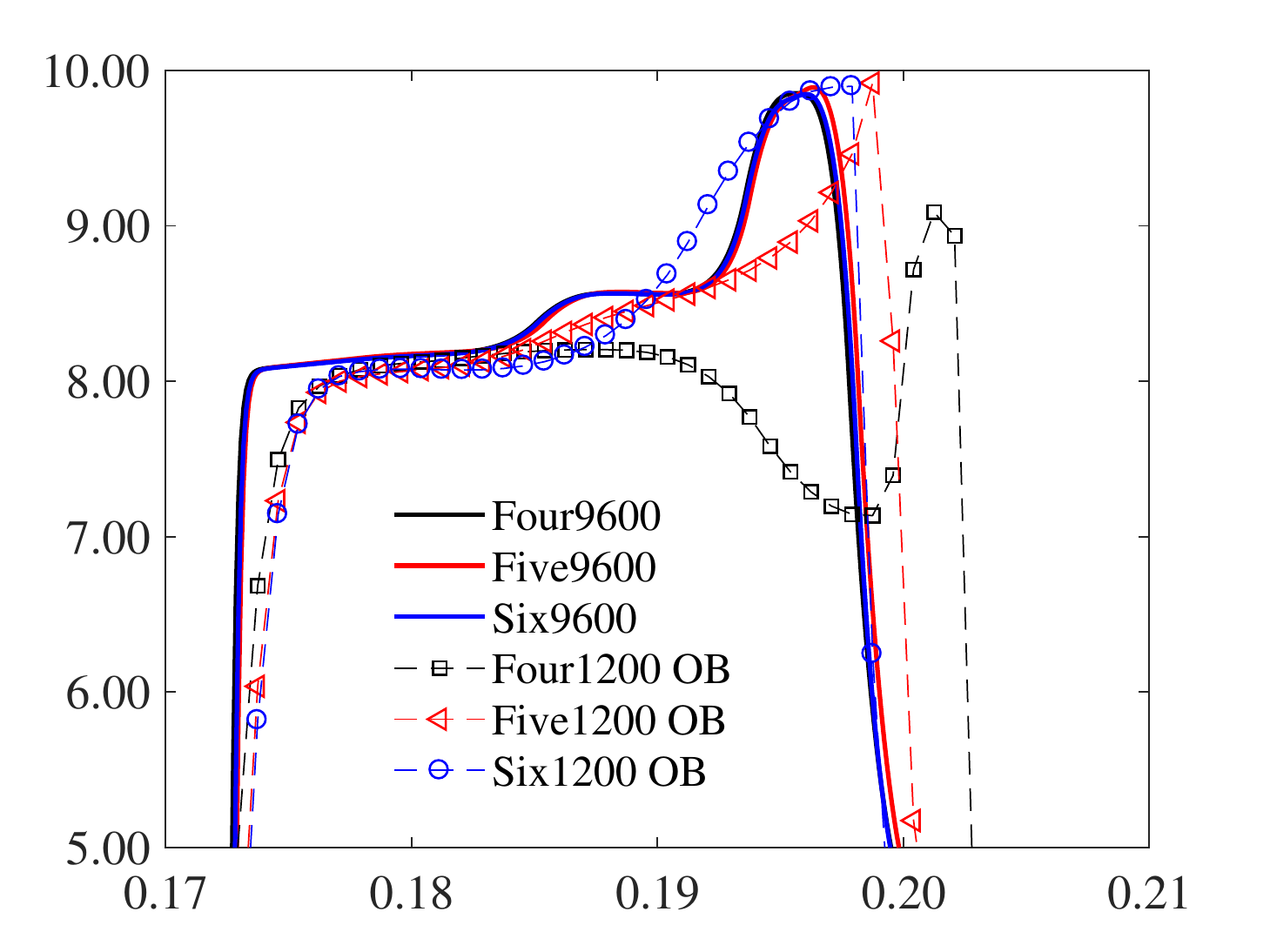}} \\
\subfloat[Volume/mass fraction]{\label{fig:PABLC}\includegraphics[width=0.5\textwidth]{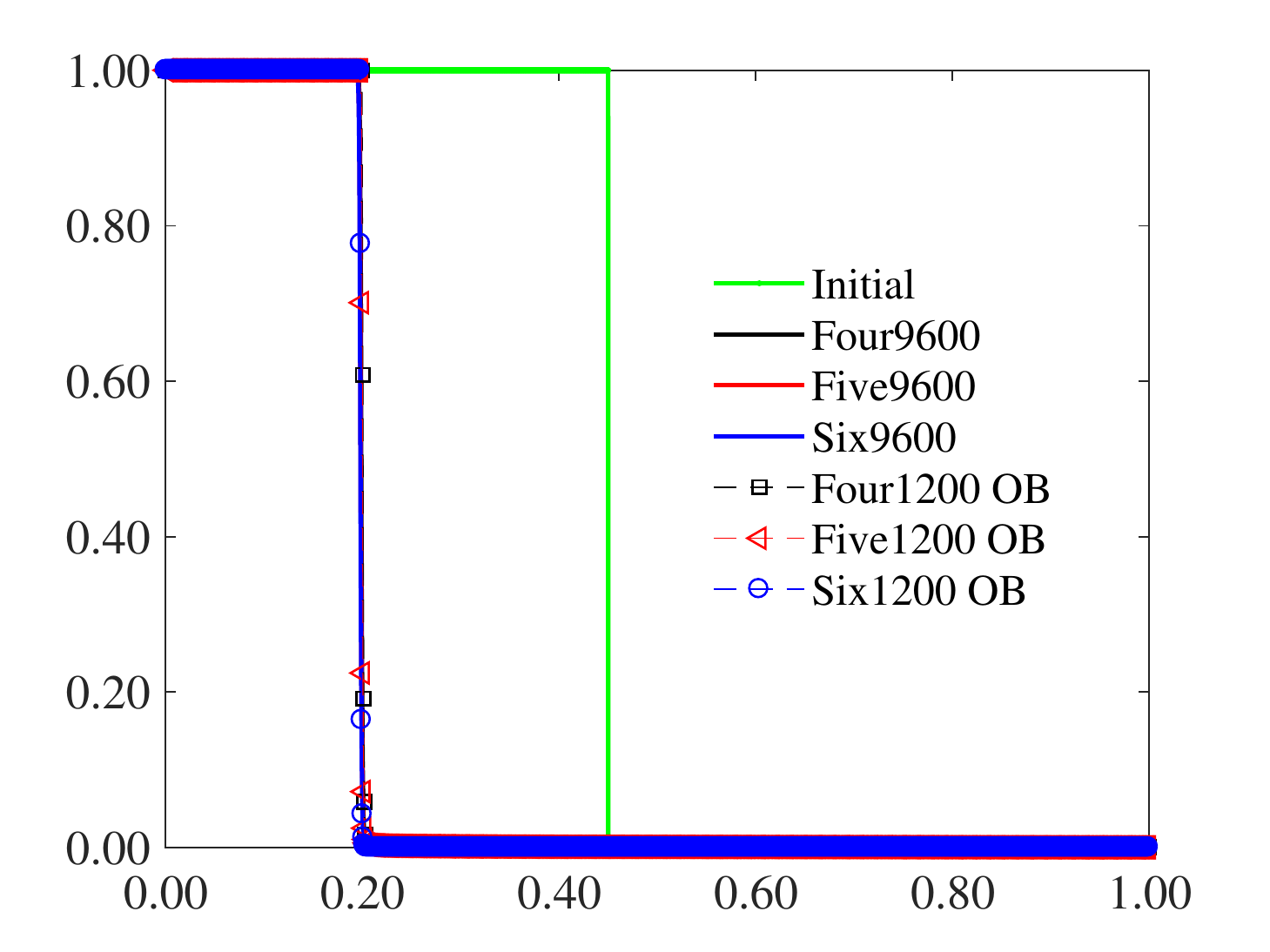}}
\subfloat[Volume/mass fraction, locally enlarged]{\label{fig:PABLD}\includegraphics[width=0.5\textwidth]{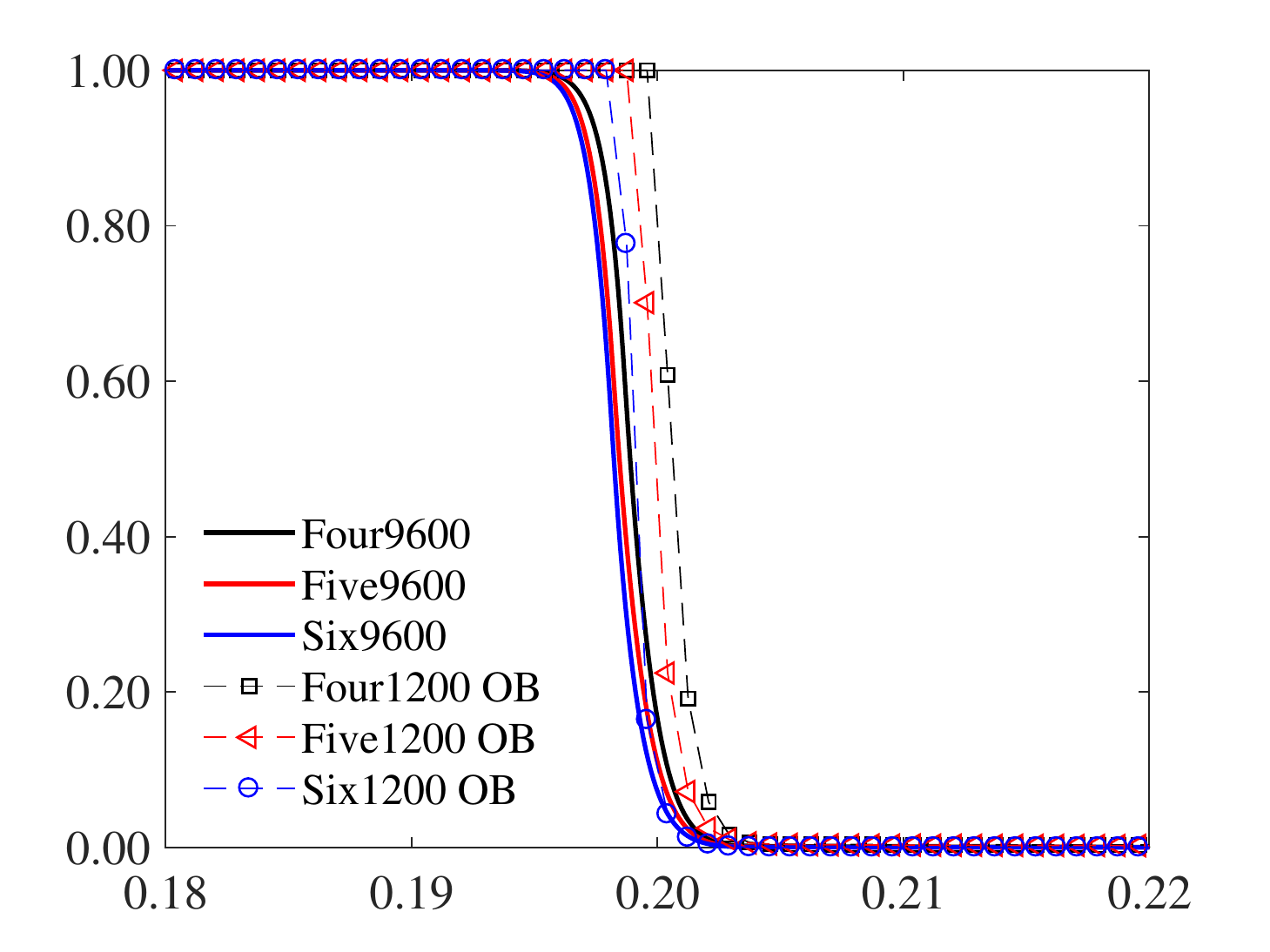}} \\
\subfloat[Temperature]{\label{fig:PABLE}\includegraphics[width=0.5\textwidth]{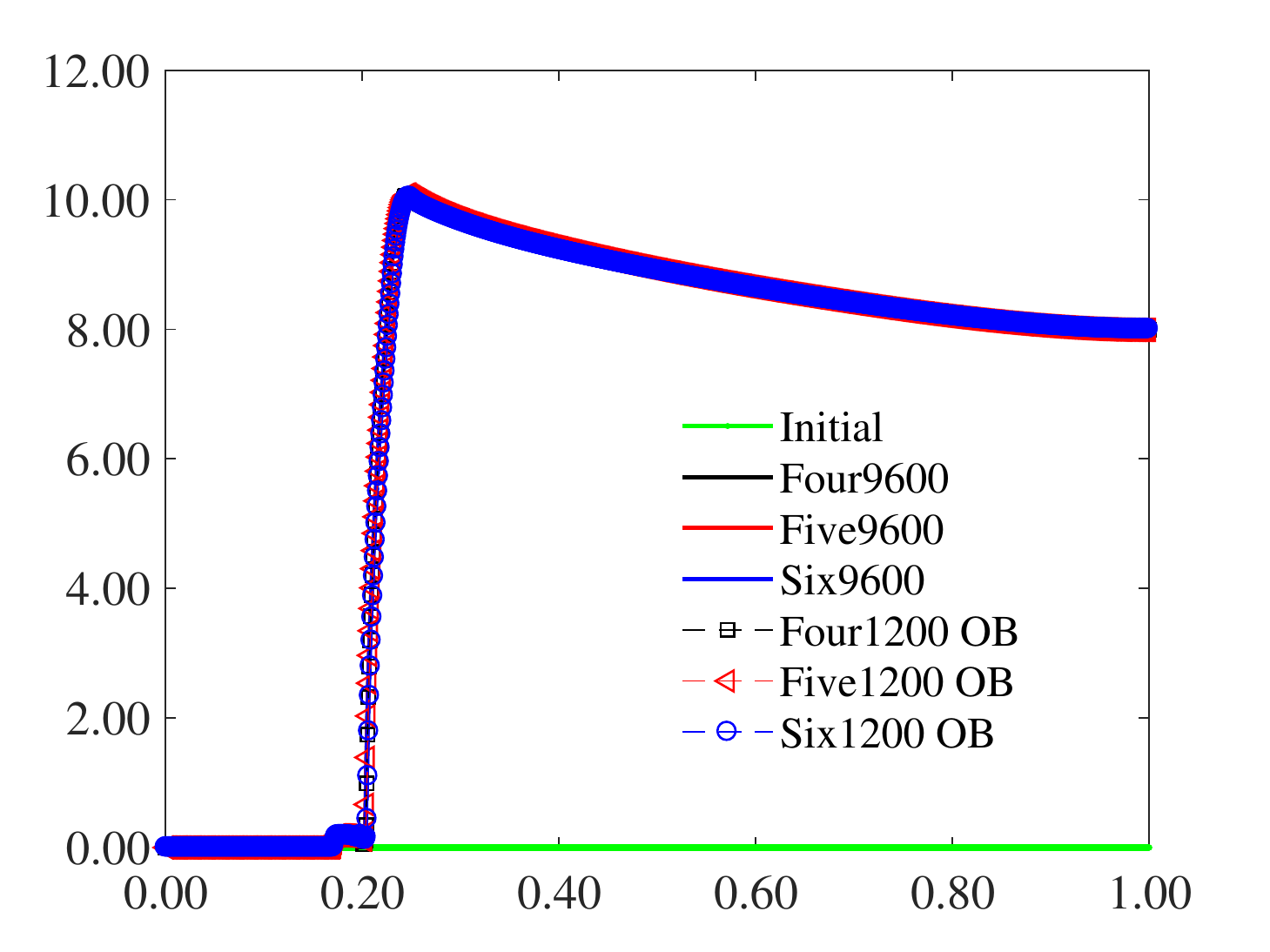}}
\subfloat[Temperature, locally enlarged]{\label{fig:PABLF}\includegraphics[width=0.5\textwidth]{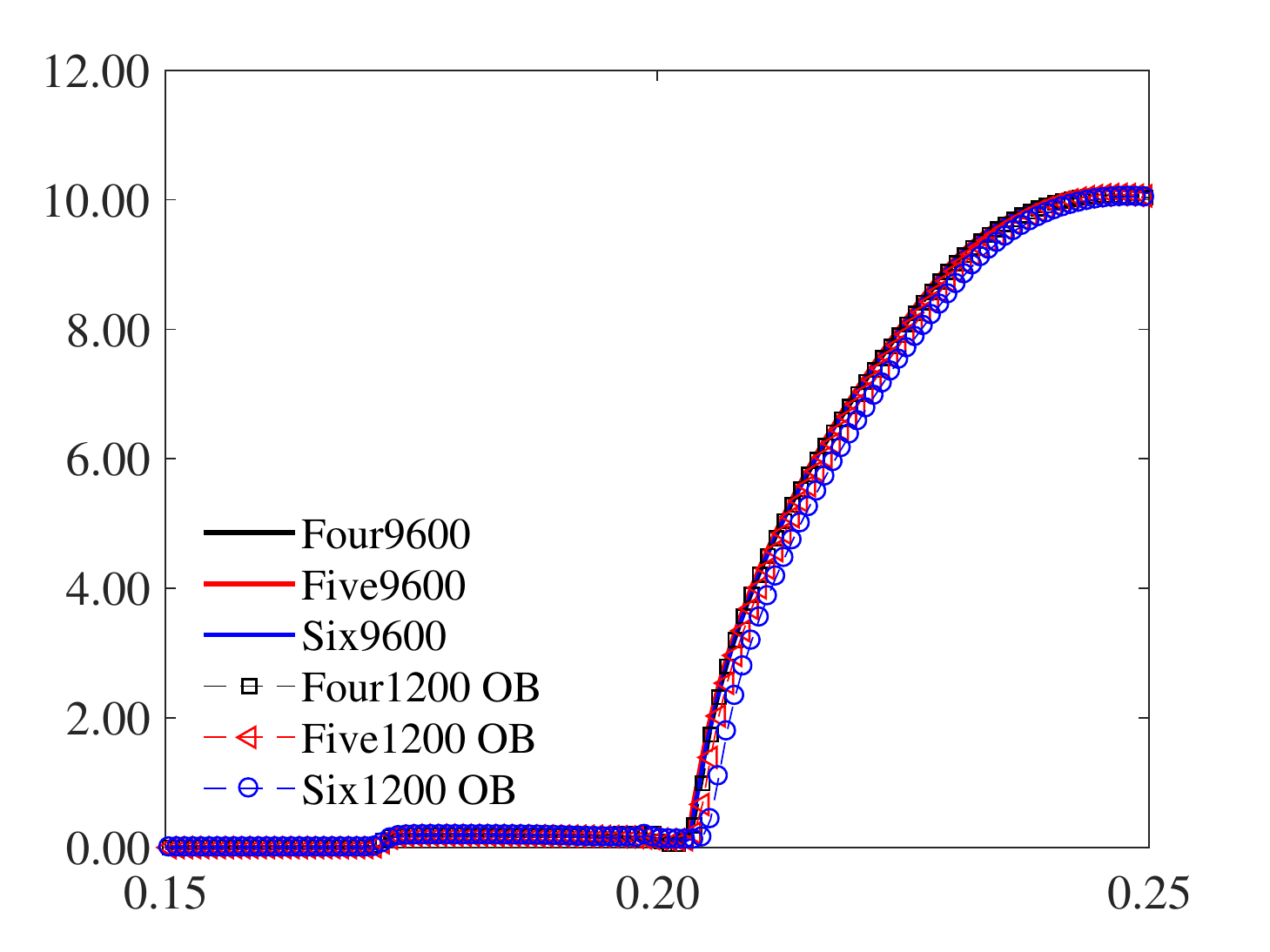}}
\caption{Numerical results for the laser ablation of a multicomponent planar target.}
\label{fig:PABL} 
\end{figure}

The numerical results of these three models at $t=6.00\times 10^{-3}$ are compared in \Cref{fig:PABL}. All three models tend to converge to the same solution.
The results show only minor differences. For example, convergence in density for the four-equation model and convergence in temperature for the six-equation model are found to be worse in comparison with the other two (\Cref{fig:PABLA,fig:PABLB,fig:PABLE,fig:PABLF}). There is also small difference in the interface velocity, as seen in \Cref{fig:PABLC,fig:PABLD}.  

\begin{figure}[htb]
\centering
\subfloat[Density, locally enlarged]{\label{fig:MUS_OB1}\includegraphics[width=0.5\textwidth]{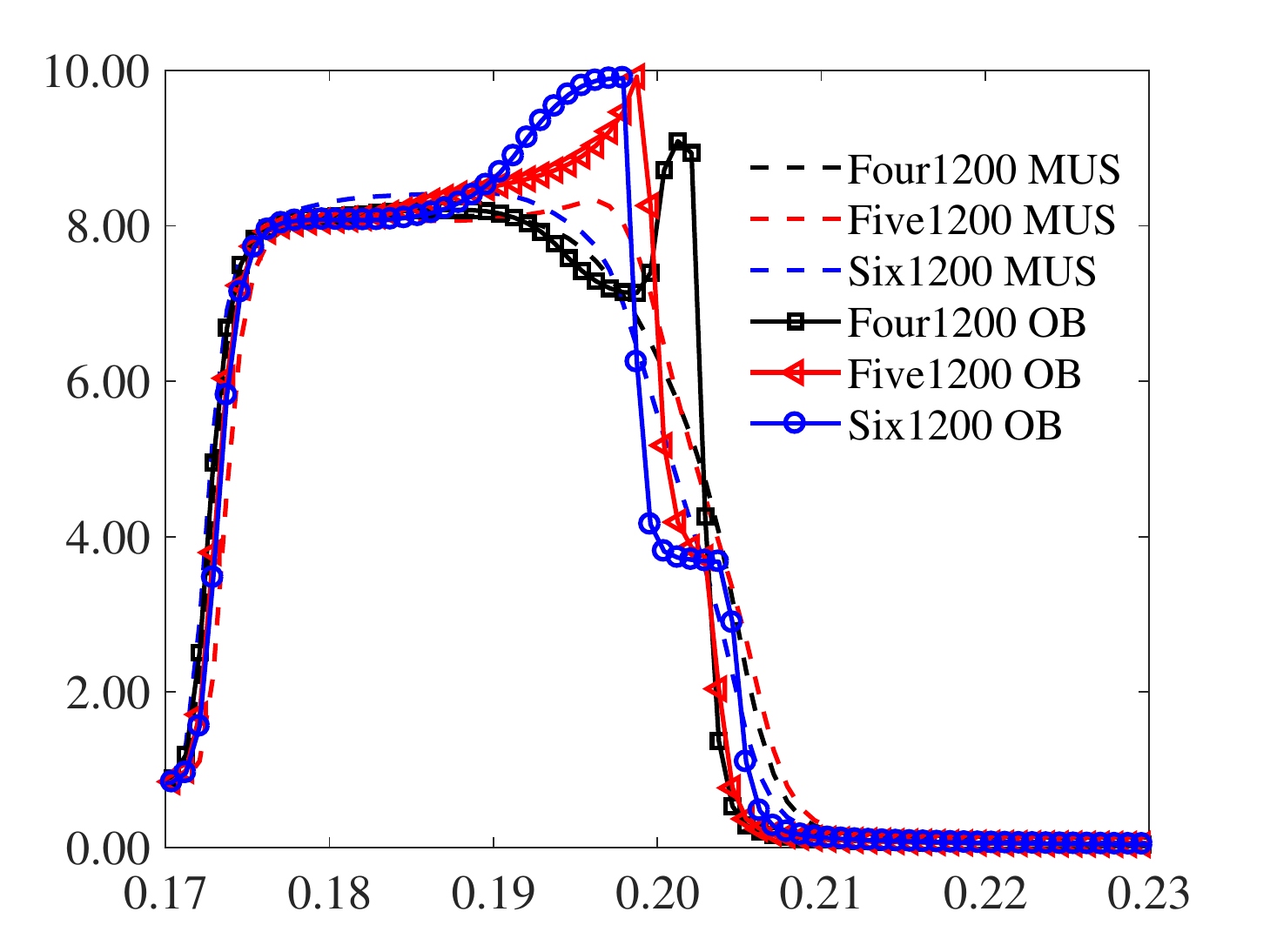}}
\subfloat[Volume/mass fraction, locally enlarged]{\label{fig:MUS_OB2}\includegraphics[width=0.5\textwidth]{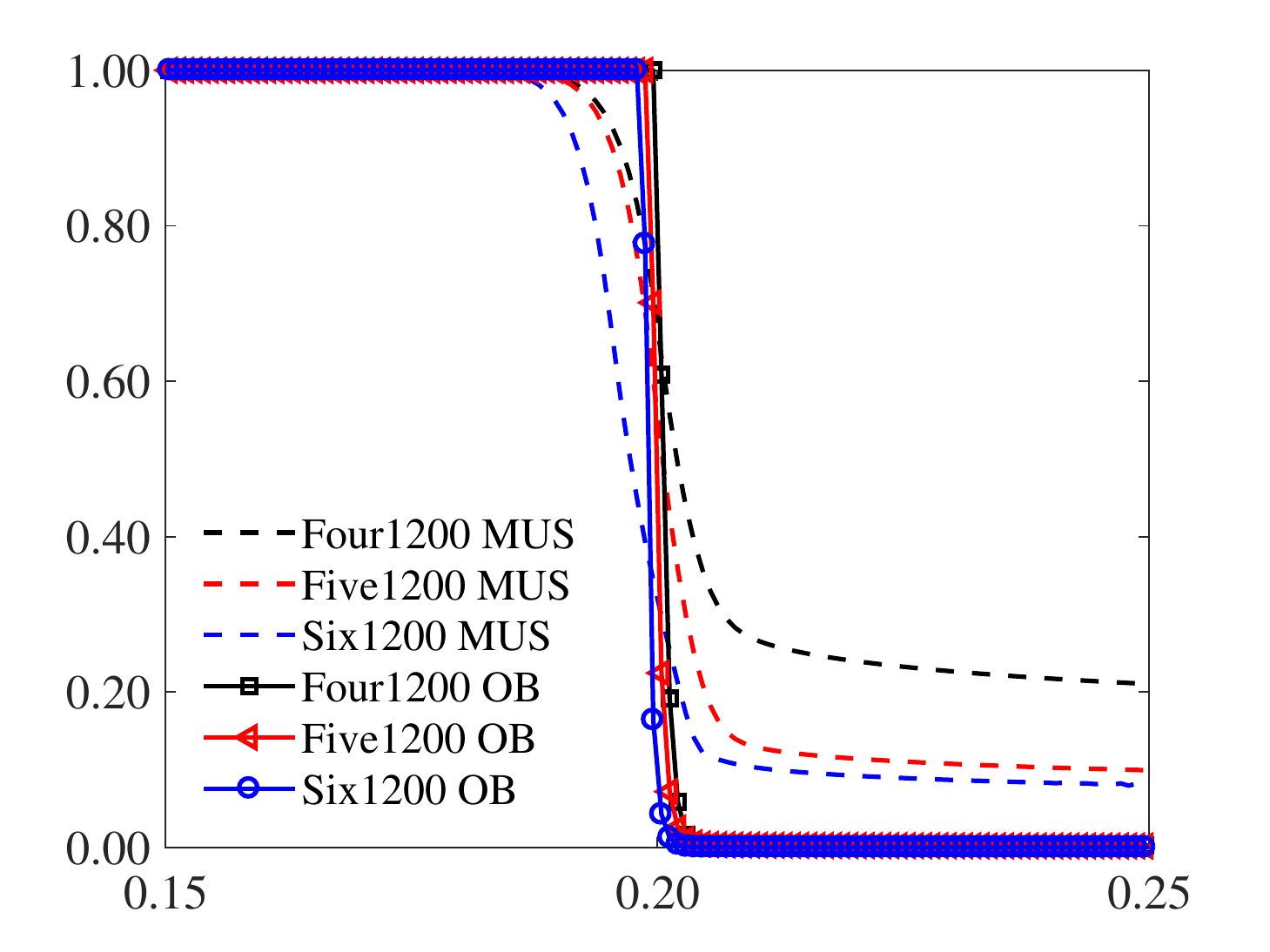}}
\caption{Comparison of the numerical results obtained the Overbee and MINMOD limiter schemes.}
\label{fig:MUS_OB} 
\end{figure}

To demonstrate the interface-sharpening effect, we compare the results obtained with the MINMOD limiter to those obtained with the Overbee scheme in \Cref{fig:MUS_OB}. One can see that with the Overbee limiter, the  diffused interface is within 2-3 computational cells, which is much less in comparison to the MINMOD scheme.

For the present problem, all materials are described with the ideal gas EOS. In this case, the effective temperature averaging procedure (\cref{eq:Taverage}) of the one-temperature five-equation model and  the temperature relaxation procedure (\cref{eq:ht_eq3_TR}) of the six-equation model yield the same result for temperature. 
The former neglects the effect of temperature relaxation on volume fraction within the diffused interface. This diffused interface is narrowed into 2-3 computational cells thanks to the interface-sharpening technique. Therefore, the advantage of the proposed model for this problem is not so evident as  that for the water-air shock tube problem in \cref{subsec:shocktubetest}.

\paragraph{Laser ablative Rayleigh–Taylor instability in a 2D thin target}
Next we consider the laser ablation problem in the 2D formulation. The interface is initially perturbed and has the following form:
\[
x_{interface} = x_R - A_m cos\left( 2 \pi y / L_y \right), 
\]
where $A_m$ is the perturbation amplitude taken as $A_m = 0.02L_y$.

\begin{figure}[htb]
\centering
{\label{fig:thinCHa}\includegraphics[width=0.5\textwidth]{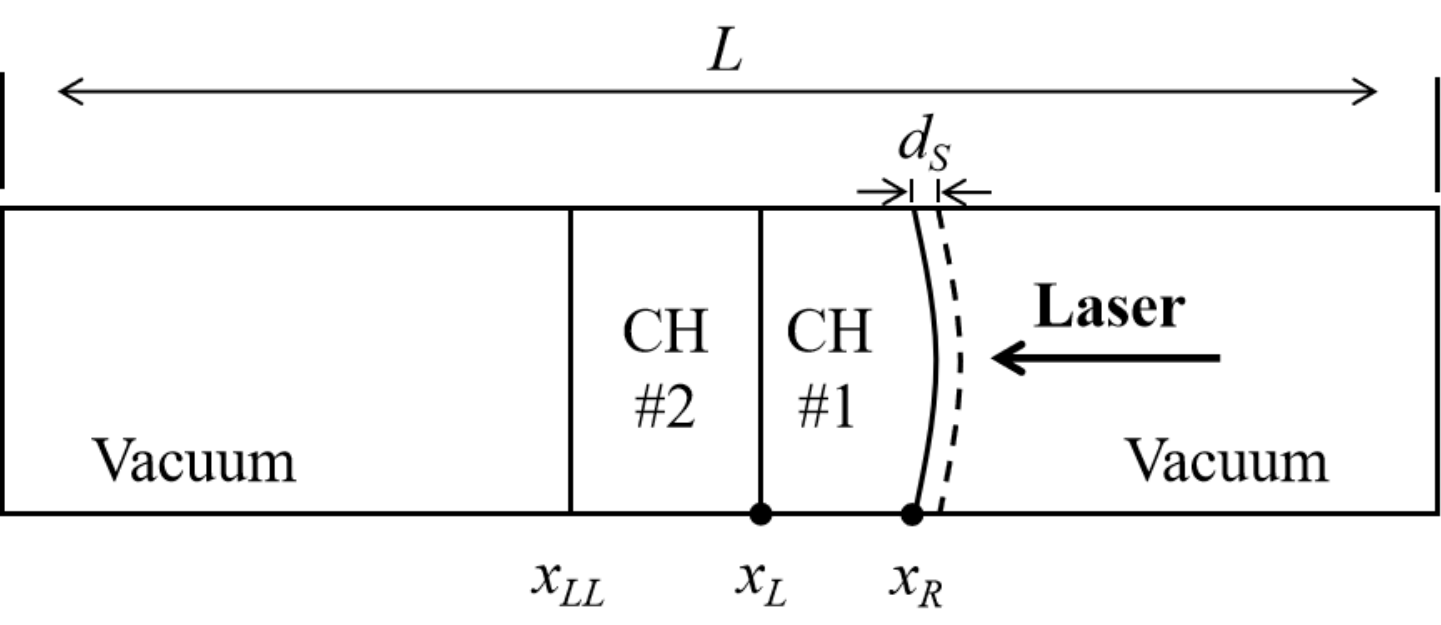}}
\caption{Schematic of the laser ablation of a thin two-phase target.}
\label{fig:thinCH} 
\end{figure}

\begin{figure}[htb]
\centering
{\label{fig:sixEvolve1}\includegraphics[width=0.75\textwidth]{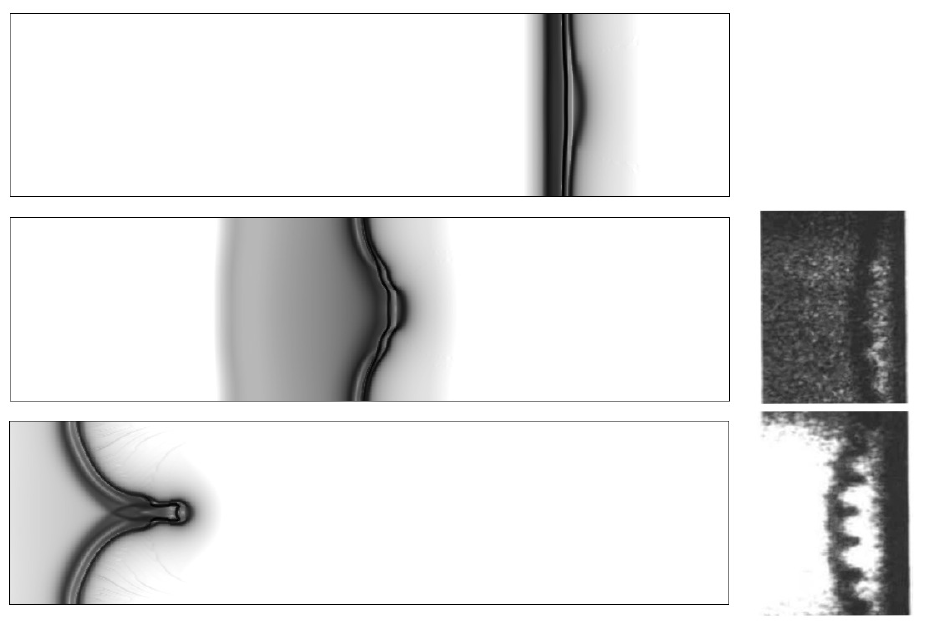}}
\caption{Evolution of the numerical Schlieren obtained with the proposed six-equation model. The figures on the right are experimental images from  \cite{Atzeni2004}.}
\label{fig:sixEvolve} 
\end{figure}

The laser ablation  of a thin target  is considered, which is accompanied with the development of Rayleigh–Taylor instability. The problem is a two-phase version of that in \cite{Li2020Numerical}. The problem set-up is displayed in \Cref{fig:thinCH}. 
The left and right ends of the target are located at $x_{LL} = 0.50 L$ and $x_{R} = 0.70L$, respectively. The two CH materials are separated by a planar interface at $x_L = 0.65L$. The evolution of the ablated target modelled with the proposed six-equation model is demonstrated in \Cref{fig:sixEvolve}. Here, the numerical Schlieren is qualitatively compared with the  experimental results for single material from GEKKO XII \cite{Atzeni2004}.

\begin{figure}[htb]
\centering
{\includegraphics[width=0.75\textwidth]{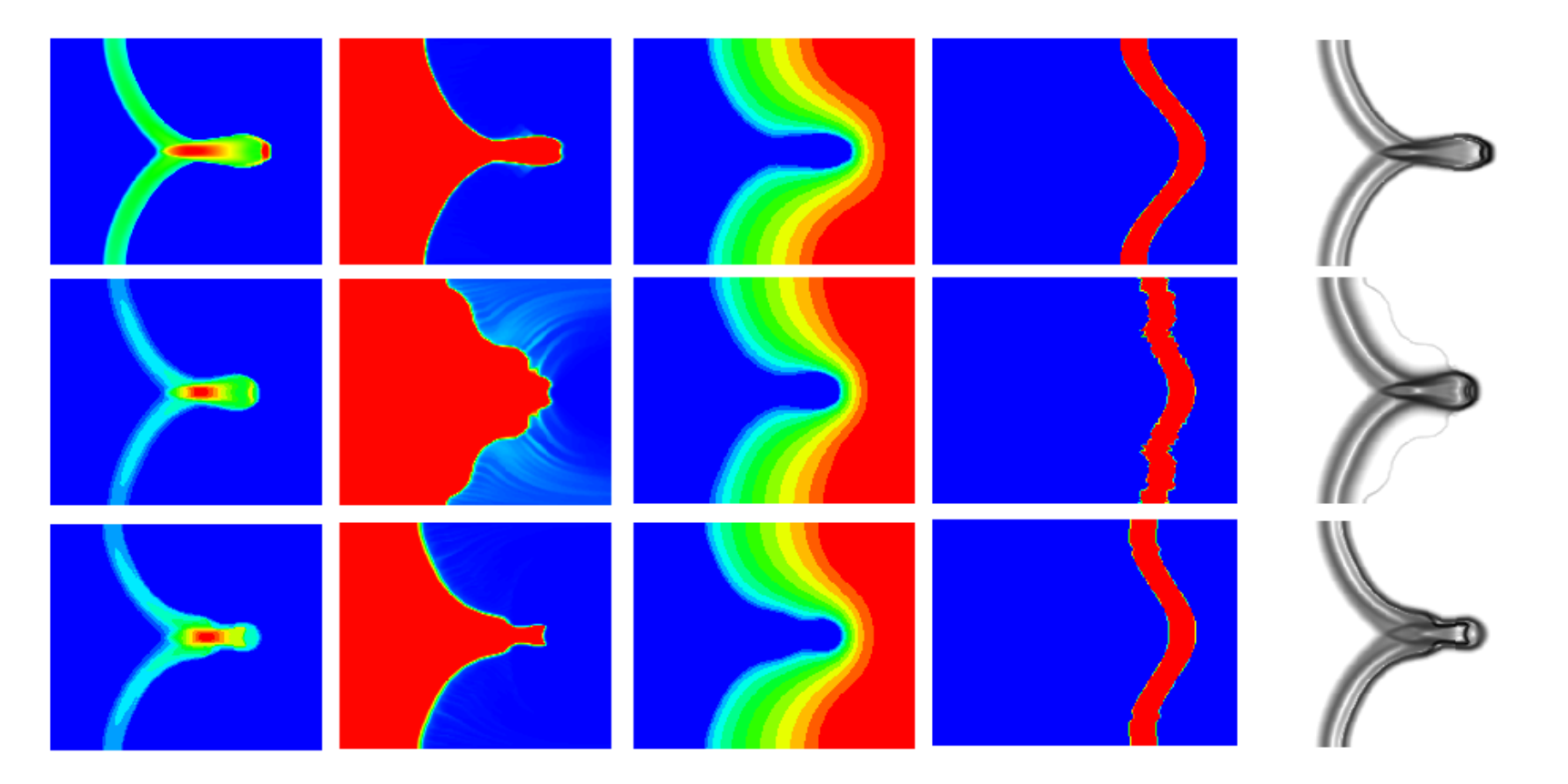}}
\caption{Comparison of the numerical results obtained with the four- (the first row), five- (the second row) and the six-equation model (the third row). Displayed results from left to right: density, volume fraction, temperature, laser absorption area, numerical Schlieren.}
\label{fig:thin85compareModels} 
\end{figure}

The numerical results obtained with different models are compared in \Cref{fig:thin85compareModels}. We can see that although the density distributions obtained with the three models are similar in appearance, the shapes of the material interfaces are different from each other. The material interface obtained with the five-equation model is more diffusive and quite different from the others. The difference in critical density distribution can be seen from the laser absorption area. Again the one-temperature five-equation model result is found to be much different from the other two, mostly due to the exceeded numerical diffusion of the material interface {\color{black}and violation of the second law of thermodynamics in the diffused zone.}

\section*{Conclusion}
In this paper we have established a temperature non-equilibrium model for modelling compressible two-phase flows with taking into account the dissipative thermal conduction and viscosity. We have proposed numerical methods based on the fractional step approach for solving the proposed model.  In this approach, the hyperbolic part  of the governing equations  is solved with the Godunov-HLLC scheme, and the parabolic part with the method of local iterations based on Chebyshev parameters. 

The proposed model have demonstrated the following advantages.
\begin{enumerate}
\item[•]It is thermodynamically consistent.
\item[•]It ensures temperature equilibrium during the heat conduction process by implementing a special phase thermal relaxation.
\item[•]It includes the effect of mechanical relaxation, thermal relaxation and heat conduction on the volume fraction.
\item[•]Numerically, it maintains the pressure, velocity and temperature equilibrium, thus avoids spurious oscillations in the vicinity of material interfaces.
\item[•]It shows superior convergence performance when compared to other models  with non-physical diffused mixture.
\item[•]Thanks to its physical consistency with the most complete Baer-Nunziato model, our model can be used for simulating two-phase flows with  both resolved and non-resolved interfaces.
\end{enumerate}

We have compared the proposed six-equation model  with the one temperature, one pressure five-equation model both analytically and numerically. Our analysis shows that this five-equation model is not consistent with the second law of thermodynamics. Numerical experiments on the laser ablation of a CH target demonstrate that the temperature-equilibrium five-equation model yields numerical results much different from those of the four-equation and six-equation models. 

In our future work we plan to include enthalpy diffusion into our model.

\section*{Acknowledgement}
The first author appreciate Professor Wenhua Ye, Professor Junfeng Wu and Dr. Shuai Wang for motivating discussions.

\bibliography{references}
\bibliographystyle{abbrvnat}


\end{document}